\documentclass[11pt,french]{amsart}
\usepackage[utf8]{inputenc}
\usepackage{latexsym}
\usepackage[all]{xy}
 
\usepackage{amsmath}
\usepackage{amssymb}
\usepackage{amsfonts}

\usepackage[useregional]{datetime2}

\usepackage{latexsym}

\usepackage{url}
\usepackage{longtable}
\usepackage{dsfont}

\usepackage{enumitem,xcolor}
\setlist{
topsep=2ex, partopsep=2ex, parsep=.75ex, itemsep=1ex, font=\normalfont, before={\itshape}
}
\newlist{myenumerate}{enumerate}{4}
\setlist[myenumerate]{label=(\Alph*), ref=\Alph*, before={\color{red}\scshape}}
\newlist{mydescription}{description}{1}
\setlist[mydescription]{font=\normalfont\normalcolor\sffamily, before*={\
scshape}}

\def\F{{\mathbb F}} 
\def\N{{\mathbb N}}
\def\Z{{\mathbb Z}}
\def\Q{{\mathbb Q}}
\def\R{{\mathbb R}}

{\newtheorem{theorem}{Theorem}[section]   
   \newtheorem{proposition}[theorem]{Proposition}
   \newtheorem{lemma}[theorem]{Lemma}  
   \newtheorem{remark}[theorem]{Remark} 
   \newtheorem{corollary}[theorem]{Corollary}

   \newtheorem*{definition}{Definition}    
{\theoremstyle{definition}   
   }
{\theoremstyle{remark}

   \newtheorem*{preuve}{Proof}

   \newtheorem*{preuvele}{Proof of the lemma}  

\author{Michel Vaqui\'e}
\address{Institut de Math\'ematiques de Toulouse UMR 5219, CNRS, Universit\'e de Toulouse,  
UPS, 118 route de Narbonne, F-31062 Toulouse Cedex 9, 
France} 
\email{vaquie@math.univ-toulouse.fr} 
\title[Cuts of an ordered abelian group]{Cuts of an ordered abelian group}

\begin{document}  


\begin{abstract}

In this article, we study the cuts of a totally ordered abelian group $\Gamma$. We begin by recalling some results on ordered sets $I$ and on the associated sets ${\bf IS}(I)$ and ${\bf FS}(I)$ of initial and final segments of $I$. 

For a totally ordered set $I$ we review the notion of an $I$-structure defined on a module over a ring $R$, and the definition of the Hahn product of a family of $R$-modules indexed by $I$. 

The set ${\bf Cv}(\Gamma )$of convex subgroups of a totally ordered group $\Gamma$ is also a totally ordered set, canonically isomorphic to the set of cuts of the subset ${\bf Pr}(\Gamma )$of principal convex subgroups. One of the first results is then to equip the group $\Gamma$ with an $I$-structure where $I$ is the set ${\bf Pr}(\Gamma )$ endowed with the opposite order.

We associate a convex subgroup with every cut of the group $\Gamma$, and conversely, we can associate a family of cuts with every convex subgroup of $\Gamma$. 
It is by looking at these subgroups, and the $I$-structure of $\Gamma$ that we can obtain a classification of the different types of cuts.
\end{abstract}  


\vskip .2cm 
\date{Novembre 2025} 
   
\maketitle   

\tableofcontents

		      \section*{Introduction}
%

In this article we propose to study the structure of totally ordered abelian groups, and more specifically the cuts, or initial and final segments of these groups. This study is motivated by the importance of the role of groups in the behavior of valuations.
We recall that a valuation $\nu$ on a ring $R$ is a map $\xymatrix{\nu : R \ar[r] & \Gamma \cup \{+\infty\}}$, satisfying 
$\nu (xy)=\nu (x)+\nu (y)$, $\nu (x+y) \geq Inf(\nu (x),\nu (y))$, and $\nu (x)=+\infty$ if and only if $x=0$, 
where $\Gamma$ is a totally ordered abelian group, and we can define the group of orders $\Gamma _{\nu}$ as the ordered subgroup generated by the image of $R^* = R \setminus \{ 0 \}$. 
Then some properties of the valuation $\nu$ are determined uniquely by properties of the group of values $\Gamma _{\nu}$, in particular the rank and rational rank of the valuation are defined as the rank and rational rank of the totally ordered group $\Gamma _{\nu}$.

Let $\nu$ be a valuation on a field $K$, and let $R_{\nu}$ be its valuation ring, $R_{\nu} =\{ x \in K \ | \ \nu (x) \geq 0 \}$, then the rank of the valuation is equal to the Krull dimension of the ring $R_{\nu}$ and there is a bijection between the set of ideals of the ring $R_{\nu}$ and the set of initial segments of the ordered group $\Gamma _{\nu}$. 

The notion of composite valuation is also completely determined by the structure of the group of values.
More precisely, with any convex subgroup $\Delta$ of $\Gamma$ we can associate a decomposition of the valuation $\nu$ in the following manner.
As $\Delta$ is a convex subgroup of $\Gamma _{\nu}$ the quotient group $\bar\Delta = \Gamma _{\nu} / \Delta$ is totally ordered and the composition of the valuation $\xymatrix{\nu : K^* \ar[r] & \Gamma _{\nu}}$ and the natural morphism $\xymatrix{\Gamma _{\nu} \ar[r] & \bar\Delta}$ is a valuation $\nu _1$ on $K$. Moreover the valuation $\nu$ induces a valuation $\bar\nu$ on the residue field $\kappa _{\nu _1}$ of the valuation $\nu _1$ with group of values equal to $\Delta$.  
We say that $\nu$ is a composite valuation, composite with the valuations $\nu _1$ and $\bar\nu$, and we write $\nu = \nu _1 \circ \bar\nu$. 
The ranks of the valuations $\nu _1$ and $\bar\nu$ is strictly lower than the rank of the valuation $\nu$, we can therefore hope that their study is simpler and that the decomposition of the valuation $\nu = \nu _1 \circ \bar\nu$ is a way to approach the study of it.

Let $(K,\nu )$ be a valued field and let $L$ be a finite cyclic extension of $K$, then to study the extensions $\mu$ of the valuation $\nu$ to $L$, we introduced the notion of an \emph{admissible family of valuations}.
The most difficult case is when the extension of valued fields $(L,\mu )/(K,\nu )$has a non-trivial defect, and in this case we showed that there appears in the admissible family associated with the extension a \emph{limit augmented valuation}. 
In constructing the admissible family, the existence of the limit augmented valuation is a consequence of the existence of initial segments without greatest element. It is therefore natural to understand the structure of these initial segments to approach the study of extension of valuations. 
We refer to the articles \cite{Va1, Va2, Va3} for the presentation of these results.

All these remarks show that the study of totally ordered groups, and more precisely the study of its initial segments, and their behavior when we consider the convex subgroups of the ordered group, are essential elements for studying valuations. 

Even though the valuations that appear most often in algebraic geometry are valuations of finite rank, we will consider the totally ordered groups in all generality, that is to say without any hypothesis of finiteness on their rank. 
To do this, we must also consider the set of convex subgroups of a totally ordered group; this set is completely ordered by inclusion, and it is the cuts of this set that allow us to better understand the structure of the group.

\vskip .5cm

 In the first part of the article we recall some elementary results on ordered sets. In particular, we study the initial and final segments as well as the cuts of a totally ordered set, we recall that an initial segment, resp. a final segment, of a totally ordered set $I$  is a part $\Sigma$ of $I$ such that any element smaller, resp. greater, than an element of $\Sigma$ also belongs to $\Sigma$, and a cut $\Lambda$ of $I$ is a partition $(\Lambda _- , \Lambda _+)$ of the set $I$ such that $\Lambda _-$ is an initial segment and $\Lambda _+$ is a final segment. 
 We recall also that an interval of $I$ is the intersection of an initial segment and a final segment. 

The sets ${\bf IS}(I)$, ${\bf FS}(I)$ and ${\bf Cp}(I)$ respectively of initial segments, final segments and cuts of a totally ordered set $I$ are also totally ordered sets, and there are natural isomorphisms of ordered sets between them. 

Any morphism $\xymatrix{u:I \ar[r] & J}$ of ordered sets induces morphisms between the corresponding sets of initial and final segments, and of cuts. 
To study the behavior of these sets by such a morphism, it is practical and natural to adopt a categorical point of view, i.e. to consider an ordered set as a category. We then see that the set ${\bf IS}(I)$ of initial segments is isomorphic to the category of presheaves on $I$. 
In this way the morphism $\xymatrix{u _!:{\bf IS}(I) \ar[r] & {\bf IS}(J)}$ induced by the morphism $\xymatrix{u :I \ar[r] & J}$ corresponds exactly to the functor between the categories of presheaves associated with the sets $I$ and $J$ seen as categories.

\vskip .2cm 

In the second part of the article we consider the notion of $I$-structure on an $R$-module, where $I$ is a totally ordered set. 
This is an increasing filtration of $R$-module $\Gamma$ by a family $\Delta _I (\Gamma ) = \bigl ( \delta (\Gamma ) _i ^- , \delta (\Gamma ) _i ^+ \bigr )$ of submodules indexed by the set $I \times \{ - , + \}$, with the order defined by $(i,-) < (i,+)$ for any $i \in I$ and $(j,+) < (i,-)$ for $i<j$ in $I$. 
With this $I$-structure we can associate the family of quotient modules $\bigl (\varepsilon (\Gamma ) _i\bigr )$, with $\varepsilon (\Gamma ) _i = \delta (\Gamma) _i ^+ / \delta (\Gamma) _i ^-$. 

We consider also the \emph{Hahn product} of a family $( \Theta _i ) _{i \in I}$ of $R$-modules indexed by a totally ordered set $I$. By definition the Hahn product $\prod _{i \in I} ^{(H)} \Theta _i$ is the submodule of the product $\prod _{i \in I} \Theta _i$ made up of the elements $\underline{x}=(x _i)$ whose support $Supp(\underline{x}) = \{ i \in I \ | \ x _i \not= 0 \}$ is a well-ordered subset of $I$. 
We can endow the Hahn product with a natural $I$-structure $\Delta _I (\Theta ) = \bigl ( \delta (\Theta ) _i ^- , \delta (\Theta ) _i ^+ \bigr )$ defined by  
$\delta (\Theta )_i^- = \prod _{j > i} ^{(H)} \Theta _j$ and $\delta (\Theta )_i^+ = \prod _{j \geq i} ^{(H)} \Theta _j$. 

Then the main result of this part is the fact that, when $R$ is a field, for any $R$-module $\Gamma$ with an $I$-structure there exists a natural morphism $\xymatrix{\underline{v}:\Gamma \ar[r] & \varepsilon (\Gamma ) = \prod _{i \in I} ^{(H)} \varepsilon (\Gamma )_i}$, where $\prod _{i \in I} ^{(H)} \varepsilon (\Gamma )_i$ is the Hahn product. Moreover this morphism in an injective \emph{immediate} morphism, which means that it is very close to being an isomorphism.

\vskip .2cm 

In the third part of the article we study the structure of totally ordered abelian groups. 
Firstly we are interested in the set of convex subgroups of a totally ordered group $\Gamma$, we recall that a convex subgroup is a subgroup that is also an interval of $\Gamma$. 
Then the set of convex subgroups ${\bf Cv}(\Gamma )$ of $\Gamma$ is totally ordered by the inclusion, and one of the most useful approaches is to look at the set of cuts of ${\bf Cv}(\Gamma )$. 

We define a principal convex subgroup of $\Gamma$ as a subgroup $\Delta ^+ _{\zeta}$ defined as the maximal convex subgroup of $\Gamma$ which contains a fixed element $\zeta$, more precisely we have the equivalence 
$$\Delta \subset \Delta ^+ _{\zeta} \ \Longleftrightarrow \ \zeta \in \Delta \ .$$ 
Then we can show that the subset ${\bf Pr}(\Gamma )$ is the set of immediate successors of the ordered set ${\bf Cv}(\Gamma )$, and by general results from the first part we see that there is a natural isomorphism of ordered sets between ${\bf Cv}(\Gamma )$ and the set of cuts of ${\bf Pr}(\Gamma )$. 
Moreover the ordered group $\Gamma$ is endowed with an ${\bf I}_{\Gamma }$-structure as defined above, where ${\bf I}_{\Gamma }$ is the set ${\bf Pr}(\Gamma )$ with the opposite total order. 

In this section we also consider the totally ordered groups $\Gamma$ obtained as Hahn product of a family $\bigl (\Gamma _i \bigr )_{i \in I}$ of totally ordered groups, where $I$ is a totally ordered set. We then describe the convex subgroups and the principal convex subgroups of such a group $\Gamma = \prod _{i \in I} ^{(H)} \Gamma _i$, for this it is necessary to consider both the convex subgroups and the principal convex subgroups of each of the factors $\Gamma _i$, and the ordered set structure of the set of indices $I$.  
We then deduce from the results of the second section of the article that any totally ordered divisible group $\Gamma$ is embeddable to a Hahn product of the form $\prod _{i \in I} ^{(H)} \Gamma _i$, where each group $\Gamma _i$ is an ordered group of rank one. 

\vskip .2cm 

The fourth part of the article is devoted to cuts, or initial segments, of totally ordered groups. 
With any non trivial cut $\Lambda =(\Lambda _-,\Lambda _+)$ of a totally ordered group $\Gamma$, or with any non trivial initial segment $\Lambda _-$ of $\Gamma$, we associate a convex subgroup $\Delta (\Lambda )$ of $\Gamma$, called the \emph{invariance subgroup} of $\Lambda$, or of $\Lambda _-$, defined as the greatest convex subgroup $\Theta$ of $\Gamma$ verifying the equalities $\Theta + \Lambda _- = \Lambda _-$.

The subgroup $\Delta (\Lambda )$ is also characterized as the largest convex subgroup $\Theta$ such that the image of the cut $\Lambda$ by the canonical morphism $\xymatrix{w _{\Theta}: \Gamma \ar[r] & \Gamma / \Theta}$ is still a non-trivial cut of the totally ordered group $\Gamma / \Theta$. 
On the other hand the cut $\Lambda$ defines a non-trivial cut of any convex subgroup $\Theta$ of $\Gamma$ satisfying $\Delta (\Gamma ) \subsetneq \Theta$. 

We can also associate with any convex subgroup $\Delta$ and with any element $\zeta$ of $\Gamma$ two cuts $\Lambda ^{\scriptscriptstyle \leq \zeta + \Delta } $ and $\Lambda ^{\scriptscriptstyle \geq \zeta + \Delta}$ of the group $\Gamma$. 
These two cuts are called \emph{relatively principal}, they are the simplest cuts to describe and the invariance subgroup of each of these two cuts is equal to $\Delta$. 
If we have the strict inequalities $\Lambda ^{\scriptscriptstyle \leq \zeta + \Delta } < \Lambda < \Lambda ^{\scriptscriptstyle \geq \zeta + \Delta}$, we have to consider the maximum symmetric intervals included in the initial segment $\Lambda _-$, or in the final segment $\Lambda _+$ of the cut, and we can arrive at a complete classification of the cuts in three different types: \emph{relatively principal}, \emph{gapped} and \emph{tightened}. 
Then if $\Lambda$ is a non trivial cut of $\Gamma$, for any convex subgroups $\Theta _1$ and $\Theta _2$ with $\Theta _1 \subset \Delta (\Lambda ) \subsetneq \Theta _2$, the non-trivial cut of the totally ordered group $\Theta _2 / \Theta _1$ induced by $\Lambda$ is of the same type as the cut $\Lambda$. 

\vskip .2cm 

This article is essentially a survey; most of the results are well known, and we refer, for instance, to F.-V. Kuhlmann's articles for a more detailed bibliography (\cite{Ku1, Ku-N, Ku-Ku}).

\vskip .2cm

		      \section{Reminders on totally ordered sets}
%

\subsection{Cuts of a totally ordered set} 

Let $I$ be a totally ordered set, where we note $``\leq"$ the order and we call it \emph{smaller or equal}, and we use the usual notations $``\geq"$, $``<"$, and $``>"$, and the corresponding terminology. 
We recall the following notations and definitions. 

For any subset $X \subset I$, we note $``{\rm sup}\, X"$ and $``{\rm inf}\, X"$ respectively the \emph{supremum} (or the least upper bound) and the \emph{infimum} (or the greatest lower bound) of $X$, and $``{\rm max}\, X"$ and $``{\rm min}\, X"$ respectively the \emph{maximum} (or the greatest element) and the \emph{minimum} (or the smallest or least element) of $X$, if they exist.  

A subset $\Sigma$ of $I$ is called an \emph{initial segment} (resp. a \emph{final segment}) if for any $j$ in $\Sigma$, any element $i$ in $I$ such that $i<j$ (resp. $i>j$) belongs to $\Sigma$. 
  
We say that $\Xi$ is an \emph{interval} if for any  $j$ and $j '$ in $\Xi$ with $j < j '$, any element $i$ of $I$ such that $j < i < j '$ belongs to $\Xi$.
An interval is the intersection of an initial segment and a final segment. 

For any $j$ in $I$ the subsets $\{i \in I \ / \ i < j \}$ and $\{i \in I \ / \ i \leq j \}$, denoted respectively $I_{<j}$ (or $] - \infty , j [$) and $I_{\leq j}$ (or $] - \infty , j ]$) are initial segments, but in general there exist other ones. 
The initial segment $\Sigma$ admits a greatest element if it is of the form $I_{\leq j}$. 
In a similar way we may define the final segments $I_{>j'}$ (or $] j ', + \infty [$) and $I_{\geq j'}$ (or $[ j ', + \infty [$), the smallest element of a final segment $I_{\geq j}$, and the intervals $] j ', j [$, $] j ', j ]$, $[ j ', j [$ and $[ j ', j ]$.

An initial segment (resp. a final segment) $\Sigma$ is \emph{trivial} if $\Sigma = \emptyset$ or if $\Sigma = I$. 

\vskip .2cm 

A \emph{cut} $\Lambda$ of a totally ordered set $I$ is a partition $(\Lambda _-,\Lambda _+)$ of $I$ such that for any $i \in \Lambda _-$ and any $j \in \Lambda _+$ we have $i < j$. It is equivalent to giving a partition $I = \Lambda _- \sqcup \Lambda _+$ of $I$ where $\Lambda _-$ is an initial segment and $\Lambda _+$ is a final segment. 
We say that the cut $\Lambda$ is \emph{non trivial}, or is a \emph{Dedekind cut}, if neither of the two subsets $\Lambda _-$ or $\Lambda _+$ is the empty set. 

We remind that there exist bijections between the set ${\bf IS}(I )$ of initial segments, the set ${\bf FS}(I )$ of final segments, and the set ${\bf Cp}(I )$ of cuts of $I$ which consist of associating with every initial segment $\Sigma$ its complement $\Sigma ^{C}= I \setminus \Sigma$, and the cut $\Lambda = \bigl ( \Lambda _- , \Lambda _+ \bigr )$ of $I$ defined by $\Lambda _- = \Sigma$ and $\Lambda _+ = \Sigma ^{C}$.   

This induces bijections between the set ${\bf IS}(I )^{\circ}$ of non trivial initial segments, the set ${\bf FS}(I )^{\circ}$ of non trivial final segments, and the set ${\bf Cp}(I )^{\circ}$ of non trivial cuts of $I$. 

\begin{remark}\label{rmq:coupure}
(1) For any initial segment $\Sigma$ we have the following equivalences 
$$i \in \Sigma \ \Longleftrightarrow \  I_{\leq i} \subset \Sigma \ \Longleftrightarrow \ I_{<i} \subsetneq \Sigma \ .$$ 

(2) For any cut $\Lambda = (\Lambda _- , \Lambda _+)$ of $I$ we have the equalities  
$$\Lambda _- = \bigcup _{j \in \Lambda _-} I_{\leq j} = \bigcap _{j \in \Lambda _+} I_{< j} 
\quad\hbox{and}\quad
\Lambda _+ = \bigcup _{j \in \Lambda _+} I_{\geq j}  = \bigcap _{j \in \Lambda _-} I_{> j}  \ .$$ 
\end{remark} 

\vskip .2cm

\begin{remark}\label{rmq:inclusion-intervals} 
\begin{enumerate} 
\item Let $(\Sigma _{\alpha})_{\alpha \in A}$ be a family of initial segments (resp. final segments) of $I$, then the sets $\bigcap _{\alpha \in A} \Sigma _{\alpha}$ and $\bigcup _{\alpha \in A} \Sigma _{\alpha}$ are initial segments (resp. final segments) of $I$. 

\item Let $(\Xi _{\alpha})_{\alpha \in A}$ be a family of intervals of $I$, then the set $\bigcap _{\alpha \in A} \Xi _{\alpha}$ is still an interval of $I$. 
Moreover if the family $(\Xi _i)_{\alpha \in A}$ is totally ordered by inclusion (i.e. for any two intervals $\Xi _{\alpha}$ and $\Xi _{\beta}$ of the family we have $\Xi _{\alpha} \subset \Xi _{\beta}$ or $\Xi _{\beta} \subset \Xi _{\alpha}$), the set $\bigcup _{\alpha \in A} \Xi _{\alpha}$ is still an interval of $I$. 
\end{enumerate}
\end{remark}

\vskip .2cm

We can define an order $\leq$ on the set ${\bf Cp}(I )$ of cuts of $I$ as follows: let $\Lambda = \bigl (\Lambda _- , \Lambda _+ \bigr )$ and $\Lambda '= \bigl (\Lambda '_- , \Lambda '_+ \bigr )$ be two cuts of $I$, we claim  
$$\Lambda \leq \Lambda '  \  \Longleftrightarrow  \  \Lambda _- \subset \Lambda ' _-  \  \Longleftrightarrow  \  \Lambda '_+ \subset \Lambda _+ \ .$$ 

This order corresponds to the order defined on the set ${\bf IS}(I )$ (resp. the set ${\bf FS}(I )$) in the following way 
$$\Sigma \leq \Sigma '  \  \Longleftrightarrow  \  \Sigma \subset \Sigma '  \quad ( {\rm resp.}\ \Sigma \leq \Sigma '   \  \Longleftrightarrow  \  \Sigma ' \subset \Sigma ) \ .$$ 

\vskip .2cm 

For any $i$ belonging to $I$ we can define two cuts $\Lambda ^{\leq i} = \Lambda ^{> i}$ and $\Lambda ^{\geq i} = \Lambda ^{< i}$ of the set $I$ associated respectively with the final segments ${(\Lambda ^{\leq i})}_+ = I_{> i}$  and ${(\Lambda ^{\geq i})}_+ = I_{\geq i}$, or the initial segments ${(\Lambda ^{\leq i})}_- = I_{\leq i}$  and ${(\Lambda ^{\geq i})}_- = I_{<i}$.  
For $i < j$ we have the inequalities: 
$$ \Lambda ^{\geq i} \leq \Lambda ^{\leq i} \leq \Lambda ^{\geq j} \leq \Lambda ^{\leq j} \ ,$$ 
and for any cut $\Lambda$ of $I$ we have the following equivalences:  
$$i \in \Lambda _- \ \Longleftrightarrow \  \Lambda  ^{\leq i} \leq \Lambda  \quad \hbox{and} \quad  i \in \Lambda _+ \ \Longleftrightarrow \  \Lambda \leq \Lambda ^{\geq i} \ .$$ 

\vskip .2cm 

We can define two embeddings of totally ordered sets  
$\xymatrix @C=5mm {\sigma ^{\leq}: I \ar@<-1.5pt>@{^{(}->}[r] & {\bf Cp}(I)}$ and $\xymatrix @C=5mm {\sigma ^{\geq}: I \ar@<-1.5pt>@{^{(}->}[r] & {\bf Cp}(I)}$ respectively by $\sigma ^{\leq}(i) = \Lambda ^{\leq i}$ and $\sigma ^{\geq}(i) = \Lambda ^{\geq i}$. 
These embedding correspond to the following morphisms of ordered sets: 
$$\xymatrix @R=0mm @C=5mm{**[l]\sigma ^{\leq}: I \ar@<-1.5pt>@{^{(}->}[rr] && **[r]{\bf IS}(I) & \hbox{and} &**[l]\sigma ^{\geq}: I \ar@<-1.5pt>@{^{(}->}[rr] && **[r]{\bf FS}(I) \\
**[l] i \ar@{|->}[rr] && I_{\leq i} && **[l] i \ar@{|->}[rr] && **[r]I_{\geq i} & \hskip -12mm  .} $$ 

Moreover if $I$ has no greatest element the image of the morphism $\sigma ^{\leq}$ is included in the set ${\bf IS}(I)^{\circ}$, and if $I$ has no smallest element the image of the morphism $\sigma ^{\geq}$ is included in the set ${\bf FS}(I)^{\circ}$.

\vskip .2cm 

We recall the following obvious result. 

\begin{remark}\label{rmq:successeur} 
Let $x$ and $y$ be two elements of a totally ordered set $I$, then the following assertions are equivalent: 

(a) $y < x$ and $I_{> y} \cap I_{< x} = \emptyset$ ($\Longleftrightarrow I_{<x} \subset I_{\leq y} \Longleftrightarrow I_{>y} \subset I_{\geq x}$).  

(b) $x = {\rm min}\, I_{> y}$. 

(c) $y = {\rm max}\, I_{< x}$.

\end{remark}

\vskip .2cm 

\begin{definition} 
Let $x$ and $y$ be two elements of $I$ which satisfy the equivalent properties of remark \ref{rmq:successeur}, then we say that $x$ is an \emph{immediate successor} of $y$, or that $y$ is an \emph{immediate predecessor} of $x$. 
We note $x = {\rm suc} (y)$ and $y = {\rm pre}(x)$.  

The subset of $I$ of all the elements $x$ which are immediate successors, or equivalently which admit an immediate predecessor, is denoted ${\rm Suc} (I)$, 
and the subset of $I$ of all elements $y$ which are immediate predecessor, or equivalently which admit an immediate successor, is denoted ${\rm Pre} (I)$. 
\end{definition}

\vskip .2cm 

\begin{lemma}\label{le:predecesseur-successeur-dans-CpI}
Let $\Lambda$ and $\Lambda '$ be two cuts of an ordered set $I$ with $\Lambda \leq \Lambda '$, then $\Lambda$ is the immediate predecessor of $\Lambda '$, or $\Lambda '$ is the immediate successor of $\Lambda$, if and only if there exists $i \in I$ such that $\Lambda = \Lambda ^{\geq i}$ and $\Lambda ' = \Lambda ^{\leq i}$. 

Then the morphisms $\sigma ^{\geq}$ and $\sigma ^{\leq}$ defined above induce isomorphisms of ordered sets: 
$$\xymatrix @R=0mm @C=10mm{**[l]\sigma ^{\geq} : I \ar[r] ^-{\sim} & **[r]{\rm Pre}({\bf Cp}(I)) \\ i \ar@{|->}[r] & **[r]\Lambda ^{\geq i} = \Lambda ^{< i}} \quad\hbox{and}\quad \xymatrix @R=0mm @C=10mm{**[l] \sigma ^{\leq } : I \ar[r] ^-{\sim} & **[r]{\rm Suc}({\bf Cp}(I)) \\ i \ar@{|->}[r] &**[r] \Lambda ^{\leq i} = \Lambda ^{> i}} \ .$$
\end{lemma}

\begin{preuve} 
We recall that we have the inequality $\Lambda \leq \Lambda '$ if and only if $\Lambda  _- \cap \Lambda '_+ = \emptyset$, and the inequality is strict if and only if $\Lambda  _+ \cap \Lambda '_- \not= \emptyset$. 
It is easy to see that $\Lambda = {\rm pre}(\Lambda ')$ if and only if there exists $i \in I$ such that we have the equality $\Lambda  _+ \cap \Lambda '_- = \{ i \}$. 

\hfill$\Box$ 
\end{preuve}

\vskip .2cm 

We recall that a totally ordered set $(I,\leq )$ is said to be \emph{well-ordered} if every nonempty subset $J$ of $I$ has a least element.

\begin{remark}\label{rmq:bien-ordonne} 
If $I$ is a well-ordered set any non trivial cut $\Lambda$ de $I$ is of the form $\Lambda = \Lambda ^{\geq i}$ for some element $i$ in $I$. 

Let $I^*=I$ if $I$ has no greatest element, and $I^*=I \setminus\{\bar\iota\}$ if $\bar\iota = {\rm max}\, I$. 
Then any element $x$ of $I^*$ admit an immediate successor, and we have a morphism of ordered sets $\xymatrix{{\rm suc}: I^* \ar[r] & I}$ that satifies $\sigma ^{\leq } = \sigma ^{\geq } \circ {\rm suc}$. 
\end{remark} 

If we choose  the embedding $\xymatrix{\sigma ^{\geq}: I \ar@<-1.5pt>@{^{(}->}[r] & {\bf Cp}(I)}$ and if we identify any $i \in I$ with the cut $\Lambda ^{\geq i} \in {\bf Cp}(I)$ (or with the final segment $I_{\geq i}$), we have the following equivalences: 
$$\begin{array}{rcl}
\Lambda \leq i & \Longleftrightarrow & i \in \Lambda _+ \ ; \\ 
\Lambda < i & \Longleftrightarrow & \Lambda _+ \cap I_{< i} \not= \emptyset \ ; \\ 
i \leq \Lambda &  \Longleftrightarrow &  I_{<i} \subset \Lambda _- \ \Longleftrightarrow \ \Lambda _- \subset I_{\leq i}  \ ; \\ 
i < \Lambda & \Longleftrightarrow & i \in \Lambda _-  \ .
\end{array}$$

In the same way if we choose the embedding $\xymatrix{\sigma ^{\leq}: I \ar@<-1.5pt>@{^{(}->}[r] & {\bf Cp}(I)}$ and identify any $i \in I$ with the cut $\Lambda ^{> i}  = \Lambda ^{\leq i} \in {\bf Cp}(I)$ (or with the initial segment $I_{\leq i}$), we have the following equivalences: 
$$\begin{array}{rcl}
i \leq \Lambda & \Longleftrightarrow & i \in \Lambda _- \ ; \\ 
i < \Lambda & \Longleftrightarrow & \Lambda _- \cap I_{> i} \not= \emptyset \ ; \\ 
\Lambda \leq i & \Longleftrightarrow & I_{> i} \subset \Lambda _+ \ \Longleftrightarrow \ \Lambda _+ \subset I_{\geq i}   \ ; \\  
\Lambda < i & \Longleftrightarrow & i \in \Lambda _+ \ .  
\end{array}$$

\vskip .2cm

For any subset $Y$ of the ordered set $I$, we can define two cuts $\Lambda ^{\scriptscriptstyle \leq Y} = \bigl ( \Lambda ^{\scriptscriptstyle \leq Y}_- , \Lambda ^{\scriptscriptstyle \leq Y} _+ \bigr )$ and $\Lambda ^{\scriptscriptstyle \geq Y} = \bigl ( \Lambda ^{\scriptscriptstyle \geq Y}_- , \Lambda ^{\scriptscriptstyle \geq Y} _+ \bigr )$ in the following way:
$$\begin{array}{lcl}
 \Lambda ^{\scriptscriptstyle \leq Y}_- = \{ j \in I \ | \ \exists \, i \in Y \ \hbox{such that} \ j \leq i \} & \hbox{and} & \Lambda ^{\scriptscriptstyle \leq Y} _+ = \{ j \in I \ | \ j>Y \} \\  
\Lambda ^{\scriptscriptstyle \geq Y} _- = \{ j \in I \ | \ j<Y \} & \hbox{and} &  \Lambda ^{\scriptscriptstyle \geq Y}_+ = \{ j \in I \ | \ \exists \, i \in Y \ \hbox{such that} \ j \geq i \} \ .
\end{array}$$ 
These cuts are also noted as follows: $\Lambda ^{\leq Y} = Y^+$ and $\Lambda ^{\geq Y} = Y^-$.  
For $Y=\{ i \}$ we recover the cuts $\Lambda ^{\leq i}$ and $\Lambda ^{\geq i}$ defined above, these cuts are sometime called \emph{principal}. 
We may notice that the cut $\Lambda ^{\scriptscriptstyle \leq Y}$, resp. $\Lambda ^{\scriptscriptstyle \geq Y}$, is principal if and only if the subset $Y$ has a greatest element, resp. a smallest element. 

A cut $\Lambda$ such that the initial segment $\Lambda _-$ has a greatest element and the final segment $\Lambda _+$ has a smallest element is called a \emph{jump} or a \emph{gap} (cf. {Ku1}).   

\begin{remark}\label{rmq:jump}
A cut $\Lambda =(\Lambda _- , \Lambda _+ )$ is a jump if and only if we have $\Lambda _- = I_{\leq i}$ and $\Lambda _+ = I_{\geq j}$ where $j$ is the immediate successor of $i$, and conversely any element $i$ of $I$ which admits an immediate successor $j$ defines a jump $\Lambda$ by $\Lambda = ( I_{\leq i} , I_{\geq j})$. 
\end{remark} 

\vskip .2cm 

If the set $Y$ is the empty set we get the trivial cuts $\Lambda ^{\scriptscriptstyle \leq \emptyset} =(\emptyset , I)$ and $\Lambda ^{\scriptscriptstyle \geq \emptyset} =(I,\emptyset  )$, otherwise we get the inclusion $\Lambda ^{\scriptscriptstyle \geq Y}_- \subset \Lambda ^{\scriptscriptstyle \leq Y}_-$, hence the relation $\Lambda ^{\scriptscriptstyle \geq Y} \leq \Lambda ^{\scriptscriptstyle \leq Y}$. 
Moreover the set $Y$ is included in $\Lambda ^{\scriptscriptstyle \geq Y}_+ \cap \Lambda ^{\scriptscriptstyle \leq Y}_- = \{ j \in I \ | \ \exists \, i , i' \in Y \ \hbox{such that} \ i \leq j \leq i' \}$, and we have the equality $Y = \Lambda ^{\scriptscriptstyle \geq Y}_+ \cap \Lambda ^{\scriptscriptstyle \leq Y}_-$ if and only if $Y$ is an interval of $I$. 

If we take $Y$ equal to the whole set $I$ we get the trivial cuts $\Lambda ^{\scriptscriptstyle \leq I} =(I, \emptyset )$ and $\Lambda ^{\scriptscriptstyle \geq I} =(\emptyset ,I )$, more generaly the cut $\Lambda ^{\scriptscriptstyle \leq Y}$ (resp. $\Lambda ^{\scriptscriptstyle \geq Y}$) is trivial if the set $Y$ is non bounded above (resp. non bounded below). 

\vskip .2cm 

In the same way for any subset $Y$ of $I$ we can define the two cuts $\Lambda ^{\scriptscriptstyle < Y}$ and $\Lambda ^{\scriptscriptstyle > Y}$ by the following 
$$\begin{array}{lcl}
 \Lambda ^{\scriptscriptstyle < Y}_- = \{ j \in I \ | \ \exists \, i \in Y \ \hbox{such that} \ j < i \} & \hbox{and} & \Lambda ^{\scriptscriptstyle < Y} _+ = \{ j \in I \ | \ j \geq Y \} \\  
\Lambda ^{\scriptscriptstyle > Y} _- = \{ j \in I \ | \ j \leq Y \} & \hbox{and} &  \Lambda ^{\scriptscriptstyle > Y}_+ = \{ j \in I \ | \ \exists \, i \in Y \ \hbox{such that} \ j > i \} \ .
\end{array}$$

\begin{remark}\label{rmq:<=leq->=geq}
If the set $Y$ has no greatest element (resp. no smallest element) the cuts $\Lambda ^{\scriptscriptstyle \geq Y}$ and $\Lambda ^{\scriptscriptstyle \geq Y}$ (resp. $\Lambda ^{\scriptscriptstyle \leq Y}$ and $\Lambda ^{\scriptscriptstyle < Y}$) are equal. 
And if the set $Y$ admits $y$ as greatest element (resp. smallest element), we have the equalities $\Lambda ^{\scriptscriptstyle \geq Y} = \Lambda ^{> y}$ and $\Lambda ^{\scriptscriptstyle \geq Y} = \Lambda ^{\geq y}$ (resp. $\Lambda ^{\scriptscriptstyle < Y} = \Lambda ^{< y}$ and $\Lambda ^{\scriptscriptstyle \leq Y} = \Lambda ^{\leq y}$). 
\end{remark} 

\vskip .2cm

\subsection{Categorical point of view} 

We can associate with any ordered set $I$ a category, also denoted $ I$, whose set of objects is the set $I$ and whose morphisms are defined in the following manner: 
for $i$ and $j$ in $I$ we have $Mor_I(i,j) = \emptyset$ if $j<i$ and $Mor_I(i,j) = \{ \star \}$ if $i \leq j$. 

\begin{proposition}\label{prop:Yoneda} 
For any totally ordered set $I$ the category associated with the set ${\bf IS}(I)$ of initial segments of $I$ is equivalent to the category of \emph{presheaves} on $I$ and the morphism 
$\xymatrix @R=-1mm @C=6mm{\sigma ^{\leq}: I \ar@<-1.5pt>@{^{(}->}[r] & {\bf IS}(I)}$ 
corresponds to the Yoneda embedding $\underline h_I$. 
\end{proposition} 

\begin{preuve} 
The initial segment $\Sigma \in {\bf IS}(I)$ corresponds to the presheaf defined by $\Sigma (j) =1$ if $j \in \Sigma$ and $\Sigma (j) =0$ if $j \notin \Sigma$.  
For any $i$ and $j$ in $I$ we have  $Mor_I(j,i) = \{ \star \}$ if and only if $j \in I_{\leq i}$, then we have $\underline h_I (i) = I_{\leq i} = \sigma ^{\leq}(i)$.  

\hfill$\Box$ 
\end{preuve}  
 
\vskip .2cm 

Let $\xymatrix{u:J \ar[r] & I}$ be a morphism of ordered sets, for any $j$ and $j '$ in $J$ we have the following relations
$$\begin{array}{rcl}
j \leq j ' & \Longrightarrow & u(j ) \leq u(j ') \\ 
u(j ) < u(j ') & \Longrightarrow & j  < j ' \ . 
\end{array}$$
The morphism $u$ corresponds to a functor between the associated categories, and it induces the pair of adjoint functors $\xymatrix{u^*: {\bf IS}(I) \ar@<-1.5pt>[r] & {\bf IS}(J): u_! \ar@<-1.5pt>[l]}$, which is defined by the following
$$\begin{array}{l} 
u^*(\Sigma ') = u^{-1}(\Sigma ') \ \hbox{for any} \ \Sigma ' \in {\bf IS}(I) \\ 
u_!(J_{\leq j}) = I_{\leq u(j)} \ \hbox{for any} \ j \in J  \hbox{, hence} \ u_!(\Sigma ) = \bigcup _{j \in \Sigma} I_{\leq u(j)} \ \hbox{for any} \ \Sigma \in {\bf IS}(J) \ .
\end{array}$$
The fact that the functor $u_!$ is left adjoint to the functor $u^*$ is equivalent with the following assertion: 
$$u_! ( \Sigma ) \subset \Sigma ' \ \Longleftrightarrow \ \Sigma \subset u^* (\Sigma ')$$
for any $\Sigma \in {\bf IS}(J)$ and any $\Sigma ' \in {\bf IS}(I)$.  
We always have the inclusion $u(\Sigma ) \subset u_!(\Sigma )$, and $u_!(\Sigma )$ is the smallest initial segment which contains $u(\Sigma )$. 

We can also define a functor $\xymatrix{u_*: {\bf IS}(J) \ar[r] & {\bf IS}(I)}$ which is right adjoint to $u^*$. For any $\Sigma \in {\bf IS}(J)$ its image is defined by $u_*(\Sigma ) = \{ i \in I \ | \ u^{-1}(I_{\leq i}) \subset \Sigma \}$, which is a consequence of the equality $u_*(\Sigma ) (i) = Hom _{{\bf IS}(J)}( u^*(I_{\leq i}) , \Sigma )$ (cf. \cite[Expos\'e I, Proposition 5.4, 2)]{SGA4}).  

\vskip .2cm 

\begin{remark}\label{rmq:dualite}
For any initial segment $\Sigma$ of $J$ we have also the equality $u_*(\Sigma ) = \bigcap _{j \in \Sigma ^{C}} I_{<u(j)}$, where $\Sigma ^{C}$ is the final segment $J \setminus \Sigma$ of $J$, in particular we have the equality $u_*(J_{<j}) = I_{<u(j)}$.

We can prove directly the equality $\{ i \in I \ | \ u^{-1}(I_{\leq i}) \subset \Sigma \} = \bigcap _{j \in \Sigma ^{C}} I_{<u(j)}$, which is equivalent to the following assertion: 
$$u^{-1} \bigl ( I_{\leq i} \bigr ) \subset \Sigma \ \Longleftrightarrow \  i < u(j) \ \forall\ j \in \Sigma ^{C} \ ,$$
but we can also deduce the equality $u_*(\Sigma ) = \bigcap _{j \in \Sigma ^{C}} I_{<u(j)}$ from the equality 
$I \setminus u_*(\Sigma ) = u_*(\Sigma ^{C})$,  
where the image $u_*(\Theta )$ of a final segment $\Theta$ of $J$ is the final segment $\bigcup _{j \in \Theta} I_{\geq u(j)}$ of $I$ (we study below the category {\bf FS}(I) of final sections and of the morphisms $\xymatrix{\underline h^I : I \ar[r] & {\bf FS}(I)}$). 
\end{remark}

\vskip .2cm 

Let $I$ be a totally ordered set, for any subset $Y$ of $I$ the least upper bound ${\rm sup}\, Y$ and the greatest lower bound ${\rm inf}\, Y$ of $Y$ in $I$ correspond respectively to the \emph{colimit} and to the \emph{limit} of the functor $\xymatrix @R=-1mm @C=8mm{u: Y\ar@<-1.5pt>@{^{(}->}[r] & I}$ induced by the inclusion.   
Then we say that a totally ordered set $I$ is \emph{cocomplete} (resp. \emph{complete}) if any subset $Y$ of $I$ admits a least upper bound (resp. a greatest lower bound), and we deduce from the correspondance stated in proposition \ref{prop:Yoneda} the following result. 

\begin{remark} \label{rmq:cocomplete=complete}
Let $I$ be a totally ordered set, then $I$ is complete if and only if it is cocomplete. 

It is enough to note that if $X$ is a subset of $I$, the element $x = {\rm inf}\, X$ is equal to the element $y = {\rm sup}\, Y$, where $Y$ is the subset of $I$ defined by 
$$Y = \{ i \in I \ | \ \forall \, j \in X \ i \leq j \} = \bigcap _{j \in X} I _{\leq j} \ .$$
\end{remark}

\begin{corollary}\label{cor:propriete-universelle} 
Let $I$ be a totally ordered set and let $\xymatrix{\underline h _I : I \ar[r] & {\bf IS}(I)}$ be the morphism defined above, then the ordered set ${\bf IS}(I)$ is cocomplete. 
Moreover for any morphism $\xymatrix{v : I \ar[r] & K}$ of ordered sets with $K$ cocomplete, there exists a unique morphism $\xymatrix{\hat v :{\bf IS}(I) \ar[r] & K}$ of ordered sets such that $\hat v \circ \underline h _I = v$. 
\end{corollary}

\begin{preuve} 
This is a direct consequence of proposition \ref{prop:Yoneda} and the classical result on categories, but we give an elementary proof. 

Let $Y$ be a subset of ${\bf IS}(I)$, we can write $Y = \{ \Sigma _{\alpha} | \alpha \in A \}$ for some set $A$, then the subset $\Sigma = \cup _{\alpha \in A} \Sigma _{\alpha}$ of $I$ is an initial segment, and it is easy to see that $\Sigma = {\rm sup}\, Y$ in ${\bf IS}(I)$.   

Let $\xymatrix{v : I \ar[r] & K}$ be a morphism with $K$ cocomplete, then we define the image of an initial segment $\Sigma \in {\bf IS}(I)$ in $K$ by $\hat v (\Sigma ) = {\rm sup}\, Y$ where $Y$ is the subset of $K$ defined by $Y = v(\Sigma )$.  

\hfill$\Box$ 
\end{preuve}  
 
\vskip .2cm 

\begin{lemma}\label{le:representable} 
An initial segment $\Sigma$ of an ordered set $I$ is \emph{representable}, i.e. is in the image of the Yoneda embedding $\xymatrix{\underline h _I : I \ar[r] & {\bf IS}(I)}$, if and only if $\Sigma$ belongs to the set ${\rm Suc} ({\bf IS}(I))$ of immediate successors of the ordered set ${\bf IS}(I)$. 
\end{lemma} 

\begin{preuve} 
This is another formulation of lemma \ref{le:predecesseur-successeur-dans-CpI}.

\hfill$\Box$ 
\end{preuve}

\vskip .2cm 

Let $\xymatrix{u: J \ar[r] & I}$ be a morphism of totally ordered sets that we consider as a functor between the associated categories, and let $\xymatrix{u^*: {\bf IS}(I) \ar@<-1.5pt>[r] & {\bf IS}(J): u_! \ar@<-1.5pt>[l]}$ and $\xymatrix{u^*: {\bf IS}(I) \ar@<+1.5pt>[r] & {\bf IS}(J): u_* \ar@<+1.5pt>[l]}$ be the two pairs of adjoint functors defined above, where the functors $u_!$ and $u_*$ are respectively the left adjoint and the right adjoint of the functor $u^*$. 
We deduce from the units and counits of these adjonctions that for any initial segment $\Sigma$ of $J$ and for any initial segment $\Sigma '$ of $I$ that we have the morphims $\xymatrix{u^* u_* (\Sigma )  \ar[r] & \Sigma   \ar[r] & u^* u_! (\Sigma )}$ and $\xymatrix{u_! u^* (\Sigma ') \ar[r] & \Sigma '  \ar[r] & u_* u^* (\Sigma ')}$, hence the inclusions: 

$$\begin{array}{l} 
u^* u_* (\Sigma ) \subset \Sigma  \subset u^* u_! (\Sigma ) \\
u_! u^* (\Sigma ') \subset \Sigma ' \subset u_* u^* (\Sigma ') \ .   
\end{array}$$

\vskip .2cm 

Let $\xymatrix{u: J \ar[r] & I}$ be an injective morphism of totally ordered sets, and we assume that the morphism $u$ is the inclusion of a subset $J$ of $I$, the morphism $\xymatrix{u^*: {\bf IS}(I) \ar[r] & {\bf IS}(J)}$ is defined by $u^*(\Sigma ') = \Sigma \cap J$ for any initial segment $\Sigma '$ of $I$.

\begin{lemma}\label{le:morphisme-injectif}
Let $\xymatrix{u: J \ar[r] & I}$ be an injective morphism of totally ordered sets, and let $\Sigma$ be an initial segment of $J$, then $u_!(\Sigma )$ is the smallest initial segment $\Sigma '$ of $I$ such that we have $\Sigma ' \cap J = \Sigma$ and $u_*(\Sigma )$ is the biggest initial segment $\Sigma '$ of $I$ such that we have $\Sigma ' \cap J = \Sigma$. 

More over an initial segment $\Sigma '$ of $I$ satisfies the equality $u^*(\Sigma ') = \Sigma$ if and only if we have the inclusions $u_!(\Sigma ) \subset \Sigma ' \subset u_*(\Sigma )$. 
\end{lemma} 

\begin{preuve} 
From the above any initial segment $\Sigma '$ of $I$ with $u^*(\Sigma ')=\Sigma$ satisfies the inclusions $u_! (\Sigma ) \subset \Sigma ' \subset u_* (\Sigma )$. 
 
As the morphism $u$ is injective, the associated functor $\xymatrix{u: J \ar[r] & I}$ is fully faithful, we deduce from proposition 5.6 Expos\'e I of \cite{SGA4} that the functors $\xymatrix{u_! , u_*: {\bf IS}(J) \ar[r] & {\bf IS}(I)}$ are also fully faithful, and we have the equalities $u^* u_* (\Sigma ) = \Sigma  = u^* u_! (\Sigma )$. 
Then the image $u^*(\Sigma ')$ of any initial segment $\Sigma '$ of $I$ with $u_! (\Sigma ) \subset \Sigma ' \subset u_* (\Sigma )$ is equal to $\Sigma$. 

\hfill$\Box$ 
\end{preuve}  

\vskip .2cm  

\begin{remark}\label{rmq:demonstration-directe}
If the morphism $\xymatrix{u:J \ar[r] & I}$ is an inclusion of totally ordered sets the functors $u_!$ and $u_*$ are defined by $u_!(\Sigma ) =\bigcup _{j \in \Sigma}I_{\leq j}$ and  $u_*(\Sigma ) =\bigcap _{j \in \Sigma ^{C}}I_{<j}$ for any $\Sigma$ in ${\bf IS}(J)$. 
Then we have directly  the equalities $u^*u_!(\Sigma ) = \bigcup _{j \in \Sigma}J_{\leq j} = \Sigma$ and $u^*u_*(\Sigma ) = \bigcap _{j \in \Sigma ^{C}}J_{<j} = \Sigma$. 
\end{remark}

\vskip .2cm 

If $\xymatrix{u: J \ar[r] & I}$ is a surjective morphism of totally ordered sets, the order $\leq$ on $I$ is completely determined by the order on $J$ as follows, 
let $i$ and $i '$ be in $I$ then 
$$i < i ' \ \Longleftrightarrow \ \forall \ j \ \hbox{and} \ j ' \in J \ \hbox{such that} \  u(j )=i  , \ u(j ') = i ' \  \hbox{we have} \ j  < j ' \ ,$$ 
hence the relation 
$$i \leq i ' \ \Longleftrightarrow \ \exists \ j \ \hbox{and} \ j ' \in I \ \hbox{such that} \  u(j )=i  , \ u(j ') = i  '\ \hbox{satisfying} \ j  \leq j ' \ .$$ 

In the case where the morphism $u$ is surjective the functor $u_!$ is defined by the image, i.e. $u(\Sigma ) = u_!(\Sigma )$. More generally we have the following result.

\begin{lemma}\label{le:morphisme-surjectif} 
Let $\xymatrix{u: J \ar[r] & I}$ be a surjective morphism of totally ordered sets.
The image $u(\Sigma )$ of any initial segment (resp. final segment, resp. interval) $\Sigma$ of $J$ by $u$ is an initial segment (resp. a final segment, resp. an interval) of $I$. 

Moreover if $\Sigma$ is an initial segment (resp. a final segment) which admits a greatest element (resp. a least element) $j$, $u(\Sigma )$ admits $i=u(j )$ as greatest element (resp. as least element).  
\end{lemma}

\begin{preuve} 
Let $\Sigma$ be an initial segment of $J$. 
Let $i$ and $i'$ be two elements of $I$ with $i ' \in u(\Sigma )$ and $i < i '$, then there exist $j$ and $j '$ in $J$ such that $j =u(i)$ and $j ' =u(i ')$, with $j ' \in \Sigma$. 
From $i < i '$ we deduce $j < j '$, and $j \in \Sigma$, hence $i \in u(\Sigma )$. 

\hfill$\Box$ 
\end{preuve}  

\vskip .2cm

In the same way, for any totally ordered set $I$, we can consider the category ${\bf FS}(I)$ as the category of \emph{copresheaves} on the category $I$, and the morphism $\xymatrix @R=-1mm @C=8mm{\sigma ^{\geq}: I \ar@<-1.5pt>@{^{(}->}[r] & {\bf FS}(I)}$ defined by $\sigma ^{\geq}(i) = I_{\geq i}$ 
corresponds to the embedding $\xymatrix{\underline h^I : I \ar[r] & {\bf FS}(I)}$ which send any element $i$ to the cosheaf $Mor_I(i,-)$. 
We can prove it directly or just notice that we have an isomorphism of ordered sets $\xymatrix@C=12mm{{\bf FS}(I) \ar[r] ^-{\simeq} & {{\bf IS}(I^{\rm op})}^{\rm op}}$. 

As for the category of initial segments, a morphism $\xymatrix{u:J \ar[r] & I}$  of ordered sets induces the pair of adjoint functors betwen the categories of final segments $\xymatrix{u^*: {\bf FS}(I) \ar@<-1.5pt>[r] & {\bf FS}(J): u_* \ar@<-1.5pt>[l]}$, where $u_*$ is the right adjoint functor of $u^*$, which is defined by the following
$$\begin{array}{l} 
u^*(\Sigma ') = u^{-1}(\Sigma ') \ \hbox{for any} \ \Sigma ' \in {\bf FS}(I) \\ 
u_*(J_{\geq j}) = I_{\geq u(j)} \ \hbox{for any} \ j \in J  \hbox{, hence} \ u_*(\Sigma ) = \bigcup _{j \in \Sigma} I_{\geq u(j)} \ \hbox{for any} \ \Sigma \in {\bf FS}(J) \ ,
\end{array}$$
and as above we can also define a functor $\xymatrix{u_!: {\bf FS}(J) \ar[r] & {\bf FS}(I)}$ which is left adjoint to $u^*$, and which is defined by $u_!(\Sigma ) = \{ i \in I \ | \ u^{-1}(I_{\geq i}) \subset \Sigma \} = \bigcap _{j \in \Sigma ^{C}} I_{> u(j)} $, for any $\Sigma \in {\bf FS}(J)$.  

\vskip .2cm 

The category of cuts ${\bf Cp}(I)$ is equivalent to the two categories ${\bf IS}(I)$ of initial segments and ${\bf FS}(I)$ of final segments, but the induced morphisms $\sigma ^{\leq}$ and $\sigma ^{\geq}$ from $I$ to ${\bf Cp}(I)$ are different, and we will define different pairs of functors. 

\begin{proposition}\label{prop:foncteurs-afjoints} 
Let $\xymatrix{u:J \ar[r] & I}$ be a morphism of totally ordered sets, then there exists one functor $\xymatrix{u^*: {\bf Cp}(I) \ar[r] &{\bf Cp}(J)}$, defined for any cut $\Lambda '= (\Lambda '_- , \Lambda '_+)$ of $I$ by:  
$$u^* ( \Lambda ') = ( u^*( \Lambda '_-) , u^*( \Lambda '_+) ) = ( u^{-1}( \Lambda '_-) , u^{-1}( \Lambda '_+) ) \ ,$$ 
and there exist two functors $\xymatrix{u_!: {\bf Cp}(J) \ar[r] &{\bf Cp}(I)}$ and $\xymatrix{u_*: {\bf Cp}(J) \ar[r] &{\bf Cp}(I)}$, defined for any cut $\Lambda = (\Lambda _- , \Lambda _+)$ of $J$ by:  
$$\begin{array}{l} 
u_! ( \Lambda ) = (u_!(\Lambda _-), u_!(\Lambda _+)) \\ 
u_* ( \Lambda ) = (u_*(\Lambda _-), u_*(\Lambda _+))  \  . 
\end{array}$$
Moreover the functors $u_!$ and $u_*$ are respectively the left adjoint and the right adjoint functors of the functor $u^*$. 
\end{proposition} 

\begin{preuve} 
Let $\Lambda _-$ be an initial segment of $J$, we recall that $u_! (\Lambda _- )$ is the initial segment of $I$ defined by $u_! (\Lambda _-) = \cup _{j \in \Lambda _-} I_{\leq u(j)}$, and the image $u_! (\Lambda _+)$ of the final segment $\Lambda _+ = J \setminus \Lambda _-$ is equal to $ \{ i \in I \ | \ u^{-1}(I_{\geq i}) \subset \Lambda _+ \}$. 

Then for any element $i$ in $I$ we have the following equivalences: 
$$\begin{array}{rcl}
i \in u_! (\Lambda _-) & \Longleftrightarrow & \exists j \in \Lambda _- \  \hbox{such that} \  i \leq u(j) \\ 
i \in u_! (\Lambda _+) & \Longleftrightarrow & \forall j \in J \  \hbox{such that} \  u(j) \geq i \  \hbox{we have} \  j \in \Lambda _+ \ , 
\end{array}$$
and the pair $(u_!(\Lambda _-), u_!(\Lambda _+))$ is a cut of $I$. 

In the same way we prove that the pair $(u_*(\Lambda _-), u_*(\Lambda _+))$ is also a cut of $I$. 

\hfill$\Box$ 
\end{preuve}  

\vskip .2cm

If the morphism $\xymatrix{u: J \ar[r] & I}$ is surjective, for any cut $\Lambda = (\Lambda _-,\Lambda _+)$ of $J$, we have the equalities $u_!(\Lambda _-)=u(\Lambda _-)$ and $u_*(\Lambda _+)=u(\Lambda _+)$ (cf. lemma \ref{le:morphisme-surjectif}), 
and we define the pair $u(\Lambda )$ of subsets of $I$ by $u(\Lambda )= (u(\Lambda _-),u(\Lambda _+))$. 

\begin{lemma}\label{le:comparaison} 
Let $\xymatrix{u: J \ar[r] & I}$ be a surjective morphism of totally ordered sets, then for any cut $\Lambda$ of $J$ we have the inequality $u_*(\Lambda ) \leq u_!(\Lambda )$ in ${\bf Cp}(I)$.  

We have equality $u_*(\Lambda ) = u_!(\Lambda )$ if and only if we have $u(\Lambda _-) \cap u(\Lambda _+) = \emptyset$, in that case the pair $u(\Lambda )$ is a cut of $I$. 
If we have $u_*(\Lambda ) \not= u_!(\Lambda )$ there exists one element $i$ in $I$ such that we have the equality $u(\Lambda _-) \cap u(\Lambda _+) = \{ i \}$.  
\end{lemma} 

\begin{preuve} 
We always have $u_!(\Lambda _-) \cup u_*(\Lambda _+) = u(\Lambda _-) \cup u(\Lambda _+) = I$, hence  $u_*(\Lambda ) \leq u_!(\Lambda )$, and $u(\Lambda )$ is a cut if $u(\Lambda _-) \cap u(\Lambda _+) =\emptyset$. 

An element $i$ belongs to the intersection $u(\Lambda _-) \cap u(\Lambda _+)$ if there exist $j_1$ and $j_2$ with $j_1 \in \Lambda _-$ and $j_2 \in \Lambda _+$ such that $u(j_1)=u(j_2)=i$. 
If $i_1$ and $i_2$ are two element of $u(\Lambda _-) \cap u(\Lambda _+)$ satisfying $i_1 < i_2$, for any $j_1$ and $j_2$ in $J$ with $u(j_1)=i_1$ and $u(j_2)=i_2$, we must have $j_1 < j_2$, which is impossible. 

\hfill$\Box$ 
\end{preuve}  

\vskip .2cm

\subsection{Complete and cocomplete ordered sets}

We assume that the ordered set $I$ is cocomplete, then we can define a morphism of ordered sets $\xymatrix{\lambda _- : {\bf Cp}(I) \ar[r] & I}$ by $\lambda _- = {\rm sup}\, \Lambda _-$ for any cut $\Lambda$, and in the same way if $I$ is complete we can define the morphism $\xymatrix{\lambda _+ : {\bf Cp}(I) \ar[r] & I}$ by $\lambda _+ = {\rm inf}\, \Lambda _+$. We recall that we have defined the two morphisms $\xymatrix{\sigma ^{\leq}: I \ar@<-1.5pt>@{^{(}->}[r]& {\bf Cp}(I)}$ and $\xymatrix{\sigma ^{\geq}: I \ar@<-1.5pt>@{^{(}->}[r]& {\bf Cp}(I)}$. 

\begin{proposition}\label{prop:composition}
{(1)} Let $I$ be a cocomplete ordered set, then we have the following equalities: 
\begin{enumerate} 
\item {$(\lambda _- \circ \sigma ^{\leq}) = id_I$;}
\item {$(\lambda _- \circ \sigma ^{\geq}) (i) = i$ if $i \notin {\rm Suc} (I)$ and $(\lambda _- \circ \sigma ^{\geq}) (i) = {\rm pre} (i)$ if $i \in {\rm Suc} (I)$.} 
\end{enumerate} 

Moreover for any cut $\Lambda$ we have the inequalities $(\sigma ^{\geq} \circ \lambda _-)(\Lambda ) \leq \Lambda \leq (\sigma ^{\leq} \circ \lambda _-)(\Lambda )$, with $\Lambda = (\sigma ^{\geq} \circ \lambda _-)(\Lambda ) = (I_{< \lambda _-},I_{\geq \lambda _-})$ if $\lambda _-= \lambda _-(\Lambda )\in \Lambda _+$ and $\Lambda = (\sigma ^{\leq} \circ \lambda _-)(\Lambda ) = (I_{\leq \lambda _-},I_{> \lambda _-})$ if $\lambda _- \in \Lambda _-$. 
In any case the cut $\Lambda$ is principal.

\vskip .2cm 

{(2)} Let $I$ be a complete ordered set, then we have the following equalities: 
\begin{enumerate} 
\item {$(\lambda _+ \circ \sigma ^{\geq}) = id_I$;}
\item {$(\lambda _+ \circ \sigma ^{\leq}) (i) = i$ if $i \notin {\rm Pre} (I)$ and $(\lambda _+ \circ \sigma ^{\leq}) (i) = {\rm suc} (i)$ if $i \in {\rm Pre} (I)$.} 
\end{enumerate}  

Moreover for any cut $\Lambda$ we have the inequalities $(\sigma ^{\geq} \circ \lambda _+)(\Lambda ) \leq \Lambda \leq (\sigma ^{\leq} \circ \lambda _+)(\Lambda )$, with $\Lambda = (\sigma ^{\geq} \circ \lambda _+)(\Lambda ) = (I_{< \lambda _+},I_{\geq \lambda _+})$ if $\lambda _+ = \lambda _+(\Lambda )\in \Lambda _+$ and $\Lambda = (\sigma ^{\leq} \circ \lambda _+)(\Lambda ) = (I_{< \lambda _+},I_{\geq \lambda _+})$ if $\lambda _+ \in \Lambda _-$. 
In any case the cut $\Lambda$ is principal.

\vskip .2cm 

{(3)} Let $I$ be a complete and cocomplete ordered set. 
then for any cut $\Lambda$ we have the inequality $\lambda _- \leq \lambda _+$, and we have $\lambda _- \not= \lambda _+$ if and only if the cut $\Lambda$ is a jump. 
In this case $\lambda _+$ is the immediate successor of $\lambda _-$, and we have $\Lambda = (I_{\leq \lambda _-} , I_{\geq \lambda _+})$. 
\end{proposition} 

\begin{preuve} 
This is obvious from the definitions.  

\hfill$\Box$ 
\end{preuve}  

\vskip .2cm

Let $\Sigma$ be an initial segment in a cocomplete ordered set $I$, and let $x$ in $I$ defined by $x = {\rm sup}\, (\Sigma )$, then we have the following possibilities: 
\begin{enumerate}
\item \emph{if $x \in {\rm Suc}(I)$, then $\Sigma = I_{\leq x}$;}  
\item \emph{if $x \notin {\rm Suc}(I)$, then either $\Sigma = I_{\leq x}$ or $\Sigma = I_{< x}$.} 
\end{enumerate} 

\vskip .2cm

\begin{definition}
Let $I$ be a totally ordered set, a subset $J$ of $I$ is \emph{dense above} in $I$ if for any elements $x_1$ and $x_2$ in $I$ with $x_1 < x_2$ there exists $y$ in $J$ such that $x_1 < y \leq x_2$.   
\end{definition} 

\begin{lemma}\label{le:dense-above} 
Let $I$ be a cocomplete totally ordered set, let $J$ be a dense above subset in $I$, then for any $x$ in $I$ we have the equality $x = {\rm sup}\, (I_{\leq x} \cap J)$.  
\end{lemma} 

\begin{preuve} 
For any $x$ in $I$ we have the inequality $x' =  {\rm sup}\, (I_{\leq x} \cap J) \leq x$, and if the inequality was strict, $x' < x$, there would exist $y$ in $J$ with $x' < y \leq x$, which is impossible. 

\hfill$\Box$ 
\end{preuve}  

\vskip .2cm 

\begin{lemma}\label{le:successeur-dense-above}
If $J$ is a dense above subset of a totally ordered set $I$, it contains the subset ${\rm Suc}(I)$ of immediate successors. 
\end{lemma} 
\begin{preuve} 
Let $x_2$ be an element of ${\rm Suc}(I)$ and let $x_1$ be the element of $I$ such that $x_2 = {\rm suc}(x_1)$, then we have $x_1 < x_2$ and $I_{> x_1} \cap I_{< x_2} = \emptyset$. By hypothesis on $J$ there exists $y \in J$ such that $x_1 < y \leq x_2$, hence $y=x_2$, and $x_2 \in J$. 

\hfill$\Box$ 
\end{preuve}  

\vskip .2cm

\begin{proposition}\label{prop:natural-isomorphism}  
Let $I$ be a cocomplete totally ordered set, let $J$ be a subset of $I$ such that $J$ is dense above in $I$ and contained in ${\rm Suc}(I)$. 
Then there is a natural isomorphism $\xymatrix{f: I \ar[r] & {\bf IS}(J)}$ compatible with the inclusion of $J$ in $I$ and the Yoneda embedding ${\underline h _J}$. 
\end{proposition} 

\begin{preuve} 
Let $\xymatrix{u: J \ar@<-1.5pt>@{^{(}->}[r] & I}$ be the injection of $J$ in $I$, we consider this morphism as a functor and with the previous notations we have the diagram 
$$\xymatrix@R=10mm @C=18mm{ 
J \ar[r] ^u \ar[d] ^{\underline h _J} & I \ar[d] ^{\underline h _I} \\ 
{\bf IS}(J) \ar@<+1.5pt>[r] ^{u_!}  & {\bf IS}(I) \ar@<+1.5pt>[l]^{u^*} }$$ 
where $u^*$ is defined by $u^*(\Sigma ) = \Sigma \cap J$ and where $u_!$ is the left adjoint functor of $u^*$. 
We have the equality ${\underline h _I} \circ u = u_! \circ {\underline h _J}$, moreover as the morphism $u$ is an injection, the associated functor is fully faithful and we have also the equality $\xymatrix{ {\underline h _J} = u^* \circ {\underline h _I} \circ u :J \ar[r] & {\bf IS}(J)} $. 

We define the morphism of ordered sets $f$ by $\xymatrix{f= u^* \circ {\underline h _I} :I \ar[r] & {\bf IS}(J)}$, which satisfies $f \circ u = {\underline h _J}$, and we want to show it is an isomorphism. 

Let $\xymatrix{g= {\bf IS}(J) \ar[r] & I}$ defined by $g(\Sigma ) = {\rm sup}\, (u_!(\Sigma ))$, with $u_!(\Sigma ) = \bigcup _{j \in \Sigma} I_{\leq u(j)}$, hence we also have the equality $g(\Sigma ) = {\rm sup}\, (\Sigma )$ (it is also a consequence of the fact that the functor $u_!$ is a left adjoint and that the ${\rm sup}$ is a colimit). 
By construction for any $i$ in $I$ we have $g \circ f (i) = {\rm sup}\, (I_{\leq i} \cap J)$, and as $J$ is dense above in $I$ we deduce from lemma \ref{le:dense-above} that $g \circ f (i) =i$. Hence we have the equality $g \circ f = {\rm id}_I$, and the morphism $f$ is injective. 

To show that the morphism $f$ is an isomorphism with $g = f^{-1}$, we must show that any initial segment $\Sigma$ in $J$ is equal to $f(x) = I_{\leq x} \cap J$, where $x$ is defined by $x=g(\Sigma ) = {\rm sup}\, (u_!(\Sigma ))$.  
We have either $u_!(\Sigma ) = I_{\leq x}$ and $\Sigma = I_{\leq x} \cap J$ if $x \in \Sigma$, or $u_!(\Sigma ) = I_{< x}$ and $\Sigma = I_{< x} \cap J$ if $x \notin \Sigma$, then for $x \in \Sigma$ or for $x \notin J$ it is obvious that we have $\Sigma = I_{\leq x} \cap J$. 

We assume that we have $x \in J$ and $x \notin \Sigma$, then by hypothesis there would exist $y \in I$ such that $x={\rm suc}(y)$, hence $I_{< x} = I_{\leq y}$, and we would have the equality $\Sigma = I_{\leq y} \cap J$, which is impossible from lemma\ref{le:dense-above}. 

\hfill$\Box$ 
\end{preuve}  

\vskip .2cm 

\begin{remark}\label{rmq:successeur-dense-above} 
We deduce from lemma \ref{le:successeur-dense-above} that a subset $J$ of a totally ordered set $I$ which satisfies the hypothesis of the proposition \ref{prop:natural-isomorphism} is exactly the set ${\rm Suc}(I)$. 
This is also a direct consequence of the isomorphism $\xymatrix{f: I \ar[r] & {\bf IS}(J)}$ of proposition \ref{prop:natural-isomorphism} and of lemma \ref{le:representable}. 
\end{remark}

\vskip .2cm 

\begin{corollary}\label{cor:restriction} 
We assume that $J$ is a subset of a totally ordered set $I$ which satisfies the hypothesis of the proposition \ref{prop:natural-isomorphism}, then the natural morphism $\xymatrix{u^* : {\bf IS}(I) \ar[r] & {\bf IS}(J)}$ which sends an initial segment $\Sigma$ of $I$ to its restriction $\Sigma \cap J$ is an epimorphism. 

Two distinct initial segments $\Sigma \subsetneq \Sigma '$ of $I$ satisfy $\Sigma \cap J = \Sigma ' \cap J$ if and only if there exists $x \notin J$ with $\Sigma = I_{< x}$ and $\Sigma ' = I_{\leq x}$. 

Moreover we have the equality ${\rm sup}\, (u^*(\Sigma )) = {\rm sup}\, (\Sigma )$. 
\end{corollary} 

\begin{preuve}
By proposition \ref{prop:natural-isomorphism}, there exists an isomorphism $\xymatrix{f :I \ar[r] & {\bf IS}(J)}$ with $f = u^* \circ {\underline h _I}$, then the morphism $u^*$ is an epimorphism.  

Let $\Sigma$ be an initial segment of $I$, let $x_1$ and $x_2$ be the elements of $I$ defined respectively by $x_1 = {\rm sup}\, (u^*(\Sigma ))$ and $x_2 = {\rm sup}\, (\Sigma )$. We deduce from the equality $u^*(\Sigma ) = \Sigma \cap J$ the inequality $x_1 \leq x_2$, and as $J$ is dense above in $I$ we cannot have the strict inequality $x_1 < x_2$.  

\hfill$\Box$ 
\end{preuve}  

\vskip .2cm 

%
%
%

\begin{remark}\label{rmq:comparaison} 
Let $I$ be a complete and cocomplete ordered set, let $J$ be a subset of $I$. 
For any cut $\Lambda =(\Lambda _-,\Lambda _+)$ of $I$, we want to compare $\lambda _- ={\rm sup}\, \Lambda _-$ and $\lambda _- ^{\scriptscriptstyle (J)} = {\rm sup} (\Lambda _- \cap J)$, and $\lambda _+ = {\rm inf}\, \Lambda _+$ and $\lambda _+^{\scriptscriptstyle (J)} ={\rm inf} (\Lambda _+ \cap J)$. 
We always have the inequalities $\lambda _- ^{\scriptscriptstyle (J)} \leq \lambda _- \leq \lambda _+ \leq \lambda _+ ^{\scriptscriptstyle (J)}$.

We assume now that the subset $J$ satisfies the hypothesis of the proposition \ref{prop:natural-isomorphism}, then we deduce from corollary \ref{cor:restriction}  the equality $\lambda _- ^{\scriptscriptstyle (J)} = \lambda _-$. 
Moreover if we have the strict inequality $\lambda _+ < \lambda _+ ^{\scriptscriptstyle (J)}$, we have $ \lambda _+ ^{\scriptscriptstyle (J)} = {\rm suc}(\lambda _+)$. 
Then there are only three possibilities: 
\begin{enumerate} 
\item $\lambda _- ^{\scriptscriptstyle (J)} = \lambda _- = \lambda _+ = \lambda _+ ^{\scriptscriptstyle (J)}$; 
\item $\lambda _- ^{\scriptscriptstyle (J)} = \lambda _- = \lambda _+ < \lambda _+ ^{\scriptscriptstyle (J)} = {\rm suc}(\lambda _+)$,  
\item $\lambda _- ^{\scriptscriptstyle (J)} = \lambda _- < \lambda _+ = \lambda _+ ^{\scriptscriptstyle (J)} = {\rm suc}(\lambda _-)$.
\end{enumerate}
\end{remark} 

\vskip .2cm 

We can also deduce from the above a characterization of the totally ordered sets which are sets of initial segments, or cuts. This can be seen as an analogue of a characterization of the categories which are categories of presheaves. 

\begin{proposition}\label{prop:caracterisation}
Let $I$ be a totally ordered set, then the following assertions are equivalent.
\begin{enumerate}
\item There exists a totally ordered set $J$ such that $I$ is isomorphic to the set ${\bf IS}(J)$ of the initial segments of $J$. 
\item The set $I$ is cocomplete and the subset ${\rm Suc}(I)$ of immediate successors of $I$ is dense above in $I$. 
\end{enumerate}  
\end{proposition}

\begin{preuve}
We deduce from proposition \ref{prop:natural-isomorphism} that if $I$ is cocomplete and if the subset $J={\rm Suc}(I)$ is dens above in $I$, there exists an isomorphism $\xymatrix{f :I \ar[r] & {\bf IS}(J)}$. 

Conversely let $I$ be the set ${\bf IS}(J)$ of initial segments of a set $J$, then we deduce from corollary \ref{cor:propriete-universelle} that $I$ is complete, and we deduce from lemma \ref{le:representable} that the set ${\rm Suc}(I)$ of immediate successors of $I$ is the set of representable initial segments, i.e. of initial segment of the form $J_{\leq j}$, for $j \in J$. 

Let $\Sigma$ and $\Sigma '$ be two initial segments of $J$  with $\Sigma \subsetneq \Sigma '$, then for any element $j$ in $\Sigma ' \setminus \Sigma$ we have the inclusions $\Sigma \subsetneq J_{\leq j} \subset \Sigma '$, hence the result. 

\hfill$\Box$ 
\end{preuve}  

\vskip .2cm

                     \section{$I$-structure on an $R$-module} 
%

\subsection{$I$-structure}

In this part we will essentially recall the results of P. Conrad in \cite{Co}. 
We give ourselves a totally ordered space $I$ and a commutative ring $R$, unless otherwise stated all the modules that we will consider in this paragraph will be modules over the ring $R$.

\begin{definition} 
An \emph{$I$-structure} on an $R$-module $\Gamma$ is the data of a family $\Delta _I (\Gamma ) = \bigl ( \delta (\Gamma) _i ^-, \delta (\Gamma ) _i ^+ \bigr ) _{i \in I}$ of submodules of $\Gamma$ satisfying: 
\begin{enumerate} 
\item $\delta (\Gamma )_i ^- \subset \delta (\Gamma )_i ^+$ for any $i$ in $I$; 
\item $\delta (\Gamma )_j^+ \subset \delta (\Gamma )_i^-$ for any $i<j$ in $I$; 
\item for any $\xi$ belonging to $\Gamma$, $\xi \not= 0$,  there exists $i$ in $I$ such that $\xi \in \delta (\Gamma )_i^+$ and $\xi \notin \delta (\Gamma )_i^-$. 
\end{enumerate} 
We then say that $\Gamma$ is a \emph{$I$-module}. 
Moreover for any $i$ in $I$ we note $\varepsilon (\Gamma )_i$ the quotient module $\varepsilon (\Gamma )_i = \delta (\Gamma )_i ^+ / \delta (\Gamma )_i ^-$, and we call the pair $\bigl ( I , (\varepsilon (\Gamma )_i) _{i \in I} \bigr )$ the \emph{skeleton} of $\Gamma$.  

If for any $i$ in $I$ we have a strict inclusion $\delta (\Gamma )_i ^- \subsetneq \delta (\Gamma )_i ^+$, i.e. $\varepsilon (\Gamma )_i \not= (0)$, we say that the $I$-structure is \emph{non degenerated}. 
\end{definition} 

A family $\Delta _I (\Gamma ) = \bigl ( \delta (\Gamma) _i ^-, \delta (\Gamma ) _i ^+ \bigr ) _{i \in I}$ of submodules of $\Gamma$ which only satisfies the properties (1) and (2) of the previous definition is called an  \emph{$I$-prestructure}. 
Moreover if for any $i \in I$ we have the equality $\delta (\Gamma )_i^- = \delta (\Gamma )_i^+$, i.e. $\varepsilon (\Gamma )_i = (0)$, we say that the $I$-prestructure is \emph{totally degenerated}. 

If the family $\Delta _I (\Gamma ) = \bigl ( \delta (\Gamma) _i ^-, \delta (\Gamma ) _i ^+ \bigr ) _{i \in I}$ is an $I$-structure on $\Gamma$, any subfamily $\Delta _J (\Gamma ) = \bigl ( \delta (\Gamma) _j ^-, \delta (\Gamma ) _j ^+ \bigr ) _{j \in J}$, where $J$ is a subset of $I$ is a $J$-prestructure on $\Gamma$.

\vskip .2cm 

Let $\Delta _I (\Gamma ) = \bigl ( \delta (\Gamma) _i ^-, \delta (\Gamma ) _i ^+ \bigr ) _{i \in I}$ be an $I$-structure on a $R$-module $\Gamma$, we define the subset $I _{\rm nd}(\Gamma )$ of $I$ by 
$$I_{\rm nd}(\Gamma ) \ = \ \{ i \in I \ | \ \varepsilon (\Gamma )_i \not= (0) \}  \ = \ \{ i \in I \ | \ \exists \xi \in \Gamma \hbox{ such that $\xi \in \delta (\Gamma )_i^+$ and $\xi \notin \delta (\Gamma )_i^-$} \} \ .$$ 
The subset $I_{\rm nd}(\Gamma )$ depends on the $I$-module $\Gamma$ but when there is no ambiguity  we denote it $I_{\rm nd}$. 
Then the subfamily $\bigl ( \delta (\Gamma) _i ^-, \delta (\Gamma ) _i ^+ \bigr ) _{i \in I_{\rm nd}}$ of $\Delta _I (\Gamma )$ defines an $I_{\rm nd}$-structure on $\Gamma$, which is non degenerated. 

The subset $I_{\rm nd}(\Gamma )$ of $I$ is the smallest subset $I'$ of $I$ such that the $I$-structure $\Delta _I(\Gamma )$ induces an $I'$-structure on the $R$-module $\Gamma$. 
More precisely if  $I'$ is a subset of $I$ the subfamily $\bigl ( \delta (\Gamma) _i ^-, \delta (\Gamma ) _i ^+ \bigr ) _{i \in I'}$ of $\Delta _I(\Gamma )$ defines an $I'$-structure on $\Gamma$ if and only if we have the inclusion $I_{\rm nd}(\Gamma ) \subset I'$.   

\vskip .2cm

\begin{definition}\label{def:R=Z-R=F}
If the ring $R$ is the ring of integers $\Z$ the category of $R$-modules is isomorphic to the category of  abelian groups and we speak of an $I$-structure on a group $\Gamma$ and say that $\Gamma $ is an \emph{$I$-group}. 

If the ring $R$ is a field $\F$ we speak of an $I$-structure on a $\F$-vector space $\Gamma$ and say that $\Gamma $ is an \emph{$I$-vector space}. 
\end{definition} 

\vskip .2cm

For any $\xi \not= 0$ belonging to $\Gamma$ the element $i$ of $I$ satisfying $\xi \in \delta (\Gamma )_i^+$ and $\xi \notin \delta (\Gamma )_i^-$ belongs to $I_{\rm nd}$ and is unique. 
We denote $\iota _{\Delta}(\xi )$ this element of $I$ and we have the following equivalences: $$\xi \in \delta (\Gamma ) _j^- \, \Longleftrightarrow \, j < \iota _{\Delta}(\xi ) \quad\hbox{and}\quad \xi \in \delta (\Gamma ) _j^+ \, \Longleftrightarrow \, j \leq \iota _{\Delta}(\xi ) \ .$$

The image of the element $\xi$ of $\Gamma$ in the quotient module $\varepsilon (\Gamma )_{\iota _{\Delta}(\xi )}$ is non zero, and is called the \emph{initial part} of $\xi$, we denote this element ${\bf in}_{\Delta}(\xi )$.

\vskip .2cm

We deduce from the definition that we have the equalities 
$$\bigcap _{i \in I} \delta (\Gamma )_i^- = (0) \quad\hbox{and}\quad \bigcup _{i \in I} \delta (\Gamma )_i^+ = \Gamma \ .$$

Moreover if $I_{\rm nd}$ has a greatest element $i_1$ we have $\delta (\Gamma )_i^- = (0)$ for any $i \geq i_1$, 
and we have the equality $\displaystyle{\bigcap _{i \in I} \delta (\Gamma )_i^+ = \delta (\Gamma )_{i_1}^+} \not= (0)$, 
otherwise $\displaystyle{\bigcap _{i \in I} \delta (\Gamma )_i^+ = (0)}$. 

In the same way if $I_{\rm nd}$ has a smallest element $i_2$ we have $\delta (\Gamma )_i^+ = \Gamma$ for any $i \leq i_2$, 
and we have the equality $\displaystyle{\bigcup _{i \in I} \delta (\Gamma )_i^- = \delta (\Gamma )_{i_2}^-} \not= \Gamma$, 
otherwise $\displaystyle{\bigcup _{i \in I} \delta (\Gamma )_i^- = \Gamma}$.

\begin{lemma}\label{le:egalites} 
For any $i$ in $I$ we have the equalities 
$$\delta (\Gamma )_i ^- = \bigcup _{j>i} \delta (\Gamma )_j ^+ \quad\hbox{and}\quad \delta (\Gamma )_i ^+ = \bigcap _{j<i} \delta (\Gamma )_j ^-\ .$$
\end{lemma} 

\begin{preuve} 
We have the following equivalences: 
$$\begin{array}{rcl} 
\xi \in \bigcup _{j>i} \delta (\Gamma )_j ^+ & \Longleftrightarrow & \exists \, j > i \quad\hbox{such that }\quad j \leq \iota _{\Delta}(\xi ) \\ 
& \Longleftrightarrow & i < \iota _{\Delta}(\xi ) \\ 
& \Longleftrightarrow & \xi \in \delta (\Gamma )_i^- \ ; \\ \\
\xi \in \bigcap _{j<i} \delta (\Gamma )_j ^- & \Longleftrightarrow & \forall \, j < i \, , \quad j < \iota _{\Delta}(\xi ) \\ 
& \Longleftrightarrow & i \leq \iota _{\Delta}(\xi ) \\ 
& \Longleftrightarrow & \xi \in \delta (\Gamma )_i^+ \ .
\end{array}$$

\hfill$\Box$ 
\end{preuve}

\vskip .2cm

\begin{remark}\label{rmq:successeur2} 
If the element $i$ in $I$ admits an immediate successor $j={\rm suc}(i)$, we deduce from the previous lemma that we have the equality $\delta (\Gamma )_i^- = \delta (\Gamma )_j ^+$. 
And in the same way if $i$ admits an immediate predecessor $j={\rm pre}(i)$, we have the equality $\delta (\Gamma )_i^+ = \delta (\Gamma )_j ^-$. 

Hence if the set $I$ is well-ordered all the modules $\delta (\Gamma )_i^-$ are of the form $\delta (\Gamma )_j ^+$ for some $j$. 
\end{remark} 

\vskip .2cm

Let $I$ and $J$ be two totally ordered sets with $I \subset J$ and let $\Gamma$ be a $R$-module endowed with an $I$-structure $\Delta _I (\Gamma ) = \bigl ( \delta (\Gamma) _i ^-, \delta (\Gamma ) _i ^+ \bigr ) _{i \in I}$.  
For any $j \in J$ we define the submodules $\delta (\Gamma ) _j ^-$ and $\delta (\Gamma ) _j ^+$ by  
$$\delta (\Gamma ) _j^- =  \bigcup _{i_1 > j} \delta (\Gamma ) _{i_1}^+ \quad\hbox{and}\quad   \delta (\Gamma ) _j^+ = \bigcap _{i_2 < j} \delta (\Gamma ) _{i_2} ^- \ ,$$ 
and we notice that according to lemma \ref{le:egalites} if $j$ is an element belonging to $I$ we find clearly the submodules $\delta (\Gamma ) _j ^-$ and $\delta (\Gamma ) _j ^+$ belonging to the $I$-structure $\Delta _I(\Gamma )$. 
We thus define a family $\Delta _J (\Gamma ) = \bigl ( \delta (\Gamma) _j ^-, \delta (\Gamma ) _j ^+ \bigr ) _{j \in J}$, that we call the \emph{extension of the family $\Delta _I(\Gamma )$ to $J$}.  

If there exists an element $j$ in $J$ satisfying $j>i$ for any $i$ belonging to $I$, we have he equalities 
$$\delta (\Gamma ) _j^- =  \bigcup _{\emptyset} \delta (\Gamma ) _{i}^+ = (0) \quad\hbox{and}\quad  \delta (\Gamma ) _j^+ = \bigcap _{I} \delta (\Gamma ) _{i} ^- = (0)\ ,$$ 
and in the same way if there exists an element $j$ in $J$ satisfying $j<i$ for any $i$ belonging to $I$, we have the equalities 
$$\delta (\Gamma ) _j^- =  \bigcup _{I} \delta (\Gamma ) _{i}^+ = \Gamma \quad\hbox{and}\quad  \delta (\Gamma ) _j^+ = \bigcap _{\emptyset} \delta (\Gamma ) _{i} ^- = \Gamma \ .$$ 

We say that the set $I$ is \emph{dense} in the set $J$ if for any pair of elements  $j<j'$ in $J$ there exists an element $i$ in $I$ satisfying $j \leq i \leq j'$. 
In particular if $j$ and $j'$ don't belong to $I$ we deduce from it that there exists $i$ in $I$ such that $j<i<j'$. 

\begin{proposition} \label{prop:extension-de-I-a-J} 
Let $I \subset J$ be two totally ordered sets with $I$ dense in $J$, let $\Delta _I (\Gamma ) = \bigl ( \delta (\Gamma) _i ^-, \delta (\Gamma ) _i ^+ \bigr ) _{i \in I}$ be an $I$-structure on a $R$-module $\Gamma$. Then the family $\Delta _J (\Gamma ) = \bigl ( \delta (\Gamma) _j ^-, \delta (\Gamma ) _j ^+ \bigr ) _{j \in J}$ constructed as the extension of the family $\Delta _I(\Gamma )$ to $J$  defines an $J$-structure on $\Gamma$.  
\end{proposition} 

\begin{preuve} 

For any $i_1$ and $i_2$ in $I$ satisfying $i_2 < j < i_1$ we have $\delta (\Gamma )_{i_1}^+ \subset \delta (\Gamma )_{i_2}^-$, therefore we have the inclusion 
$\delta (\Gamma ) _j^- \subset \delta (\Gamma ) _j^+ $. 

From the definition of the submodules  $\delta (\Gamma ) _j ^-$ and $\delta (\Gamma ) _j ^+$ we check that for any pair $(i,j)$ with $i \in I$ and $j \in J$ we have $\delta (\Gamma ) _j^+ \subset \delta (\Gamma ) _i ^-$ for $i<j$ and $\delta (\Gamma ) _i ^+ \subset \delta (\Gamma ) _j^- $ for $i>j$.  
We deduce from it that for any $j<j'$ in $J$, as $I$ is dense in $J$ there exists $i$ in $I$ belonging to the interval $[j,j']$ we still have the inclusion $\delta (\Gamma ) _j^+ \subset \delta (\Gamma )_{j'}^-$. 

For any $\xi$ belonging to $\Gamma$ there exists $i$ in $I$ such that $\xi \in \delta (\Gamma ) _i ^+$ and $\xi \notin \delta (\Gamma )_i^-$. 

\hfill$\Box$ 
\end{preuve} 

\vskip .2cm

\begin{remark}\label{rmq:extension} 
The $J$-structure obtained as extension of $\Delta _I(\Gamma )$ satisfies the equality $\delta (\Gamma ) _j ^- = \delta (\Gamma ) _j ^+$ for any $j$ belonging to $J \setminus I$. 
Especially we have the equality $J_{\rm nd}(\Gamma ) = I_{\rm nd}(\Gamma )$. 
\end{remark} 

\vskip .2cm

\begin{corollary}\label{cor:extension-de-I-structure} 
Let $I$ be a totally ordered set and ${\bf Cp}(I)$ be the ordered set of cuts of $I$,  and let $\Delta _I (\Gamma ) = \bigl ( \delta (\Gamma) _i ^-, \delta (\Gamma ) _i ^+ \bigr ) _{i \in I}$ be an $I$-structure on a $R$-module $\Gamma$. 
Then we can define two ${\bf Cp}(I)$-structures on $\Gamma$ from the two embeddings $\xymatrix @C=5mm {\sigma ^{\leq}: I \ar@<-1.5pt>@{^{(}->}[r] & {\bf Cp}(I)}$ and $\xymatrix @C=5mm {\sigma ^{\geq}: I \ar@<-1.5pt>@{^{(}->}[r] & {\bf Cp}(I)}$, noted respectively $\Delta ^{\scriptscriptstyle (\leq )}_{{\bf Cp}(I)} (\Gamma ) = \bigl ( \delta ^{\scriptscriptstyle (\leq )}(\Gamma) _{{\bf L}} ^-, \delta ^{\scriptscriptstyle (\leq )}(\Gamma ) _{{\bf L}} ^+ \bigr ) _{{\bf L} \in {\bf Cp}(I)}$ and $\Delta ^{\scriptscriptstyle (\geq )} _{{\bf Cp}(I)} (\Gamma ) = \bigl ( \delta ^{\scriptscriptstyle (\geq )} (\Gamma) _{{\bf L}} ^-, \delta ^{\scriptscriptstyle (\geq )} (\Gamma ) _{{\bf L}} ^+ \bigr ) _{{\bf L} \in {\bf Cp}(I)}$. 

Moreover for any cut ${\bf L} = ({\bf L}_- , {\bf L}_+)$ of $I$, these structures satisfy the following equalities: 
$$\delta ^{\scriptscriptstyle (\leq )}(\Gamma ) _{{\bf L}}^- =  \bigcup _{i \in {\bf L} _+} \delta (\Gamma ) _{i}^+ \quad\hbox{and}\quad   
\delta ^{\scriptscriptstyle (\leq )}(\Gamma ) _{{\bf L}}^+  = \bigcap _{i \in {\bf L}_-} \delta (\Gamma ) _{i} ^+ \ ;$$ 
$$\delta ^{\scriptscriptstyle (\geq )}(\Gamma ) _{{\bf L}}^- =  \bigcup _{i \in {\bf L} _+} \delta (\Gamma ) _{i}^- \quad\hbox{and}\quad   
\delta ^{\scriptscriptstyle (\geq )}(\Gamma ) _{{\bf L}}^+  = \bigcap _{i \in {\bf L}_-} \delta (\Gamma ) _{i} ^- \ .$$ 

\end{corollary} 

\begin{preuve} 
We consider the first embedding $\xymatrix @R=-1mm @C=8mm{\sigma ^{\leq} : I \ar@<-1.5pt>@{^{(}->}[r] & {\bf Cp}(I)}$, and we identify $I$ with its image in ${\bf Cp}(I)$, then for ${\bf L} = ({\bf L} _- , {\bf L} _+)$ in ${\bf Cp}(I)$ and for $i \in I$ we remember the following relations: 
\begin{enumerate} 
\item\emph{$i < {\bf L}$ if and only if $i \in \bigcup _{j \in {\bf L} _-}I_{<j}$;}  
\item\emph{$i = {\bf L}$ if and only if ${\bf L} _- = I_{\leq i}$;}  
\item\emph{$i > {\bf L}$ if and only if $i \in {\bf L} _+$.}  
\end{enumerate} 
For any ${\bf L} \in {\bf Cp}(I)$ we define the submodules $\delta ^{\scriptscriptstyle (\leq )}(\Gamma ) _j ^-$ and $\delta ^{\scriptscriptstyle (\leq )}(\Gamma ) _j ^+$ by  
$$\delta ^{\scriptscriptstyle (\leq )}(\Gamma ) _{{\bf L}}^- =  \bigcup _{i > {\bf L}} \delta (\Gamma ) _{i}^+ =  \bigcup _{i \in {\bf L} _+} \delta (\Gamma ) _{i}^+ \quad\hbox{and}\quad   
\delta ^{\scriptscriptstyle (\leq )}(\Gamma ) _{{\bf L}}^+  = \bigcap _{i < {\bf L}} \delta (\Gamma ) _{i} ^- = \bigcap _{j \in {\bf L} _-} \bigcap _{i < j} \delta (\Gamma ) _{i} ^- = \bigcap _{j \in {\bf L}_-} \delta (\Gamma ) _{j} ^+ \ .$$ 

In the same way if we consider the second embedding $\xymatrix @R=-1mm @C=8mm{\sigma ^{\geq} : I \ar@<-1.5pt>@{^{(}->}[r] & {\bf Cp}(I)}$, and we identify $I$ with its image in ${\bf Cp}(I)$, then for ${\bf L} = ({\bf L} _- , {\bf L} _+)$ in ${\bf Cp}(I)$ and for $i \in I$ we remember the following relations: 
\begin{enumerate} 
\item\emph{$i < {\bf L}$ if and only if $i \in {\bf L} _-$;}  
\item\emph{$i = {\bf L}$ if and only if ${\bf L} _+ = I_{\geq i}$;}  
\item\emph{$i > {\bf L}$ if and only if $i \in \bigcup _{j \in {\bf L} _+}I_{>j}$.}  
\end{enumerate} 
For any ${\bf L} \in {\bf Cp}(I)$ we define the submodules $\delta ^{\scriptscriptstyle (\geq )}(\Gamma ) _j ^-$ and $\delta ^{\scriptscriptstyle (\geq )}(\Gamma ) _j ^+$ by  
$$\delta ^{\scriptscriptstyle (\geq )}(\Gamma ) _{{\bf L}}^- =  \bigcup _{i > {\bf L}} \delta (\Gamma ) _{i}^+  = \bigcup _{j \in {\bf L} _+}  \bigcup _{i > j} \delta (\Gamma ) _{i}^+ =  \bigcup _{j \in {\bf L} _+} \delta (\Gamma ) _{j}^- \quad\hbox{and}\quad   
\delta ^{\scriptscriptstyle (\geq )}(\Gamma ) _{{\bf L}}^+  = \bigcap _{i < {\bf L}} \delta (\Gamma ) _{i} ^-  = \bigcap _{i \in {\bf L}_-} \delta (\Gamma ) _{i} ^- \ .$$ 

\hfill$\Box$ 
\end{preuve} 

\vskip .2cm

\begin{remark}\label{rmq:totalement-degeneree} 
For an $R$-module $\Gamma$ endowed with an $I$-structure $\Delta _I (\Gamma ) = \bigl ( \delta (\Gamma) _i ^-, \delta (\Gamma ) _i ^+ \bigr ) _{i \in I}$, we can define another family $\Delta _{{\bf Cp}(I)} (\Gamma ) = \bigl ( \delta (\Gamma) _{{\bf L}} ^-, \delta (\Gamma ) _{{\bf L}} ^+ \bigr ) _{{\bf L} \in {{\bf Cp}(I)}}$ of $R$-submodules of $\Gamma$ by 
$$\delta (\Gamma ) _{{\bf L}}^- =  \bigcup _{j \in {\bf L} _+} \delta (\Gamma ) _j^+ \quad\hbox{and}\quad  \delta (\Gamma ) _{{\bf L}}^+ = \bigcap _{i \in {\bf L} _-} \delta (\Gamma ) _i ^- \ ,$$ 
and this family  $\Delta _{{\bf Cp}(I)} (\Gamma )$ is an ${\bf Cp}(I)$-prestructure on $\Gamma$. 
 
If we consider the two cuts ${\bf L} ^{\geq i} = ( I_{<i}, I_{\geq i})$ and ${\bf L} ^{> i} = ( I_{\leq i}, I_{> i})$ associated with an element $i$ of $I$ we have the following eqqualities: 
$$\delta (\Gamma ) ^- _{{\bf L} ^{>i}} = \delta (\Gamma ) ^+ _{{\bf L} ^{>i}} = \delta (\Gamma )_i ^- 
\quad\hbox{and}\quad 
\delta (\Gamma ) _{{\bf L} ^{\geq i}} ^- =  \delta (\Gamma ) _{{\bf L} ^{\geq i}} ^+ = \delta (\Gamma )_i ^+ \ .$$

Moreover if the set $I$ is well-ordered, by remark \ref{rmq:bien-ordonne} we know that any cut ${\bf L}$ of $I$ is of the form ${\bf L} = {\bf L} ^{\geq i}$ for some $i \in I$. Then for any $I$-module $\Gamma$ the ${\bf Cp}(I)$-prestructure induced on $\Gamma$ is totally degenerated. 
\end{remark} 

\vskip .2cm

\subsection{Immediate morphism and Hahn-product} 

\begin{definition} 
Let $\Gamma$ and $\Theta$ be two $I$-modules, an \emph{$I$-morphism} $\phi$ between $\Gamma$ and $\Theta$ is a morphism of $R$-modules $\xymatrix{\phi : \Gamma \ar[r] & \Theta}$ satisfying 
$\phi (\delta (\Gamma )_i ^-) \subset \delta (\Theta )_i ^-$ and $\phi (\delta (\Gamma )_i ^+) \subset \delta (\Theta )_i ^+$ for any $i$ in $I$. 
\end{definition} 

Let $\xymatrix{\phi : \Gamma \ar[r] & \Theta}$ be an $I$-morphism between two $I$-modules, then for any $i$ in $I$ the morphism $\phi$ induces a morphism of modules: 
$\xymatrix{\varepsilon (\phi )_i : \varepsilon (\Gamma )_i \ar[r] & \varepsilon (\Theta )_i}$.  

\vskip .2cm

Let $\Gamma$ be a $R$-module and $\Gamma '$ be a submodule of $\Gamma$, then an $I$-structure on $\Gamma$ induces an $I$-structure on $\Gamma '$ with the equalities $\delta (\Gamma ') _i ^- =  \delta (\Gamma ) _i ^- \cap \Gamma '$ and $\delta (\Gamma ') _i ^+ = \delta (\Gamma ) _i ^+ \cap \Gamma '$. 
We have an injective $I$-morphism $\xymatrix{\phi : \Gamma ' \ar[r] & \Gamma}$ and the induced morphisms $\xymatrix{\varepsilon (\phi )_i : \varepsilon (\Gamma ')_i \ar[r] & \varepsilon (\Gamma )_i}$ are injective.  

\begin{definition} 
An injective $I$-morphism $\xymatrix{\phi : \Gamma \ar@<-1.5pt>@{^{(}->}[r] & \Theta}$ is \emph{immediate} if for any $i$ in $I$ we have the equalities 
$\phi (\delta (\Gamma )_i ^-) = \delta (\Theta )_i ^- \cap \phi (\Gamma )$ and $\phi (\delta (\Gamma )_i ^+) = \delta (\Theta )_i ^+ \cap \phi (\Gamma )$, and if the induced morphisms 
$\xymatrix{\varepsilon (\phi )_i : \varepsilon (\Gamma )_i \ar[r] & \varepsilon (\Theta )_i}$ are isomorphisms. 
\end{definition} 

If the $I$-morphism $\xymatrix{\phi : \Gamma \ar@<-1.5pt>@{^{(}->}[r] & \Theta}$ is immediate, the skeletons of the $I$-modules $\Gamma$ and $\Theta$ are the same. 

\vskip .2cm 

Let $\xymatrix{\phi : \Gamma \ar[r] & \Theta}$ be a morphism of $I$-modules, then for any $\xi$ in $\Gamma$ we have the inequality $\iota _{\Delta} (\xi ) \leq \iota _{\Delta}(\phi (\xi ))$, 
if moreover we assume that the morphism is immediate we have the equality $\iota _{\Delta} (\xi ) = \iota _{\Delta}(\phi (\xi ))$ and the image by the morphism $\varepsilon (\phi ) _{\iota _{\Delta} (\xi )}$ of the initial part of $\xi$ is the initial part of $\phi (\xi )$.

\vskip .2cm

\begin{definition} 
Let $(\Theta _i)_{i \in I}$ be a family of $R$-modules indexed by a totally ordered set $I$, the \emph{Hahn product} $\prod _{i \in I} ^{(H)} \Theta _i$ is the submodule of the $R$-module $\prod _{i \in I} \Theta _i$ made up of the elements $\underline x = (x_i)$ whose support $Supp(\underline x) = \{ i \in I \ | \ x_i \not= 0 \}$ is a well-ordered subset of $I$. 
\end{definition} 

We have the inclusions: $\bigoplus _{i \in I}  \Theta _i \subset \prod _{i \in I} ^{(H)} \Theta _i \subset \prod _{i \in I} \Theta _i $, with the equality $\prod _{i \in I} ^{(H)} \Theta _i = \prod _{i \in I} \Theta _i $ if the set $I$ is well-ordered and the equalities $\bigoplus _{i \in I}  \Theta _i = \prod _{i \in I} ^{(H)} \Theta _i = \prod _{i \in I} \Theta _i $ if $I$ is finite. 

\vskip .2cm 

Let $\Theta =\prod _{i \in I} ^{(H)} \Theta _i$ be a $R$-module obtained as the Hahn product of the family $(\Theta _i)_{i \in I}$ indexed by the totally ordered set $I$. 
For any cut ${\bf L} = ({\bf L} _- , {\bf L} _+)$ of $I$, we can define a submodule $\Theta _{{\bf L}}$ of $\Theta$ as follows:  
$$\Theta _{{\bf L}} = \{ \underline x = (x_i) \in \Theta \ | \ x_i=0 \ \hbox{for any} \ i \in {\bf L} _- \} \ .$$

For any $\underline x = (x_i)$ belonging to $\Theta$ we note $\iota _{\underline x}$ the smallest element of $Supp(\underline x)$, i.e. the element of $I$ defined by 
$$\forall i , \ i< \iota _{\underline x} , \, x_i =0 \quad\hbox{and}\quad x_{\iota _{\underline x}} \not= 0 \ ,$$ 
then we have 
$$\Theta _{{\bf L}} = \{ \underline x = (x_i) \in \Theta \ | \ \iota _{\underline x} \in {\bf L} _+ \} \ .$$

\vskip .2cm 

\begin{remark} 
The submodule $\Theta _{{\bf L}}$ is isomorphic to the module $\prod _{i \in {\bf L} _+} ^{(H)} \Theta _i$ and the quotient module $\Theta / \Theta _{{\bf L}}$ is isomorphic to the module $\prod _{i \in {\bf L} _-} ^{(H)} \Theta _i$. 

If ${\bf L}$ and ${\bf L} '$ are two cuts of $I$ satisfying ${\bf L} \leq {\bf L} '$, we have the inclusion $\Theta _{{\bf L} '} \subset \Theta _{{\bf L}}$. 
\end{remark}

\vskip .2cm 

\begin{proposition} 
Let $\Theta =\prod _{i \in I} ^{(H)} \Theta _i$ be a $R$-module obtained as the Hahn product of the family $(\Theta _i)_{i \in I}$, then the family $\Delta _I(\Theta ) = \bigl ( \delta (\Theta) _i ^-, \delta (\Theta ) _i ^+ \bigr ) _{i \in I}$ defined by  
$\delta (\Theta )_i^- = \Theta _{{\bf L} ^{>i}} = \prod _{j > i} ^{(H)} \Theta _j$ and $\delta (\Theta )_i^+ = \Theta _{{\bf L} ^{\geq i}} = \prod _{j \geq i} ^{(H)} \Theta _j$ 
defines an $I$-structure on $\Theta$. 

Moreover for any element $\underline x = (x_i)$ of $\Theta$ we have the equality $\iota _{\Delta}(\underline x) = \iota _{\underline x}$, and its initial part ${\bf in}_{\Theta}(x)$ is the element $x_{\iota _{\underline x}}$.  

The skeleton of $\Theta =\prod _{i \in I} ^{(H)} \Theta _i$ is the pair $\bigl ( I , (\Theta _i) _{i \in I} \bigr )$
\end{proposition}

\begin{preuve} 
For any $i$ belonging to $I$ the submodules $\Theta _{{\bf L} ^{>i}}$ and $\Theta _{{\bf L} ^{\geq i}}$ associated respectively with the cuts ${\bf L} ^{>i}$ and ${\bf L} ^{\geq i}$ of the set $I$ are defined by   
$$\Theta _{{\bf L} ^{>i}} = \{ \underline x \in \Theta \ |  \ \iota _{\underline x} > i \} \quad\hbox{and}\quad \Theta _{{\bf L} ^{\geq i}} = \{ \underline x \in \Theta \ | \ \iota _{\underline x} \geq  i \} \ ,$$
and satisfy $\Theta _{{\bf L} ^{>i}} \subset \Theta _{{\bf L} ^{\geq i}}$ and $\Theta _{{\bf L} ^{\geq i}} / \Theta _{{\bf L} ^{> i}} \simeq \Theta _i$. 

For $i<j$ we have the inclusion ${\bf L} ^{>i} \subset {\bf L} ^{\geq j}$, hence $\Theta _{{\bf L} ^{\geq j}} \subset \Theta _{{\bf L} ^{>i}}$. 

For any $\underline x = (x_i)$ belonging to $\prod _{i \in I} ^{(H)} \Theta _i$ we have $\underline x  \in \Theta _{{\bf L} ^{\geq \iota _{\underline x}}}$ and $\underline x  \notin \Theta _{{\bf L} ^{> \iota _{\underline x}}}$.

\hfill$\Box$ 
\end{preuve}

Let $\Theta$ be the $R$-module obtained as the Hahn product of a family $\bigl ( \Theta _i \bigr ) _{i \in I}$, we always assume it is endowed with the previously defined  $I$-structure $\Delta _I$. 

For any cut ${\bf L}$ of $I$ we have the following equivalences:  
$$\begin{array}{l}
i \in {\bf L} _- \ \Longleftrightarrow \  {\bf L}  ^{>i} \leq {\bf L}  \ \Longleftrightarrow \ \Theta _{{\bf L}} \subset \Theta _{{\bf L} ^{>i}} \\ 
i \in {\bf L} _+ \ \Longleftrightarrow \  {\bf L} \leq {\bf L} ^{\geq i} \ \Longleftrightarrow \ \Theta _{{\bf L} ^{\geq i }} \subset \Theta _{{\bf L}} \ ,
\end{array}$$ 
we deduce from remark \ref{rmq:coupure} \emph{(2)} the following equalities:   
$$\bigcup _{i \in {\bf L} _+} \Theta _{{\bf L} ^{\geq i }} = \Theta _{{\bf L}} = \bigcap _{i \in {\bf L} _-} \Theta _{{\bf L} ^{> i}} \ .$$ 

\vskip .2cm

Let $\Gamma$ be a $R$-module and let $\Delta _I (\Gamma ) = \bigl ( \delta (\Gamma) _i ^-, \delta (\Gamma ) _i ^+ \bigr ) _{i \in I}$ be an $I$-structure on $\Gamma$. 
Hence we can define the module $\varepsilon (\Gamma )$ obtained as the Hahn product of the skeleton of  $\Gamma$, i.e. defined by 
$\varepsilon (\Gamma ) = \prod _{i \in I} ^{(H)} \varepsilon (\Gamma ) _i$. 
By construction the skeletons of the $I$-modules $\Gamma$ and $\varepsilon (\Gamma )$ are identical and equal to $\bigl ( I , (\varepsilon (\Gamma )_i) _{i \in I} \bigr )$.  

We want to find conditions on the $I$-structure of $\Gamma$ allowing us to construct an immediate $I$-morphism $\xymatrix{\underline{v} : \Gamma \ar[r] & \varepsilon (\Gamma )}$. 

If such a morphism exists we remark that for any $\xi$ belonging to $\Gamma$ its image $\underline{v}(\xi )$ in $\varepsilon (\Gamma )$ is of the form $\underline{v}(\xi ) = ( v_i(\xi ))_{i \in I}$, with $v_i(\xi )$ belonging to the module $\varepsilon (\Gamma )_i$. 
As the le morphism $\underline v$ is an $I$-morphism we have the inequality $\iota _{\Delta} (\xi ) \leq \iota _{\Delta}(\underline v(\xi )) = \iota _{\underline v(\xi )}$, we deduce from this that we have $v_j(\xi ) =0$ for $j < \iota _{\Delta} (\xi )$, and as the morphism $\underline v$ is immediate we have the equality $\iota _{\Delta} (\xi ) = \iota _{\underline v(\xi )}$. 
We say that the morphism $\xymatrix{\underline{v} : \Gamma \ar[r] & \varepsilon (\Gamma )}$ is a \emph{natural immediate $I$-morphism} if moreover we have the equality $v_{\iota _{\Delta} (\xi )}(\xi ) = v_{\iota _{\underline v(\xi )}}(\xi ) = {\bf in}_{\Delta}(\xi )$.

\vskip .2cm

\begin{definition} 
Let $\Gamma$ be an $I$-module, a \emph{$I$-decomposition} of $\Gamma$ is the data of a family $\Psi (\Gamma )= \bigl ( \psi (\Gamma ) _i \bigr ) _{i \in I}$ of submodules of $\Gamma$ satisfying: 
\begin{enumerate} 
\item $\delta (\Gamma )_i ^- = \delta (\Gamma )_i^+ \cap \psi (\Gamma ) _i$ for any $i$ in $I$; 
\item $\Gamma = \delta (\Gamma )_i^+ + \psi (\Gamma )_i$ for any $i$ in $I$; 
\item for any $\xi$ in $\Gamma$ the subset $\{ i \in I \, | \, \xi \notin \psi (\Gamma ) _i \}$ is a well-ordered set. 
\end{enumerate} 
\end{definition}

\begin{proposition} \label{prop:equivalence-decomposition-morphisme} 
Let $\Gamma$ be an $I$-module, then there exists a bijection between the set of the natural immediate  $I$-morphisms $\xymatrix{\underline{v} : \Gamma \ar[r] & \varepsilon (\Gamma )}$ and the set of the  $I$-decompositions of $\Gamma$. 
\end{proposition} 

\begin{preuve} 
%
The data of a morphism $\xymatrix{\underline{v} : \Gamma \ar[r] & \varepsilon (\Gamma ) = \prod _{i \in I} ^{(H)} \varepsilon (\Gamma ) _i}$ is equivalent to the data of a family $\bigl ( v_i \bigr )_{i \in I}$ of morphisms of modules $\xymatrix{v_i : \Gamma \ar[r] & \varepsilon (\Gamma )_i} = \delta (\Gamma )_i^+ / \delta (\Gamma )_i^-$ such that for any $\xi$ belonging to $\Gamma$ the subset $\{ i \in I \, | \, v_i(\xi )\not= 0 \} = \{ i \in I \, | \, \xi \notin Ker(v_i) \}$ is a well-ordered set.

For the morphism $\underline v$ to be a natural immediate $I$-morphism we must further assume that for any $i \in I$ the restriction of the morphism $\xymatrix{v_i : \Gamma \ar[r] & \varepsilon (\Gamma )_i}$ to the submodule $\delta (\Gamma )_i^+$ is the canonical morphsim $\xymatrix{\delta (\Gamma )_i^+ \ar@{->>}[r] & \delta (\Gamma )_i^+ / \delta (\Gamma )_i^-}$. 

The proposition is then a consequence of the following lemma.  

\hfill$\Box$ 
\end{preuve} 

\vskip .2cm 

\begin{lemma} 
Let $\Gamma$ be a $R$-module and $\Delta ^-$ and $ \Delta ^+$ be two submodules of $\Gamma$ with $\Delta ^- \subset \Delta ^+$, then there exists a bijection $\psi$ between the two following sets.
\begin{enumerate}
\item The set ${\bf Mor }_{\Gamma}$ of morphisms of modules $\xymatrix{v:\Gamma \ar[r] & \Delta ^+ / \Delta ^-}$ such that the restriction of $v$ to the submodule $\Delta ^+$ is the natural morphism $\xymatrix{\Delta ^+ \ar@{->>}[r] & \Delta ^+ / \Delta ^-}$. 
\item The set ${\bf Mod}_{\Gamma}$ of submodule $\Psi$ of $\Gamma$ satisfying: 
\begin{enumerate} 
\item $\Delta ^- = \Delta^+ \cap \Psi$; 
\item $\Gamma = \Delta^+ + \Psi$.  
\end{enumerate}  
\end{enumerate} 

This bijection is define by $\psi (v) = \Psi = Ker (v)$. 
\end{lemma} 

\begin{preuve} 
Let us assume that there exists a morphism $\xymatrix{v:\Gamma \ar[r] & \Delta ^+ / \Delta ^-}$ that satisfies the properties stated above, then the morphism $v$ is surjective and if we denote $\Psi$ its kernel the inclusion of $\Delta ^+$ in $\Gamma$ induces an isomorphism $\xymatrix{\bar u:\Delta ^+ / \Delta ^- \ar[r] & \Gamma / \Psi}$  such that the inverse morphism ${\bar u}^{-1}$ is the morphism $\xymatrix{\bar v:\Gamma / \Psi \ar[r] & \Delta ^+ / \Delta ^-}$ induced by $v$.   
We deduce the properties (a) and (b) of the submodule $\Psi$ of $\Gamma$. 

Conversely if we have a submodule $\Psi$ of $\Gamma$ which satisfies the properties (a) and (b) we have an isomorphism $\xymatrix{\bar u:\Delta ^+ / \Delta ^- \ar[r] & \Gamma / \Psi}$ and the morphism $\xymatrix{v:\Gamma \ar[r] & \Delta ^+ / \Delta ^-}$ obtained as the composite morphism ${\bar u}^{-1} \circ w$, where $w$ is the natural morphism of $\Gamma$ to $\Gamma / \Psi$, satisfies the properties we're looking for.  

\hfill$\Box$ 
\end{preuve} 

\vskip .2cm

Let $\Gamma$ be an $I$-module, let $\Gamma '$ be a submodule of $\Gamma$ provided with the induced $I$-structure and let $\Psi (\Gamma ') = \bigl ( \psi (\Gamma ') _i \bigr ) _{i \in I}$ be a $I$-decomposition of $\Gamma '$, then an \emph{extension} of the $I$-decomposition $\Psi (\Gamma ')$ to $\Gamma$ is a $I$-decomposition $\Psi (\Gamma ) = \bigl ( \psi (\Gamma ) _i \bigr ) _{i \in I}$ of $\Gamma$ such that for any $i$ in $I$ we have $\psi (\Gamma ')_i = \psi (\Gamma )_i \cap \Gamma '$. 

\vskip .2cm

\begin{theorem}\label{th:$I$-decomposition} 
Let $\Gamma$ be an $I$-vector space on $\F$ and let $\Gamma '$ be a vector subspace of $\Gamma$ provided with the induced $I$-structure, then any $I$-decomposition of $\Gamma '$ admits an extension to $\Gamma$. 
\end{theorem}

\begin{preuve} 
The result is proved by transfinite recursion. We recall that if  $\Gamma '$ is a vector subspace of $\Gamma$ there exists a well-ordered family of vector subspaces $\Pi (\alpha )$ such that  
\begin{enumerate} 
\item \emph{$\Pi (0) = \Gamma '$;}  
\item \emph{$\Pi (\alpha +1) = \Pi (\alpha ) \oplus \F \zeta$, for an element $\zeta$ of $\Gamma$;}  
\item \emph{$\Pi (\beta ) = \bigcup _{\alpha < \beta} \Pi (\alpha )$ if $\beta$ is a limit ordinal;}  
\item \emph{$\Pi (\gamma ) = \Gamma$ for an ordinal $\gamma $.} 
\end{enumerate}

We will show by recursion that any $I$-decomposition $\Psi (\Gamma ') = \bigl ( \psi (\Gamma ') _i \bigr ) _{i \in I}$ of $\Gamma '$, admits an extension $\Psi (\Gamma )$ to $\Gamma$. 
Let us assume that for any $\alpha < \beta$ we have a $I$-decomposition $\Psi (\Pi (\alpha ))$ of the vector space $\Pi (\alpha )$ and we want to extend it to a $I$-decomposition of $\Pi (\beta )$. 

By hypothesis for any $\alpha$ the $I$-structure on $\Pi (\alpha )$ is the one induced by the $I$-structure on $\Gamma$, then we have the following equalities
$\delta (\Pi (\alpha )) _i ^- =  \delta (\Gamma ) _i ^- \cap \Pi (\alpha )$ and $\delta (\Pi (\alpha )) _i ^+ = \delta (\Gamma ) _i ^+ \cap \Pi (\alpha )$. 

\vskip .2cm 

\emph{Cas} I. $\beta$ is a limit ordinal.  

The subspace $\Pi (\beta )$ is of the form $\Pi (\beta ) = \bigcup _{\alpha < \beta} \Pi (\alpha )$ and for any $i$ in $I$ we have also the equalities $\delta (\Pi (\beta ))_i^- = \bigcup _{\alpha < \beta} \delta (\Pi (\alpha ))_i^-$ and $\delta (\Pi (\beta ))_i^+ = \bigcup _{\alpha < \beta} \delta (\Pi (\alpha ))_i^+$. 

We then define for any $i$ in $I$: 
 $$\hbox{$\psi (\Pi (\beta ))_i = \bigcup _{\alpha < \beta} \psi (\Pi (\alpha )) _i$  ,}$$   
and it is easy to check that we define in this way a $I$-decomposition $\Psi (\Pi (\beta ))$ of $\Pi (\beta )$. 

\vskip .2cm 

\emph{Cas} II.  $\beta$ is equal to $\alpha +1$.   

The subspace $\Pi (\beta )$ is of the form $\Pi (\beta ) = \Pi (\alpha ) \oplus \F \zeta$ and for any $i$ in $I$ we have the following equalities 
$$\xymatrix @R=1mm @C=1mm{
\hbox{$\delta (\Pi (\beta ))_i^- = \delta (\Pi (\alpha ))_i^- \oplus \F \zeta$ if $\zeta$ belongs to $\delta (\Gamma )_i^-$} & \hbox{and $\delta (\Pi (\beta ))_i^- = \delta (\Pi (\alpha ))_i^-$ otherwise; } \\ 
\hbox{$\delta (\Pi (\beta ))_i^+ = \delta (\Pi (\alpha ))_i^+ \oplus \F \zeta$ if $\zeta$ belongs to $\delta (\Gamma )_i^+$} & \hbox{and $\delta (\Pi (\beta ))_i^+ = \delta (\Pi (\alpha ))_i^+$ otherwise.} 
}$$
 
\vskip .2cm

\begin{lemma} 
Let $\Gamma$ be a vector space and $\Delta ^- \subset \Delta ^+$ and $\Psi (\alpha )\subset \Pi (\alpha )$ four subspaces of $\Gamma$ satisfying 
$\Delta ^+ \cap \Psi (\alpha )= \Delta ^- \cap \Pi (\alpha )$ and $(\Delta ^+ \cap \Pi (\alpha )) + \Psi (\alpha )= \Pi (\alpha )$.

Then for any element $\zeta$ of $\Gamma$ with $\zeta \notin \Pi$, if we claim $\Pi (\beta ) = \Pi (\alpha ) \oplus \F \zeta$, there exists $\Psi (y)\subset \Pi (\beta )$ satisfying 
$\Delta ^+ \cap \Psi (\beta )= \Delta ^- \cap \Pi (\beta )$ and $(\Delta ^+ \cap \Pi (\beta )) + \Psi (\beta )= \Pi (\beta )$.  
\end{lemma} 

\begin{preuvele}
If $\zeta$ belongs to $\Delta ^-$ or if $\zeta$ doesn't belong to $\Delta ^+$ we claim $\Psi (\beta ) = \Psi (\alpha ) + \F \zeta$ and if $\zeta$ belongs to $\Delta ^+ \setminus \Delta ^-$ we claim $\Psi (\beta ) = \Psi (\alpha )$. 

\hfill$\Box$
\end{preuvele}

We deduce from the previous lemma that if we set
\begin{itemize}
\item\emph{$\psi (\Pi (\beta ))_i = \psi (\Pi (\alpha )) _i\oplus \F \zeta$ for any $i$ such that  $\zeta \in \delta (\Gamma )_i^-$ or $\zeta \notin \delta (\Gamma )_i^+$,  i.e. for $i \not= \iota _{\Delta}(\zeta)$;}    
\item $\psi (\Pi (\beta ))_i = \psi (\Pi (\alpha ))_i$ for any $i$ such that  $\zeta \in \delta (\Gamma )_i^+ \setminus \delta (\Gamma )_i^-$, i.e. for $i= \iota _{\Delta}(\zeta)$. 
\end{itemize} 
we find a family of vector subspaces of $\Pi (\beta )$ which, for any $i$ in $I$, satisfies the following properties  
\begin{enumerate} 
\item $\delta (\Pi (\beta ))_i ^- = \delta (\Pi (\beta ))_i^+ \cap \psi (\Pi (\beta )) _i$; 
\item $\Pi (\beta ) = \delta (\Pi (\beta ))_i^+ + \psi (\Pi (\beta ))_i$. 
\end{enumerate} 

It remains to be verified that for any $\xi$ belonging to $\Pi (\beta )$ the subset $\{ i \in I \, | \, \xi \notin \psi (\Pi (\beta )) _i \}$ is well-ordered. 

We can write $\xi = \xi ' + f \zeta$ with $\xi ' \in \Pi (\alpha )$ and $f \in \F$, and we deduce from the above that the subset $\{ i \in I \, | \, \xi \notin \psi (\Pi (\beta )) _i \}$ is included in $\{ i \in I \, | \, \xi ' \notin \psi (\Pi (\alpha )) _i \} \cup \{ \iota _{\Delta}(\zeta ) \}$, hence is well-ordered by induction hypothesis.  

\hfill$\Box$ 
\end{preuve} 

\vskip .2cm 

\begin{corollary}\label{cor:theoreme-de-Hahn} 
Let $\Gamma$ be a $\F$-vector space endowed with an $I$-structure, then there exists a natural immediate $I$-morphism $\xymatrix{\underline v : \Gamma \ar[r] & \varepsilon (\Gamma ) = \prod _{i \in I} ^{(H)} \varepsilon (\Gamma )_i}$. 
\end{corollary} 

\begin{preuve} 
It is enough  to apply the previous theorem to the trivial structure on $\Gamma ' =(0)$. 

\hfill$\Box$ 
\end{preuve}

\vskip .2cm 

\begin{remark}\label{rmq:filtration-graduation} 
We may consider that an $I$-structure on an $R$-module $\Gamma$ is an analogue of a filtration of the module by a totally ordered set $I$, and the module $\varepsilon (\Gamma )$ is the analogue of the associated graduated module. 
\end{remark}  

\vskip .2cm

		      \section{Totally ordered abelian groups}
%

\subsection{Convex subgroups} 

We suppose now that $\Gamma$ is an \emph{ordered group}, i.e. an abelian group endowed with a total ordering which is compatible with the law group.  We denote this group law by the addition. 
An interval $\Xi$ of $\Gamma$ is \emph{symmetric} if the opposite $- \xi$ of any $\xi$ belonging to $\Xi$ belongs also to $\Xi$, we consider especially the three symmetric intervals $\Xi = \emptyset$, $\Xi = \{ 0 \}$ and $\Xi = \Gamma$. 
 
\vskip .2cm 

\begin{remark}\label{rmq:inclusion-intervalles-symetriques}
The set ${\bf Sym}(\Gamma )$ of symmetric intervals of an ordered group $\Gamma$ is totally ordered by inclusion. 
\end{remark} 

\vskip .2cm

A cut ${\Lambda} = \bigl ({\Lambda} _- , {\Lambda} _+ \bigr )$ of a totally ordered group $\Gamma$ is called \emph{positive} (resp. \emph{non negative}) if ${\Lambda} _+$ is included in the subset $\Gamma _{> 0}$ of the positive elements (resp. in the subset $\Gamma _{\geq  0}$ of the non negative elements) of $\Gamma$, we denote by ${\bf Cp}^{++}(\Gamma )$ (resp. by ${\bf Cp}^{+}(\Gamma )$)  the set of positive (resp. non negative) cuts of $\Gamma$: 
$${\Lambda} \in {\bf Cp}^{++}(\Gamma )  \  \Longleftrightarrow \  {\Lambda} ^{> 0} = {\Lambda} ^{\leq 0} \leq {\Lambda} \quad\hbox{and}\quad {\Lambda} \in {\bf Cp}^{+}(\Gamma )  \  \Longleftrightarrow \  {\Lambda} ^{\geq 0} \leq {\Lambda} \ .$$ 
We can identify the set ${\bf Cp}^{+}(\Gamma )$ of non-negative cuts on $\Gamma$ with the set ${\bf Cp}(\Gamma _{\geq 0})$ of cuts on $\Gamma _{\geq  0}$. 

\vskip .2cm 

With any non negative cut ${\Lambda}$ we can associate a symmetric interval $\Xi _{{\Lambda}}$ defined as the complementary set of ${\Lambda} _+ \cup - {\Lambda} _+$, which is also equal to the subset ${\Lambda} _- \cap - {\Lambda} _-$. 
This interval is non empty if the cut ${\Lambda}$ is positive. 

A positive cut ${\Lambda}$ is \emph{irreducible} if the subset ${\Lambda} _- \cap \Gamma _{\geq 0}$ is stable by the addition, i.e. if for any $\zeta$ and $\xi$ in $\Gamma$ we have 
$$\zeta \in {\Lambda} _- \cap \Gamma _{\geq 0}  \  \hbox{and} \ \xi \in {\Lambda} _- \cap \Gamma _{\geq 0}  \ \Longrightarrow \  \zeta + \xi \in {\Lambda} _- \cap \Gamma _{\geq 0} \ .$$  
This is equivalent to the following condition 
$$\zeta \in {\Lambda} _- \cap \Gamma _{\geq 0} \ \Longrightarrow \  2 \zeta = \zeta + \zeta \in {\Lambda} _- \cap \Gamma _{\geq 0} \ .$$

\vskip .2cm 

We remind the following definition. 

\begin{definition} 
A \emph{convex} subgroup of $\Gamma$ is a subgroup $\Delta$ which is also a symmetric interval.
\end{definition} 

\vskip .2cm 

\begin{lemma}\label{le:coupure-irreductible} 
The positive cut ${\Lambda} = \bigl ({\Lambda} _- , {\Lambda} _+ \bigr )$ is irreducible if and only if the interval $\Xi _{{\Lambda}}$ associated with ${\Lambda}$ is a convex subgroup of $\Gamma$. 
\end{lemma}

\begin{preuve} 
As  the interval $\Xi _{{\Lambda}}$ is symmetric it is enough to check that it is stable by the addition, which is an immediate consequence of the definition of $\Xi _{{\Lambda}}$ as the subset ${\Lambda} _- \cap - {\Lambda} _-$. 

\hfill$\Box$ 
\end{preuve} 

\vskip .2cm 

If $\Delta$ is a convex subgroup of $\Gamma$ the quotient group $\Gamma / \Delta$ can be endowed with a structure of an ordered group such that the natural morphism $\xymatrix{w _{\Delta} :\Gamma \ar[r] & \Gamma / \Delta}$ is a morphisme of ordered groups.  
Then for $\xi$ and $\zeta$ in $\Gamma$ we have the relation: 
$$
w(\xi ) \leq w(\zeta ) \  \Longleftrightarrow \  \exists \ \delta \in \Delta \ \hbox{such that } \ \xi + \delta \leq \zeta \ .
$$
Conversely if $\xymatrix{u:\Gamma \ar[r] & \Gamma '}$ is a morphism of ordered groups the kernel $\Delta = Ker(u)$ is a convex subgroup of $\Gamma$. 

\vskip .2cm 

We note ${\bf Cv}(\Gamma )$ the set of proper convex subgroups of $\Gamma$ and $\widehat{\bf Cv}(\Gamma )$ the set ${\bf Cv}(\Gamma ) \cup\{\Gamma\}$ of convex subgroups of $\Gamma$. These two sets are totally ordered by inclusion. 

If $\Delta _1$ and $\Delta _2$ are two distinct non-trivial convex subgroups of $\Gamma$ with $\Delta _1 \subset \Delta _2$ we denote by $\xymatrix{w_{\Delta _1}: \Gamma \ar[r] & \Gamma / \Delta _1}$ and $\xymatrix{w_{\Delta _2}: \Gamma \ar[r] & \Gamma / \Delta _2}$ the induced morphisms of ordered groups, the subgroup $\Delta _2 / \Delta _1$ is a convex subgroup of $\Gamma / \Delta _1$ and we denote by $\xymatrix{w_{\Delta _2 / \Delta _1}: \Gamma / \Delta _1 \ar[r] & \Gamma / \Delta _2}$ the induced morphism.

\begin{remark}\label{rmq:sous-groupes-isoles-du-quotient}
(1) For any convex subgroup $\Delta$ of $\Gamma$ the morphism $\xymatrix{w _{\Delta} :\Gamma \ar[r] & \bar\Gamma = \Gamma / \Delta}$ induces a morphism of totally ordered sets $\xymatrix{{\bf w} _{\Delta} :{\bf Sym}(\Gamma ) \ar[r] & {\bf Sym}(\bar\Gamma )}$ between the sets of symmetric intervals, such that the image of the subset ${\bf Cv}(\Gamma )$ of proper convex subgroups of $\Gamma$ is the subset ${\bf Cv}(\bar\Gamma )$ of proper convex subgroups of $\bar\Gamma$.  

(2) We have a natural isomorphisms of ordered sets between the sets of convex subgroups of the group $\Delta _2 / \Delta _1$ and the subsets of convex subgroups of $\Gamma$ between $\Delta _1$ and $\Delta _2$:
$$\begin{array}{l} 
{\bf Cv}(\Delta _2 / \Delta _1 ) \simeq \{ \Delta \in {\bf Cv}(\Gamma ) \ | \  \Delta _1 \subset \Delta \subsetneq \Delta _2 \} \\ 
\widehat{\bf Cv}(\Delta _2 / \Delta _1 ) \simeq \{ \Delta \in {\bf Cv}(\Gamma ) \ | \  \Delta _1 \subset \Delta \subset \Delta _2 \} \ .
\end{array}$$
\end{remark}

The order-type of the set ${\bf Cv}(\Gamma )$ is called the \emph{rank} of the group $\Gamma$ and is denoted ${\rm rk}(\Gamma )$. 

\vskip .2cm

Let $(\Theta _i)_{i \in I}$ be a family of ordered groups, indexed by a totally ordered set $I$, we remind the definition of the Hahn product $\Theta = \prod _{i \in I} ^{(H)} \Theta _i$: this is the subgroup of the abelian group $\prod _{i \in I} \Theta _i$ of the elements $\underline x = (x_i)$ whose support $Supp(\underline x) = \{ i \in I \ | \ x_i \not= 0 \}$ is a well-ordered subset of $I$.   

For any element $\underline x = (x_i)$ of the Hahn product $\Theta$, with $\underline x \not= 0$, we denote $\iota _{\underline x}$ the smallest element of $Supp(\underline x)$, and ${\bf in}_{\Theta}(\underline x)$ the element $x_{\iota _{\underline x}}$ of the group $\Theta _{\iota _{\underline x}}$.

We recall that a subset $\Theta ^+$ of an abelian group $\Theta$ defines a structure of totally ordered group on $\Theta$ by $\Theta ^+ = \{ \xi \in \Theta \ | \ \xi > 0 \}$, i.e. $\Theta ^+$ is the subset of positive elements of $\Theta$, if and only if it satisfies the following properties: 
\begin{enumerate} 
\item\emph{$(\Theta ^+ , - \Theta ^+)$ is a partition of $\Theta \setminus \{ 0 \}$;} 
\item\emph{$\Theta ^+$ is stable by the addition.} 
\end{enumerate}

\begin{proposition} 
The Hahn product $\Theta = \prod _{i \in I} ^{(H)} \Theta _i$, equipped with the lexicographic order defined by the subset $\Theta ^+ = \{ \underline x \in \Theta \ | \ {\bf in}_{\Theta}(\underline x) >0 \}$, is an ordered group.  
\end{proposition} 

\begin{preuve} 
It is obious that the subset $\Theta ^+$ satisfies the previous properties. 

\hfill$\Box$ 
\end{preuve}

We can also define the lexisographic order on the Hahn product $\Theta =\prod _{i \in I} ^{(H)} \Theta _i$ in the following way:  
let $\underline x = (x_i)$ and $\underline y = (y_i)$ be two distinct elements of $\Theta$, then the set $\{ i \in I \ |  x_i \not= y_i \}$ is a non-empty subset of $Supp(\underline x) \cup Supp(\underline y)$, hence it has a smallest element  $i_0$ and we claim $\underline x < \underline y$ if $x_{i_0} < y_{i_0}$.

\vskip .2cm 

As above we define for any cut ${\bf L} = ({\bf L} _- , {\bf L} _+)$ of  $I$ the subgroup  $\Theta _{{\bf L}}$ of $\Theta$ by the following: 
$$\Theta _{{\bf L}} \  = \  \{ \underline x = (x_i) \in \Theta \ | \ x_i=0 \ \hbox{for any} \ i \in {\bf L} _- \} \ = \  \{ \underline x = (x_i) \in \Theta \ | \ \iota _{\underline x} \in {\bf L} _+ \} \ .$$
 
\vskip .2cm 

\begin{lemma} 
The subgroup $\Theta _{{\bf L}}$ is a convex subgroup of $\Theta$, it is isomorphic to the ordered group $\prod _{i \in {\bf L} _+} ^{(H)} \Theta _i$ and the quotient group $\Theta / \Theta _{{\bf L}}$ is isomorphic to the ordered group $\prod _{i \in {\bf L} _-} ^{(H)} \Theta _i$. 
\end{lemma} 

\begin{preuve} 
It is enough to show that if we have $\underline x$ and $\underline y$ satisfying $0 < \underline y < \underline x$ with $\underline x$ belonging to $\Theta _{{\bf L}}$, then $\underline y$ belongs to $\Theta _{{\bf L}}$. 

By hypothesis we have ${\bf in}_{\Theta}(\underline y) = y _{\iota _{\underline y}} > 0$, $y_i =0$ for $i < \iota _{\underline y}$, ${\bf in}_{\Theta}(\underline x) = x _{\iota _{\underline x}} >0$, $x_i =0$ for $i < \iota _{\underline x}$, $\iota _{\underline x} \in {\bf L} _+$ and there exists $\iota \in I$ with $y_i = x_i$ for $i < \iota$ and $y_{\iota} < x _{\iota}$, then we have $\iota \geq Inf(\iota _{\underline x}, \iota _{\underline y})$. 
Let us assume that we have the inequality $\iota _{\underline y} < \iota _{\underline x}$, then we should have $y_{\iota _{\underline y}} >0$ and $x_{\iota _{\underline y}}=0$, hence $\iota _{\underline y} > \iota$, which is impossible. 
Then we deduce that we have the inequality $\iota _{\underline y} \geq \iota _{\underline x}$, hence $\iota _{\underline y}$ is also in ${\bf L} _+$, and $\underline y$ belongs to the subgroup $\Theta _{{\bf L}}$. 

\hfill$\Box$ 
\end{preuve}

\vskip .2cm

For any $i$ belonging to $I$ the subgroups $\Theta _{{\bf L} ^{>i}}$ and $\Theta _{{\bf L} ^{\geq i}}$, respectively associated with the cuts ${\bf L} ^{>i}$ and ${\bf L} ^{\geq i}$ of the set $I$, are convex subgroups of $\Theta$, and we deduce from remark \ref{rmq:sous-groupes-isoles-du-quotient} an isomorphism of ordered sets 
$$\xymatrix{ {\bf F}_i : \widehat{\bf Cv}(\Theta _i ) \ar[r] & \{ \Delta \in {\bf Cv}(\Gamma ) \ | \ \Theta _{{\bf L}^{>i}} \subset \Delta \subset  \Theta _{{\bf L} ^{\geq i}} \}} $$
between the set $\widehat{\bf Cv}(\Theta _i)$ of convex subgroups of the quotient group $\Theta _i = \Theta _{{\bf L} ^{\geq i}} / \Theta _{{\bf L} ^{> i}}$ and the subset $\{ \Delta \in {\bf Cv}(\Gamma ) \ | \ \Theta _{{\bf L} ^{>i}} \subset \Delta \subset  \Theta _{{\bf L} ^{\geq i}} \}$.

We can consider the group $\Theta _i$ as an ordered subgroup of $\Theta$ by the injective morphism $\xymatrix{\vartheta _i : \Theta _i \ar@<-1.5pt>@{^{(}->}[r] & \Theta}$ defined by $\vartheta _i (x_i) = \underline y$ where the element $\underline y = ( y_j) _{j \in I}$ satisfies $y_j = 0$ for $j \not= i$ and $y_i = x_i$. 
Then  the image by ${\bf F}_i$ of a convex subgroup $\Delta _i$ of $\Theta _i$ is the smallest convex subgroup of $\Theta$ that contains $\vartheta _i (\Delta _i)$, more precisely we have the equality 
${\bf F}_i(\Delta _i) = \{ \underline x \in \Theta \ | \ \iota _{\underline x} \geq i \ \hbox{and} \ x_i \in \Delta _i \}$, and the image $\Delta _i$ of a convex subgroup $\Delta$ of $\Theta$ by the inverse isomorphisms ${\bf G} _i$ of ${\bf F}_i$ is defined by ${\bf G}_i(\Delta ) = \vartheta _i ^{-1}(\Delta )$. 

\vskip .2cm

More generally let ${\bf L}$ and ${\bf L}'$ be two cuts of $I$, we have the inequality ${\bf L} \leq {\bf L}'$ if and only if we have the inclusion $\Theta _{{\bf L}'} \subset \Theta _{\bf L}$, then the quotient group $\Theta _{\bf L} / \Theta _{{\bf L}'}$ is isomorphic to the group $\prod _{i \in {\bf L} _+ \cap {\bf L}' _-} ^{(H)} \Theta _i$ and there is an isomorphism between the set $\widehat{\bf Cv}(\Theta _{\bf L} / \Theta _{{\bf L}'})$ and the subset $\{ \Delta \in {\bf Cv}(\Gamma ) \ | \ \Theta _{{\bf L}'} \subset \Delta \subset \Theta _{{\bf L} } \}$.

\vskip .2cm 

\begin{remark}\label{rmq:comparaison-coupures}
For any cut ${\bf L} = ({\bf L} _- , {\bf L} _+)$ of  $I$ we have the following equivalences:
$$i \in {\bf L} _+ \ \Longleftrightarrow \ \Theta _{{\bf L} ^{\geq i}} \subset \Theta _{\bf L} \quad\hbox{and}
\quad i \in {\bf L} _- \ \Longleftrightarrow \ \Theta _{\bf L} \subset \Theta _{{\bf L}  ^{> i}} \ .$$
\end{remark} 

\vskip .2cm

Let $\Delta$ be a convex subgroup of $\Theta =\prod _{i \in I} ^{(H)} \Theta _i$, and let ${\bf L} ^{\scriptscriptstyle \geq (\Delta)} = ({\bf L} ^{\scriptscriptstyle > (\Delta)} _- , {\bf L} ^{\scriptscriptstyle > (\Delta)} _+)$ and ${\bf L} ^{\scriptscriptstyle \geq (\Delta)} = ({\bf L} ^{\scriptscriptstyle \geq (\Delta)} _- , {\bf L} ^{\scriptscriptstyle \geq (\Delta)} _+)$ be the two cuts of $I$ defined by the following: 
$$\begin{array}{l} 
{\bf L} ^{\scriptscriptstyle > (\Delta)} _- = \{ i \in I \ | \ \Delta \subset \Theta _{{\bf L} ^{>i}} \} \quad \hbox{and}\quad  
{\bf L} ^{\scriptscriptstyle > (\Delta)} _+ = \{ i \in I \ | \ \Theta _{{\bf L} ^{>i}} \subsetneq \Delta \} \\ 
{\bf L} ^{\scriptscriptstyle \geq (\Delta)} _- = \{ i \in I \ | \ \Delta \subsetneq \Theta _{{\bf L} ^{\geq i}} \} \quad \hbox{and}\quad  
{\bf L} ^{\scriptscriptstyle \geq (\Delta)} _+ = \{ i \in I \ | \  \Theta _{{\bf L} ^{\geq i}}  \subset \Delta \} \ . 
\end{array}$$ 

Let $\Delta$ and $\Delta '$ be two convex subgroups of $\Theta$ with the inclusion $\Delta ' \subset \Delta$, then we have the inequalities 
${\bf L} ^{\scriptscriptstyle > (\Delta)} \leq {\bf L} ^{\scriptscriptstyle > (\Delta ')}$ and ${\bf L} ^{\scriptscriptstyle \geq (\Delta)} \leq {\bf L} ^{\scriptscriptstyle \geq (\Delta ')}$, hence the inclusions: 
$$\Theta _{{\bf L} ^{\scriptscriptstyle > (\Delta ')}} \subset \Theta _{{\bf L} ^{\scriptscriptstyle > (\Delta)}} \quad\hbox{and}\quad 
\Theta _{{\bf L} ^{\scriptscriptstyle \geq (\Delta ')}} \subset \Theta _{{\bf L} ^{\scriptscriptstyle \geq (\Delta)}} \ .$$

\vskip .2cm 

\begin{proposition}\label{prop:coupures-associees} 
(1) Let $\Delta$ be a convex subgroup of $\Theta$, then we have the inclusions
$$\Theta _{{\bf L} ^{\scriptscriptstyle \geq (\Delta)}} \subset \Delta \subset \Theta _{{\bf L} ^{\scriptscriptstyle > (\Delta)}} \ ,$$ 
hence the inequality ${\bf L} ^{\scriptscriptstyle > (\Delta)} \leq {\bf L} ^{\scriptscriptstyle \geq (\Delta)}$.  

(2) If the convex subgroup $\Delta$ is the subgroup $\Theta _{\bf L}$ associated with a cut ${\bf L}$ of $I$ we have the equalities 
$$ {\bf L} ^{\scriptscriptstyle > (\Delta)} = {\bf L} ^{\scriptscriptstyle \geq (\Delta)}= {\bf L} \ .$$ 

(3) If the convex subgroup $\Delta$ is not a subgroup associated with a cut, there exists $i \in I$ such that we have the equalities 
$${\bf L}^{\scriptscriptstyle > (\Delta)} = {\bf L}^{\geq i} \ \hbox{and} \ 
{\bf L}^{\scriptscriptstyle \geq (\Delta)} = {\bf L}^{>i} \ ,$$ 
and we have the strict inclusions
$\Theta _{{\bf L} ^{>i}} \subsetneq \Delta \subsetneq \Theta _{{\bf L} ^{\geq i}}$.  

\end{proposition} 

\begin{preuve} 
(1) For any element $\underline x$ of $\Theta$ let $i=\iota _{\underline x}$ be the smallest element of $Supp(\underline x)$, then we have:

$$\begin{array}{lcl}
x \in \Theta _{{\bf L} ^{\scriptscriptstyle \geq (\Delta)}} & \Longleftrightarrow & i \in {\bf L} ^{\scriptscriptstyle \geq (\Delta)} _+ \\ 
&  \Longleftrightarrow & \Theta _{{\bf L}^{\geq i}} \subset \Delta \ \Longrightarrow \ \underline x \in \Delta \ \hbox{because} \ x \in \Theta _{{\bf L}^{\geq i}} \ ,
\end{array}$$

$$\begin{array}{lcl}
x \notin \Theta _{{\bf L} ^{\scriptscriptstyle > (\Delta)}} & \Longleftrightarrow & i \in {\bf L} ^{\scriptscriptstyle > (\Delta)} _- \\ 
&  \Longleftrightarrow & \Delta \subset  \Theta _{{\bf L}^{> i}}  \ \Longrightarrow \ \underline x \notin \Delta \ \hbox{because} \ x \notin \Theta _{{\bf L}^{> i}} \ .
\end{array}$$
\vskip .2cm 

(2) If we have the equality $\Delta = \Theta _{\bf L}$, we deduce from the definitions that we have: 
$$\begin{array}{l} 
i \in {\bf L} ^{\scriptscriptstyle > (\Delta)}_- \ \Longleftrightarrow \  \Theta _{\bf L} \subset \Theta _{{\bf L}^{> i}} \ \Longleftrightarrow \  {\bf L}^{> i} \leq {\bf L} \\
i \in {\bf L} ^{\scriptscriptstyle \geq (\Delta)}_- \ \Longleftrightarrow \  \Theta _{\bf L} \subsetneq \Theta _{{\bf L}^{\geq i}} \ \Longleftrightarrow \  {\bf L}^{\geq i} < {\bf L} \ ,
\end{array}$$ 
hence the equalities ${\bf L}^{\scriptscriptstyle > (\Delta)}_- = {\bf L}^{\scriptscriptstyle \geq (\Delta)}_- = {\bf L}_-$ from remark \ref{rmq:coupure} \emph{(1)}. 
\vskip .2cm 

(3) If the subgroup $\Delta$ is not associated with a cut of $I$, we have a strict inequality ${\bf L} ^{\scriptscriptstyle > (\Delta)} < {\bf L} ^{\scriptscriptstyle \geq (\Delta)}$, this is equivalent to say that the set ${\bf L} ^{\scriptscriptstyle > (\Delta)} _+ \cap {\bf L} ^{\scriptscriptstyle \geq (\Delta)} _- = \{ i \in I \ | \ \Theta _{{\bf L} ^{>i}} \subsetneq \Delta \subsetneq  \Theta _{{\bf L} ^{\geq i}} \}$ is non-empty, then it is reduced to a unique element $i$. 

Then the cut ${\bf L} ^{\scriptscriptstyle > (\Delta)}$ is the immediate predecessor of the cut ${\bf L} ^{\scriptscriptstyle \geq (\Delta)}$, and by lemma \ref{le:predecesseur-successeur-dans-CpI}, this is equivalent to the fact that there exists $i \in I$ such that ${\bf L} ^{\scriptscriptstyle > (\Delta)} = {\bf L} ^{\geq i}$ and ${\bf L} ^{\scriptscriptstyle \geq (\Delta)} = {\bf L} ^{> i}$.

\hfill$\Box$ 
\end{preuve}

\vskip .2cm

\begin{theorem}\label{th:sous-groupe-isole-du-produit}  
Any convex subgroup $\Delta$ of the ordered group $\Theta =\prod _{i \in I} ^{(H)} \Theta _i$ is of one of the following two forms: 
\begin{enumerate}
\item $\Delta$ is the convex subgroup $\Theta _{{\bf L}}$ associated with a cut ${\bf L}$ of $I$, where the whole group $\Theta$ corresponds to the cut $(\emptyset, I)$ and $(0)$ corresponds to the cut $(I,\emptyset )$;  
\item there exists a non-trivial proper convex subgroup $\Delta _i$ of the group $\Theta _i$, for some $i \in I$, such that $\Delta = {\bf F}_i(\Delta _i)$. 
\end{enumerate} 
\end{theorem}

\begin{preuve} 
We deduce from proposition \ref{prop:coupures-associees} that if the convex subgroup $\Delta$ is not associated with a cut, there exists $i \in I$ such that we have the strict inclusions 
$\Theta _{{\bf L}^{> i}} \subsetneq \Delta \subsetneq \Theta _{{\bf L}^{\geq i}}$. 
Then there exists a non-trivial proper convex subgroup $\Delta _i$ of $\Theta _i$ such that we have the equality $\Delta = {\bf F}_i(\Delta _i)$.    

\hfill$\Box$ 
\end{preuve} 

\vskip .2cm

\begin{remark}\label{rmq:sous-groupe-isole-du-produit}  
(1) Let $\Delta$ be a convex subgroup of the Hahn product $\Theta =\prod _{i \in I} ^{(H)} \Theta _i$, then we are in one of the following two cases, which are not mutually exclusive: 
\begin{enumerate} [label=(\alph*)]
\item $\Delta$ is the convex subgroup $\Theta _{{\bf L}}$ associated with a cut ${\bf L}$ of $I$, then we have 
$$\textstyle{\Delta =\prod _{i \in {\bf L}_+} ^{(H)} \Theta _i};$$
\item there exists $i$ in $I$ and a convex subgroup $\Delta _i$ of $\Theta _i$ such that $\Delta = {\bf F}_i(\Delta _i)$, i.e. we have 
$$\textstyle{\Delta = \Delta _i \times \prod _{j > i} ^{(H)} \Theta _j}.$$
For $\Delta _i =(0)$ we have $\Delta = \Theta _{{\bf L}^{> i}}$, 
and for $\Delta _i = \Theta _i$ we have $\Delta = \Theta _{{\bf L}^{\geq i}}$.  

\noindent For $(0) \subsetneq \Delta _i \subset \Theta _i$ we have ${\bf L}^{\scriptscriptstyle > (\Delta)} = {\bf L}^{\geq i}$, and for $(0) \subset \Delta _i \subsetneq \Theta _i$ we have ${\bf L}^{\scriptscriptstyle \geq (\Delta)} = {\bf L}^{> i}$.
\end{enumerate} 

\vskip .2cm 

(2) We can also give the following descriptions: 
\begin{enumerate} [label=(\alph*)]
\item If the cut ${\bf L}^{\scriptscriptstyle > (\Delta)}$ is not of the form ${\bf L}^{\geq i}$ for any $i$ in $I$, the convex subgroup $\Delta$ is the convex subgroup $\Theta _{{\bf L}^{\scriptscriptstyle > (\Delta)}}$ associated with the cut ${\bf L}^{\scriptscriptstyle > (\Delta)}$. 

\item If there exists $i$ in $I$ such that we have the equality ${\bf L}^{\scriptscriptstyle > (\Delta)} = {\bf L}^{\geq i}$, the convex subgroup $\Delta$ is equal to ${\bf F}_i(\Delta _i)$, where $\Delta _i$ is the non-trivial convex subgroup of $\Theta _i$ defined by $\Delta _i = \vartheta _i ^{-1}(\Delta )$. 
\end{enumerate} 

and 

\begin{enumerate} [label=(\alph*)']
\item If the cut ${\bf L}^{\scriptscriptstyle \geq (\Delta)}$ is not of the form ${\bf L}^{> i}$ for any $i$ in $I$, the convex subgroup $\Delta$ is the convex subgroup $\Theta _{{\bf L}^{\scriptscriptstyle \geq (\Delta)}}$ associated with the cut ${\bf L}^{\scriptscriptstyle \geq (\Delta)}$. 

\item If there exists $i$ in $I$ such that we have the equality ${\bf L}^{\scriptscriptstyle \geq (\Delta)} = {\bf L}^{> i}$, the convex subgroup $\Delta$ is equal to ${\bf F}_i(\Delta _i)$, where $\Delta _i$ is the proper convex subgroup of $\Theta _i$ defined by $\Delta _i = \vartheta _i ^{-1}(\Delta )$. 
\end{enumerate} 
  
\end{remark}

\vskip .2cm

\begin{corollary} \label{cor:sous-groupe-isole-du-produit} 
If $\Theta =\prod _{i \in I} ^{(H)} \Theta _i$ is the Hahn product of a family of ordered groups of rank one, the only convex subgroups of $\Theta$ are the groups associated with the cuts of the set $I$. 
More precisely we have an isomorphism of totally ordered set between the sets $\widehat{\bf Cv}(\Theta )$ and ${\bf Cp}(I)^{\rm op}$.  
\end{corollary}

\begin{preuve}  
We have just to notice that if we have the inequality ${\bf L} \leq {\bf L} '$ in ${\bf Cp}(I)$, we have the inclusion ${\bf L} '_+ \subset {\bf L} _+$, hence the inclusion $\Theta _{{\bf L} '} \subset \Theta _{{\bf L}}$. 

\hfill$\Box$ 
\end{preuve}

\vskip .2cm 

In the case of a finite family $(\Theta _1, \ldots , \Theta _n)$ of ordered groups of rank one, the Hahn product $\Theta =\prod _{1 \leq i \leq n} ^{(H)} \Theta _i$ is equal to the direct sum equipped with the lexicographic order: 
$$ \Theta = \{ \underline x = (x_1, \ldots , x_n) \ | \ x_i \in \Theta _i \} \ ,$$
with $\underline x \leq \underline y$ if and only of there exists $i$, $0 \leq i \leq n-1$, such that  $x_1 = y_1 , \ldots x_i = y_i$ and $x_{i+1} \leq y_{i+1}$.

The convex subgroups of $\Theta$ are the subgroups $\Delta _i$, $0 \leq i \leq n-1$, defined by  
$$\Delta _i = \{ \underline x = (x_1, \ldots , x_n) \ | \ x_j = 0  \  \hbox{for} \ 1 \leq j \leq n-i \} \ .$$

\vskip .2cm 

\subsection{Principal convex subgroups}  

Let $\Gamma$ be an ordered group, as above we denote ${\bf Cv}(\Gamma )$ the totally ordered set of its proper convex subgroups, and $\widehat{\bf Cv}(\Gamma )$ the set of all the convex subgroups, the inclusion ${\bf Cv}(\Gamma )\subset \widehat{\bf Cv}(\Gamma )$ induces an injective morphism  $\xymatrix{u_! : {\bf Cp}({\bf Cv}(\Gamma )) \ar@<-1.5pt>@{^{(}->}[r] & {\bf Cp}(\widehat{\bf Cv}(\Gamma ))}$ defined by $u_!(({\bf L} _- , {\bf L} _+)) = ({\bf L} _- , {\bf L} _+ \cup \{ \Gamma \})$. 
The only cut of $\widehat{\bf Cv}(\Gamma )$ which is not in the image of $u_!$ is the trivial cut $(\Gamma , \emptyset )$. 

\begin{lemma}\label{le:cv-complet-et-cocomplet} 
The totally ordered sets ${\bf Cv}(\Gamma )$ and $\widehat{\bf Cv}(\Gamma )$ are complete and the set $\widehat{\bf Cv}(\Gamma )$ is cocomplete. 
\end{lemma}  

\begin{preuve} 
We deduce from remark \ref{rmq:inclusion-intervals} that the intersection and the union of any family of convex subgroups of $\Gamma$ are also convex subgroups. 

Then we can define the supremum and the infimum of any subset $X$ of ${\bf Cv}(\Gamma )$ or of $\widehat{\bf Cv}(\Gamma )$ respectively by the union and the intersection: 
$${\rm sup}\, X = \bigcup _{\Delta \in X} \Delta \quad\hbox{and}  
\quad {\rm inf}\, X = \bigcap _{\Delta \in X} \Delta \ .$$

\hfill$\Box$ 
\end{preuve}

\begin{proposition} \label{prop:principal-convex-subgroup}
(1) Let ${\bf L} = ( {\bf L} _-, {\bf L} _+)$ be a non-trivial cut of ${\bf Cv}(\Gamma )$, then the sets $\Delta _{{\bf L}}^{-}$ and $\Delta _{{\bf L}}^{+}$ defined by  
$$\Delta _{{\bf L}} ^- = \bigcup _{\Delta \in {\bf L} _-} \Delta \quad\hbox{and}\quad \Delta _{{\bf L}} ^+ = \bigcap _{\Delta \in {\bf L} _+} \Delta$$ 
are two convex subgroups of $\Gamma$ satisfying $\Delta _{{\bf L}}^{-} \subset\Delta _{{\bf L}}^{+}$. We have $\Delta _{{\bf L}}^{-} \not= \Delta _{{\bf L}}^{+}$ if and only if the cut ${\bf L}$ is a jump, and in that case the quotient group $\Theta _{{\bf L}} = \Delta _{{\bf L}}^+ / \Delta _{{\bf L}}^-$ is an ordered group of rank one.

(2) The subgroups $\Delta _{{\bf L}}^{-}$ and $\Delta _{{\bf L}}^{+}$ are distinct if and only if there exits an element $\xi$ in $\Gamma$, with $\xi \not= 0$, such that the cut ${\bf L}$ is defined by  
${\bf L} _+ =  {\bf Cv}_{\{\xi\}}(\Gamma ) = \{ \Delta \in {\bf Cv}(\Gamma ) \ | \ \xi \in \Delta \}$ and ${\bf L} _- = {\bf Cv}_{\{\xi\}}(\Gamma ) ^C = \{ \Delta \in {\bf Cv}(\Gamma ) \ | \ \xi \notin \Delta \}$.  
\end{proposition}

\begin{preuve} 
(1) We deduce from lemma \ref{le:cv-complet-et-cocomplet} that the subgroups $\Delta ^- _{{\bf L}}$ and $\Delta ^+ _{{\bf L}}$ correspond respectively to the element $\lambda _-({\bf L})$ and $\lambda _+({\bf L})$ defined in paragraph 1.3.
The inclusion $\Delta _{{\bf L}}^{-} \subset\Delta _{{\bf L}}^{+}$ is a consequence of the inequality $\lambda _-({\bf L}) \leq \lambda _-({\bf L})$ of proposition \ref{prop:composition} (3). 

Moreover if we have inequality $\Delta _{{\bf L}}^{-} \not=\Delta _{{\bf L}}^{+}$ the convex subgroup $\Delta _{{\bf L}}^{-}$ is an immediate predecessor of $\Delta _{{\bf L}}^{+}$, hence there doesn't exist any convex subgroup strictly between $\Delta _{{\bf L}}^-$ and $\Delta _{{\bf L}}^+$, and we deduce from remark \ref{rmq:sous-groupes-isoles-du-quotient} (2) that the the quotient group $\Theta _{{\bf L}}$ is an ordered group of rank one. 

\vskip .2cm 

(2) If the cut ${\bf L}$ is defined by 
${\bf L} _+ =  {\bf Cv}_{\{\xi\}}(\Gamma )$ and ${\bf L} _- =   {\bf Cv}_{\{\xi\}}(\Gamma ) ^C$,  
the element $\xi$ belongs to $\Delta _{{\bf L}}^+$ and doesn't belong to $\Delta _{{\bf L}}^-$.  
Conversely if we assume that there exists an element $\xi$ of $\Gamma$ that belongs to $\Delta _{{\bf L}}^+$ and that doesn't belong to $\Delta _{{\bf L}}^-$, 
for any $\Delta$ in ${\bf L} _+$ we have $\xi \in \Delta$, hence we have the inclusion ${\bf L} _+ \subset {\bf Cv}_{\{\xi\}}(\Gamma )$, 
in the same way for any $\Delta$ in ${\bf L} _-$ we have $\xi \notin \Delta$, hence we have the inclusion ${\bf L} _- \subset {\bf Cv}_{\{\xi\}}(\Gamma )^C$.    

\hfill$\Box$ 
\end{preuve}  

\begin{remark}\label{rmq:jump2}
(1) If ${\bf L}$ and ${\bf L} '$ are two cuts of ${\bf Cv}(\Gamma )$ with ${\bf L} \leq {\bf L} '$ we have the inclusions $\Delta _{{\bf L}}^{-} \subset \Delta _{{\bf L} '}^{-}$ and $\Delta _{{\bf L}}^{+} \subset \Delta _{{\bf L} '}^{+}$. 

(2) For the trivial cut ${\bf L} = (\emptyset , {\bf Cv}(\Gamma ))$ we can define the subgroups $\Delta _{{\bf L}}^{-}$ and $\Delta _{{\bf L}}^{+}$ in the same way, and we have $\Delta _{{\bf L}}^{-} = \Delta _{{\bf L}}^{+} =(0)$. 
For the trivial cut ${\bf L} = ({\bf Cv}(\Gamma ), \emptyset ))$ there are two cases to consider: 
if the set  ${\bf Cv}(\Gamma)$ of proper convex subgroups admits a greatest element $\Theta$  we have $\Delta _{{\bf L}}^{-} = \Delta _{{\bf L}}^{+} = \Theta$, and if there doesn't exist a greatest element we have $\Delta _{{\bf L}}^{-} = \Gamma$ but the subgroup $\Delta _{{\bf L}}^+$ is not well defined. By analogy we sometimes pose $\Delta _{{\bf L}}^+ = \Gamma$ in that case.  
\end{remark}

\vskip .2cm 

Let $\xi$ be an element of $\Gamma$, $\xi \not= 0$, and ${\bf L} = ({\bf L} _- , {\bf L} _+)$ be the cut of ${\bf Cv}(\Gamma )$ defined by ${\bf L} _+ =  {\bf Cv}_{\{\xi\}}(\Gamma )$ and ${\bf L} _- = {\bf Cv}_{\{\xi\}}(\Gamma ) ^C$, then we denote $\Delta _{{\bf L}} ^- = \Delta _{\{\xi\}} ^-$, $\Delta _{{\bf L}} ^+ = \Delta _{\{\xi\}} ^+$ and $ \Theta _{{\bf L}} = \Theta _{\{\xi\}}$.   

We have $\Delta _{\{\xi\}}^+ = \Delta _{\{-\xi\}}^+$, and if we assume the inequality  $\xi > 0$, we have  
$$\Delta _{\{\xi\}}^+ \ = \ \{ \zeta \in \Gamma \ | \ \exists n \in \N^* \ \hbox{with} \  -n \xi \leq \zeta \leq n \xi  \} \ .$$

For $\xi =0$ we can define the set ${\bf Cv} _{\{0\}} (\Gamma)$ as above and we have ${\bf Cv} _{\{0\}} (\Gamma) = {\bf Cv} (\Gamma)$, then the convex subgroup $\Delta _{\{0\}}^+$ is the subgroup $(0)$, but we cannot define a convex subgroup $\Delta _{\{0\}}^-$. If $\xi$ is an element of $\Gamma$ that doesn't belong to any proper convex subgroup of $\Gamma$ we have ${\bf Cv} _{\{\xi\}} (\Gamma) = \emptyset$ and $\Delta _{\{\xi\}}^+$ is equal to the group $\Gamma$. 

\vskip .2cm 

\begin{remark}\label{rmq:autre-definition}
By definition for $\xi \not= 0$ we have the following equalities: 
$$\Delta _{\{\xi\}} ^- \ = \ \bigcup _{\xi \notin \Delta} \Delta \quad\hbox{and}\quad  \Delta _{\{\xi\}} ^+ \ = \ \bigcap _{\xi \in \Delta} \Delta \ , $$ 
the subgroup $\Delta _{\{\xi\}} ^+$ is the smallest convex subgroup of $\Gamma$ that contains $\xi$ and the subgroup $\Delta _{\{\xi\}}^-$ is the biggest convex subgroup of $\Gamma$ that doesn't contain $\xi$. 

Hence for any element $\xi$ and any convex subgroup $\Delta$ we have the following equivalence:  
$$\xi \in \Delta \ \Longleftrightarrow \ \Delta _{\{\xi\}}^+ \subset \Delta \ . $$
\end{remark}

\vskip .2cm 

Let $\xi$ and $\zeta$ be elements of $\Gamma \setminus \{ 0 \}$, we put 
$$\xi \simeq _{Ar} \zeta \ \Longleftrightarrow \ \Delta _{\{\xi\}}^+ = \Delta _{\{\zeta\}}^+ \ \Longleftrightarrow \ \Delta _{\{\xi\}}^- = \Delta _{\{\zeta\}}^- \ ,$$
and we say that the two elements $\xi$ and $\zeta$ of $\Gamma$ are \emph{archimedean equivalent}. 
If the two elements $\xi$ and $\zeta$ are positive, this is also equivalent to say that there exist positive natural numbers $m$ and $n$ such that we have the relation $\xi < n \zeta < m \xi$.  

\vskip .2cm 

\begin{lemma} \label{le:inclusion-stricte}
Let $\xi$ and $\zeta$ in $\Gamma$ such that $\xi$ doesn't belong to $\Delta ^+ _{\zeta}$, then we have the inclusion $\Delta ^+_{\zeta} \subset \Delta ^- _{\xi}$. 
\end{lemma} 

\begin{preuve} 
This is a consequence of the equality $\Delta _{\{\xi\}} ^- \ = \ \bigcup _{\xi \notin \Delta} \Delta$ (c.f. remark \ref{rmq:autre-definition}). 

\hfill$\Box$ 
\end{preuve} 

\vskip .2cm 

\begin{definition} 
The convex subgroups $\Delta _{\{\xi\}} ^+$, for $\xi \in \Gamma$, $\xi \not= 0$, are called \emph{principal} convex subgroups of $\Gamma$. 
We note ${\bf Pr}(\Gamma )$ the set of principal convex subgroups of $\Gamma$, it is a subset of the totally ordered set $\widehat{\bf Cv}(\Gamma )$. 
\end{definition}  

\vskip .2cm 

\begin{remark}\label{rmq:immediate predecessor}
By proposition \ref{prop:composition} (3) we see that a convex subgroup $\Delta$ is principal if and only if $\Delta$ admits an immediate predecessor in ${\bf Cv}(\Gamma )$. 
In particular $(0)$ is never a principal convex subgroup.  
\end{remark} 

\vskip .2cm 

\begin{remark}\label{rmq:immediate predecessor2}
A convex subgroup $\Delta$ of $\Gamma$ is an immediate predecessor if and only if there exits a convex subgroup $\Delta '$ with $\Delta \subset \Delta '$ and the quotient group $\Delta ' / \Delta$ is of rank one. 

A convex subgroup $\Delta$ of $\Gamma$ is an immediate successor if and only if there exits a convex subgroup $\Delta '$ with $\Delta ' \subset \Delta$ and the quotient group $\Delta / \Delta '$ is of rank one. 
\end{remark} 

\vskip .2cm 

\begin{lemma}\label{le:dense-above2} 
Let $\Delta$ and $\Delta '$ be two convex subgroups of $\Gamma$, with $\Delta \subsetneq \Delta '$, then there exists $\xi \in \Gamma$ such that we have $\Delta \subsetneq \Delta _{\{\xi\}} ^+ \subset \Delta '$. 
\end{lemma} 

\begin{preuve} 
It is enough to choose an element $\xi$ that belongs to $\Delta '$ and not to $\Delta$. 

\hfill$\Box$ 
\end{preuve} 

We can state the previous lemma by saying that the subset ${\bf Pr}(\Gamma )$ of principal subgroups is dense above in the set $\widehat{\bf Cv}(\Gamma )$ of convex subgroups, and we deduce from lemma \ref{le:dense-above} that for any convex subgroup $\Delta$ we have the following equality:  
$$\Delta \ = \ \bigcup _{\Delta _{\{\xi\}} ^+ \subset \Delta} \Delta _{\{\xi\}} ^+ 
\ = \ \bigcup _{\xi \in \Delta} \Delta _{\{\xi\}} ^+  \ .$$

\vskip .2cm 

\begin{proposition} \label{prop:cv=cpct} 
The subset ${\bf Pr}(\Gamma )$ of principal convex subgroups is equal to the set ${\rm Suc}(\widehat{\bf Cv}(\Gamma ))$ of immediate successors of the set $\widehat{\bf Cv}(\Gamma )$ of all convex subgroups of $\Gamma$, and there is a natural isomorphism $\xymatrix{f:\widehat{\bf Cv}(\Gamma ) \ar[r] & {\bf Cp}({\bf Pr}(\Gamma ))}$, compatible with the inclusion of ${\bf Pr}(\Gamma )$ in $\widehat{\bf Cv}(\Gamma )$. 
\end{proposition} 

\begin{preuve} 
We deduce from remark \ref{rmq:immediate predecessor} that the subset ${\bf Pr}(\Gamma )$ is equal to ${\rm Suc}(\widehat{\bf Cv}(\Gamma ))$ and from lemma \ref{le:dense-above2} that it is dense above in $\widehat{\bf Cv}(\Gamma )$, then the proposition is a consequence of remark \ref{rmq:successeur-dense-above} and of proposition \ref{prop:natural-isomorphism}. 

\hfill$\Box$ 
\end{preuve} 

The existence of the isomorphism $\xymatrix{f:\widehat{\bf Cv}(\Gamma ) \ar[r] & {\bf Cp}({\bf Pr}(\Gamma ))}$ means that for any cut ${\bf L} = ( {\bf L}_- , {\bf L}_+ )$ of ${\bf Pr}(\Gamma )$, there exists a unique convex subgroup $\Delta$ of $\Gamma$ such that the initial segment ${\bf L}_- $ of ${\bf Pr}(\Gamma )$ is equal to the subset $\{ \Theta \in {\bf Pr}(\Gamma ) \ | \ \Theta \subset \Delta \}$.  
Moreover this convex subgroup $\Delta$ is defined by $\Delta = \bigcup _{\Theta \in {\bf L}_-} \Theta$.

\vskip .2cm 

\begin{corollary}\label{successeur-predecesseur} 
Let $\Delta$ be a convex subgroup of $\Gamma$, then $\Delta$ is an immediate successor, $\Delta = {\rm suc}(\Delta ')$, if and only if there exists $\xi$ in $\Gamma \setminus \{0\}$ with $\Delta = \Delta _{\{\xi\}} ^+$, and $\Delta$ is an immediate predecessor, $\Delta = {\rm pre}(\Delta '')$,  if and only if there exists $\xi$ in $\Gamma \setminus \{0\}$ with $\Delta = \Delta _{\{\xi\}} ^-$. 
\end{corollary}

\begin{preuve} 
It is a direct consequence of the first part of proposition \ref{prop:cv=cpct}. 

\hfill$\Box$ 
\end{preuve} 

\vskip .2cm 

\begin{remark}\label{rmq:contientT} 
Any subset $T$ of $\Gamma$ defines a cut ${\bf L} =( {\bf L} _-, {\bf L} _+)$ of ${\bf Cv}(\Gamma )$ by the following ${\bf L} _+ = {\bf Cv}_T(\Gamma ) = \{ \Delta \in {\bf Cv}(\Gamma ) \ | \ T \subset \Delta \}$ and ${\bf L} _- = {\bf Cv}_T(\Gamma )^C$. 
Then we have $\Delta _{{\bf L}} ^- \not= \Delta _{{\bf L}}^+$ if and only if there exists an element $\xi$ in $\Gamma$ such that we have the equality ${\bf Cv}_T(\Gamma ) ={\bf Cv}_{\{\xi\}}(\Gamma )$, i.e. the smallest convex subgroup containing $T$ is a principal convex subgroup of $\Gamma$.  
\end{remark} 

\vskip .2cm 

Any symmetric interval $S$ of $\Gamma$ defines a cut ${\bf L}^{(S)} = \bigl ( {\bf L}^{(S)} _- , {\bf L}^{(S)} _+ \bigr )$ of ${\bf Cv}(\Gamma )$ in the following way: 
$${\bf L}^{(S)} _- = \{ \Delta \ | \ \Delta \subset S \} \quad\hbox{and}\quad {\bf L}^{(S)} _+ = \{ \Delta \ | \ S \subsetneq \Delta \} \ ,$$ 
which is non trivial if $S \not= \emptyset$ or $S \not= \Gamma$.  

The convex subgroup $\Delta _{{\bf L}^{(S)}} ^- = \bigcup _{\Delta \in {\bf L}^{(S)}_-} \Delta$ belongs to ${\bf L}^{(S)}_-$, it is the greatest convex subgroup of $\Gamma$ which is included in $S$, and we have ${\bf L}^{(S)} _- = \{ \Delta \ | \ \Delta \subset \Delta _{{\bf L}^{(S)}} ^- \}$. 

\begin{proposition}\label{prop:coupure-intervalle-symetrique} 
If $S$ is not a convex subgroup of $\Gamma$ the cut ${\bf L}^{(S)}$ is a jump. 
The subgroup $\Delta _{{\bf L}^{(S)}} ^+$ is the smallest convex subgroup of $\Gamma$ which contains $S$, it is a principal convex subgroup, and the subgroup $\Delta _{{\bf L}^{(S)}} ^-$ is its immediate predecessor in ${\bf Cv}(\Gamma )$. 
\end{proposition} 

\begin{preuve} 
If  $S$ is not a convex subgroup, we have the strict inclusion $\Delta _{{\bf L}^{(S)}} ^- \subsetneq S$. Let $\xi$ be an element of $S \setminus \Delta _{{\bf L}^{(S)}} ^-$, then $\xi$ belongs to any subgroup $\Delta$ in ${\bf L}^{(S)} _+$, hence $\xi \in \Delta _{{\bf L}^{(S)}} ^+ = \bigcap _{\Delta \in {\bf L}^{(S)}_+} \Delta$, and we get $\Delta _{{\bf L}^{(S)}} ^- = \Delta _{\{\xi\}} ^- \not= \Delta _{\{\xi\}} ^+ = \Delta _{{\bf L}^{(S)}} ^+$. 

\hfill$\Box$ 
\end{preuve} 

\vskip .2cm 

\begin{lemma}\label{le:phi-et-delta+}
For any symmetric interval $S$ of $\Gamma$ and any $\zeta$ in $\Gamma$, we have the equivalence: 
$$\zeta \in \Delta _{{\bf L}^{(S)}} ^- \ \Longleftrightarrow \ \Delta^+ _{\{\zeta\}} \subset S 
\ \Longleftrightarrow \  n \zeta \in S \ \hbox{for all} \ n \in \N ^* \ .$$
\end{lemma} 

\begin{preuve} 
We have the following equivalences: 
$$\begin{array}{rcl}
\zeta \in \Delta _{{\bf L}^{(S)}} ^- & \Longleftrightarrow & \Delta^+ _{\{\zeta\}} \subset \Delta _{{\bf L}^{(S)}} ^-  ,\ \hbox{by remark \ref{rmq:autre-definition}} \\ 
& \Longleftrightarrow & \Delta ^+ _{\{\zeta\}} \subset S \ . 
\end{array} $$
Moreover, if $\zeta$ is a positive element of $\Gamma$, the group $ \Delta ^+ _{\{\zeta\}}$ is defined by 
$$\Delta _{\{\zeta\}}^+ \ = \ \{ \gamma \in \Gamma \ | \ \exists n \in \N^* \ \hbox{with} \  -n \zeta \leq \gamma \leq n \zeta  \} \ ,$$ 
hence the result. 

\hfill$\Box$ 
\end{preuve} 

\vskip .2cm 

For any subset $X$ of the group $\Gamma$, and for any integer $n \in \N ^*$ we can define the subset $(1/n) X$ of $\Gamma$ by $(1/n) X = \{ \gamma \in \Gamma \ | \ n \gamma \in X \}$, and the subset $X _{\rm tor}$ by 
$$X _{\rm tor} \ = \ \bigcap _{n \in \N ^*} (1/n) X \ = \  \{ \gamma \in \Gamma \ | \ \forall n \in \N ^*\ , n \gamma \in X \} \ .$$
Then we deduce from lemma \ref{le:phi-et-delta+} the equality $\Delta _{{\bf L}^{(S)}} ^- = S _{\rm tor}$. 

\vskip .2cm 

\begin{remark}\label{rmq:initial-segment-convex-subgroup}  
As $S$ is symmetric interval, for any $n \geq m$ we have the inclusion $(1/n)S \subset (1/m) S$, then we can write $S _{\rm tor} = \varprojlim (1/n) S$. 
\end{remark}

\vskip .2cm 

We consider now the totally ordered group $\Theta = \prod _{i \in I} ^{(H)} \Theta _i$ obtained as the Hahn product of a family of totally ordered groups, indexed by a totally ordered set $I$. 
From remark \ref{rmq:sous-groupe-isole-du-produit}, we know that a convex subgroup $\Delta$ of $\Theta$ is either a convex subgroup $\Theta _{{\bf L}}$ associated with a cut ${\bf L}$ of $I$ that is not of the form ${\bf L} ^{\geq i}$ for some $i \in I$, or corresponds to a non-trivial convex subgroup $\Delta _i$ of the group $\Theta _i$ for some $i \in I$.

\begin{proposition}\label{prop:predecesseur-dans-produit-de-Hahn} 
Let $\Delta$ and $\Delta '$ be two convex subgroups of the Hahn product $\Theta =\prod _{i \in I} ^{(H)} \Theta _i$ such that $\Delta$ is an immediate predecessor of $\Delta '$, then there exists $i$ in $I$ and two convex subgroups $\Delta _i$ and $\Delta _i'$ of $\Theta _i$ such that $\Delta _i$ is an immediate predecessor of $\Delta _i'$. 

For $\Delta _i = (0)$, we have $\Delta = \Theta _{{\bf L}^{>i}}$ and $\Delta '_i$ is a convex subgroup of $\Theta _i$ of rank one, with possibly $\Delta '_i=\Theta _i$. 

For $\Delta '_i = \Theta _i$, we have $\Delta ' = \Theta _{{\bf L}^{\geq i}}$ and $\Delta _i$ is a convex subgroup of $\Theta _i$ such that the quotient group $\Theta _i/\Delta _i$ of rank one. 
\end{proposition}

\begin{preuve} 
Let ${\bf L} ^{\scriptscriptstyle > (\Delta)}$ and ${\bf L} ^{\scriptscriptstyle \geq (\Delta)}$ and $ {\bf L} ^{\scriptscriptstyle > (\Delta ')}$ and ${\bf L} ^{\scriptscriptstyle \geq (\Delta ')}$ be the cuts of $I$ associated respectively with the convex subgroups $\Delta$ and $\Delta '$, and by proposition \ref{prop:coupures-associees} we have the inclusions $\Theta _{{\bf L} ^{\scriptscriptstyle \geq (\Delta)}} \subset \Delta \subset \Theta _{{\bf L} ^{\scriptscriptstyle > (\Delta)}}$ and $\Theta _{{\bf L} ^{\scriptscriptstyle \geq (\Delta ')}} \subset \Delta ' \subset \Theta _{{\bf L} ^{\scriptscriptstyle > (\Delta ')}}$. 

If $\Delta$ is the convex subgroup $\Theta _{\bf L}$ associated with a cut ${\bf L}$, we have the equality ${\bf L} ^{\scriptscriptstyle > (\Delta)} = {\bf L} ^{\scriptscriptstyle \geq (\Delta)} = {\bf L}$, then we get the inclusions $\Delta = \Theta _{{\bf L} ^{\scriptscriptstyle \geq (\Delta )}} \subset \Theta _{{\bf L} ^{\scriptscriptstyle \geq (\Delta ')}} \subset \Delta '$. 

As $\Delta$ is the immediate predecessor of $\Delta '$ we must have 
\begin{enumerate}[label=(\alph*)]
\item \emph{either $\Delta = \Theta _{{\bf L} ^{\scriptscriptstyle \geq (\Delta ')}} \subsetneq \Delta '$ and the equality ${\bf L} ^{\scriptscriptstyle \geq (\Delta ')}=  {\bf L} ^{\scriptscriptstyle \geq (\Delta)} = {\bf L}$,} 
\item \emph{or $\Delta \subsetneq \Theta _{{\bf L} ^{\scriptscriptstyle \geq (\Delta ')}} = \Delta '$ and $\Theta _{\bf L} = \Theta _{{\bf L} ^{\scriptscriptstyle \geq (\Delta)}}$ is the immediate predecessor of $\Theta _{{\bf L} ^{\scriptscriptstyle \geq (\Delta ')}}$}
\end{enumerate} 

(a) By proposition \ref{prop:coupures-associees} we deduce from the strict inclusion $\Theta _{{\bf L} ^{\scriptscriptstyle \geq (\Delta ')}} \subsetneq \Delta '$ that there exists $i$ in $I$ such that we have the equality ${\bf L} = {\bf L} ^{\scriptscriptstyle \geq (\Delta ')} = {\bf L} ^{> i}$, and we have $\Delta = \prod _{j>i} ^{(H)} \Theta _j$. 

Moreover $\Delta '$ corresponds to a proper subgroup $\Delta '_i$ of $\Theta _i$ and as $\Delta '$ is the immediate successor of $\Delta$, which corresponds to the subgroup $(0)$ of $\Theta _i$, the rank of $\Delta '_i$ is equal to one. 

(b) If $\Theta _{\bf L}$ is the immediate predecessor of $\Theta _{{\bf L} ^{\scriptscriptstyle \geq (\Delta ')}}$, we deduce from the injective morphisms of ordered sets $\xymatrix{{\bf Cp}(I)^{\rm op} \ar@<-1.5pt>@{^{(}->}[r]& \widehat{\bf Cv}(\Theta )}$ defined by $\xymatrix{{\bf L} \ar@{|->}[r] & \Theta _{\bf L}}$, that the cut ${\bf L}$ is an immediate successor in ${\rm Cp}(I)$, hence by lemma \ref{le:predecesseur-successeur-dans-CpI} there exists $i$ in $I$ such that the cut ${\bf L}$ is equal to ${\bf L} ^{> i}$.  

In that case the subgroup $\Delta '_i$ is equal to the group $\Theta_i$, and the rank of $\Theta _i$ is equal to one.  

\vskip .2cm 

If the convex subgroup $\Delta$ of $\Theta$ is not associated with a cut, there exists $i$ in $I$ such that we have $\Theta _{{\bf L} ^{> i}} \subsetneq \Delta \subsetneq \Theta _{{\bf L} ^{\geq i}}$, and $\Delta$ and $\Delta '$ are the image by the morphism ${\bf F}_i$ of non-trivial convex subgroups $\Delta _i$ and $\Delta '_i$ of the group $\Theta _i$. 
Then $\Delta$ is the immediate predecessor in $\widehat{\bf Cv}(\Theta )$ of $\Delta '$ if and only if $\Delta _i$ is the immediate predecessor in $\widehat{\bf Cv}(\Theta _i)$ of $\Delta s'-i$.    

\hfill$\Box$ 
\end{preuve}

\vskip .2cm

\begin{corollary}\label{cor:sous-groupe-principal-de-produit-de-hahn}
A convex subgroup $\Delta$ of  the Hahn product $\Theta = \prod _{i \in I} ^{(H)} \Theta _i$ is principal if and only if 
there exists $i$ in $I$ such that we have the equality ${\bf L} ^{\scriptscriptstyle > (\Delta)} = {\bf L}^{\geq i}$ and the convex subgroup $\Delta _i$ of $\Theta _i$ defined by $\Delta _i = \vartheta _i ^{-1}(\Delta )$  is a principal convex subgroup of $\Theta _i$.  
\end{corollary}

\begin{preuve} 
We know by remark \ref{rmq:immediate predecessor} that a convex subgroup $\Delta$ of $\Theta$ is principal if and only if there exists a convex subgroup $\Delta '$ such that $\Delta$ is the immediate successor of $\Delta '$. 
Then the result is a consequence of proposition \ref{prop:predecesseur-dans-produit-de-Hahn}.

\hfill$\Box$ 
\end{preuve} 

From corollary \ref{cor:sous-groupe-principal-de-produit-de-hahn} we deduce that a convex subgroup $\Delta$ of $\Theta = \prod _{i \in I} ^{(H)} \Theta _i$ associated with a cut ${\bf L}$ is principal if and only if there exists $i$ in $I$ such that ${\bf L} = {\bf L}^{\geq i}$ and the group $\Theta _i$ is a principal subgroup of it self.

\vskip .2cm 

\begin{remark}\label{rmq:predecesseur-produit-de-hahn}  
In a similar way we can show that a convex subgroup $\Delta$ of  the Hahn product $\Theta = \prod _{i \in I} ^{(H)} \Theta _i$ is an immediate predecessor if and only if 
there exists $i$ in $I$ such that we have the equality ${\bf L} ^{\scriptscriptstyle \geq (\Delta)} = {\bf L}^{> i}$ and the convex subgroup $\Delta _i$ of $\Theta _i$ defined by $\Delta _i = \vartheta _i ^{-1}(\Delta )$  is an immediate predecessor in the set of convex subgroups of $\Theta _i$.  
\end{remark} 

\vskip .2cm

\begin{corollary}\label{cor:non-predecesseur} 
Let $\Delta$ be a convex subgroup of the Hahn product $\Theta =\prod _{i \in I} ^{(H)} \Theta _i$ , then $\Delta$ is not an immediate predecessor if and only if we are in the following cases: 
\begin{enumerate}
\item if $\Delta$ is not associated to a cut of $I$, the convex subgroup $\Delta _i$ of a group $\Theta _i$ associated with $\Delta$ is not an immediate predecessor;  

\item if $\Delta$ is associated with a cut $\bf L$ of $I$, $\Delta = \Theta _{\bf L}$, this cut is not an immediate successor in ${\bf Cp}(I)$, i.e. ${\bf L}$ is not of the form ${\bf L}^{>i}$. 
\end{enumerate}
\end{corollary}

\begin{preuve} 
It is an immediate consequence of proposition \ref{prop:predecesseur-dans-produit-de-Hahn}. 

\hfill$\Box$ 
\end{preuve}

\vskip .2cm

\begin{corollary}\label{cor:sous-groupe-principal-de-produit-de-hahn-de-rang-un} 
If $\Theta$ is the Hahn product $\prod _{i \in I} ^{(H)} \Theta _i$ of a family of totally ordered groups of rank one, there exists an isomorphism of ordered set $\xymatrix{ I^{\rm op} \ar[r] & {\bf Pr}(\Theta )}$, defined by $\xymatrix{i \ar@{|->}[r] & \Theta _{{\bf L}^{\geq i}}}$. 
\end{corollary}

\begin{preuve} 
It is an immediate consequence of corollary \ref{cor:sous-groupe-principal-de-produit-de-hahn}. 

\hfill$\Box$ 
\end{preuve}

We can also recover this result from corollary \ref{cor:sous-groupe-isole-du-produit} and proposition \ref{prop:cv=cpct}. 
 
\vskip .2cm

\begin{remark}\label{rmq:sous-groupe-principal-associe}
Let $\Delta$ be a principal convex subgroup of the Hahn product $\Theta = \prod _{i \in I} ^{(H)} \Theta _i$, and let $\Delta _i$ be the principal convex subgroup of $\Theta _i$ associated with it.  
If $\Delta _i$ is the principal convex subgroup of the group $\Theta _i$ associated with an element $\xi$ of $\Theta _i$, i.e. if we have the equality $\Delta _i = \Delta ^+ _{\{\xi\}}$, $\Delta$ is the principal convex subgroup of $\Theta$ associated with any element $\underline x = ( x_j ) _{j \in I}$ of $\Theta$ satsifying $x _j =0$ for $j < i$ and $x_i = \xi$. 
\end{remark} 

\vskip .2cm 

We define the totally ordered set ${\bf I}_{\Gamma}$ as the set ${\bf Pr}(\Gamma )$ of principal convex subgroups of $\Gamma$ endowed with the opposite total order induced by the order on ${\bf Cv}(\Gamma )$, i.e. we have ${\bf Pr}(\Gamma )= \{ \Delta _i^+ \ ; \ i \in {\bf I}_{\Gamma} \}$ and $i \leq j$ if and only if $\Delta _j^+ \subset \Delta _i^+$.  
The smallest element $0$ of ${\bf I}_{\Gamma}$, if it exists, corresponds to the group $\Gamma$, $\Delta _0^+= \Gamma$, this element doesn't exist if and only if any element $\xi$ of $\Gamma$ belongs to a proper subgroup of $\Gamma$. 
We define the totally ordered set $\widehat{\bf I}_{\Gamma}$ obtained by adding an element $\infty$ with $i < \infty$ for any $i \in {\bf I}_{\Gamma}$, this element corresponds to the subgroup $(0)$, and then we can identify the quotient set $\Gamma / \simeq _{Ar}$ with the set $\widehat{\bf Pr}(\Gamma ) = {\bf Pr} (\Gamma ) \cup \{ (0) \}$.  

For any $i \in {\bf I}_{\Gamma}$ there exists $\xi$ belonging to $\Gamma$ such that we have $\Delta _i^+ = \Delta ^+ _{\{\xi\}}$, then we denote $\Delta _i^-$ the convex subgroup $\Delta _i^- = \Delta ^- _{\{\xi\}}$ and $\Theta _i$ the ordered group of rank one defined by $\Theta _i = \Theta _{\{\xi\}} = \Delta ^+ _i / \Delta ^- _i$. 

\vskip .2cm 

\begin{theorem}\label{th:structure-sur-groupe-totalement-ordonne} 
Let $\Gamma$ be an ordered group and let ${\bf I}_{\Gamma}$ be the totally ordered set defined above, then the family $\Delta _{{\bf I}_{\Gamma}} = \bigl ( \Delta _i^- , \Delta _i^+ \bigr ) _{i \in I _{\Gamma}}$ defines a non-degenerate ${\bf I}_{\Gamma}$-structure on $\Gamma$.  
\end{theorem} 

\begin{preuve} 
By construction for any $\xi$ belonging to $\Gamma$ we have the inclusions $ \Delta _{\{\xi\}}^- \subset \Delta _{\{\xi\}}^+$, hence for any $i$ in ${\bf I}_{\Gamma}$ the inclusions $\Delta _i^- \subset \Delta _i^+$. 

We deduce from lemma \ref{le:inclusion-stricte} that for any $i<j$ in ${\bf I}_{\Gamma}$ we have the inclusion $\Delta _j^+ \subset \Delta _i^-$. 

For any $\xi$ belonging to $\Gamma$ the convex subgroup $\Delta _i^+ = \Delta _{\{\xi\}}^+$ satisfies $\xi \in \Delta _{\{\xi\}}^+$ and $\xi \notin \Delta _{\{\xi\}}^-$.  
Conversely any subgroup $\Delta _i^+$ is of the form $\Delta _i^+ = \Delta _{\{\xi\}}^+$, hence we have always $\Delta _i^- \subsetneq \Delta _i^+$.   

\hfill$\Box$ 
\end{preuve} 

The skeleton of the ordered group $\Gamma$ for this ${\bf I}_{\Gamma}$-structure is the pair $\bigl ( {\bf I}_{\Gamma} , (\Theta _i) _{i \in {\bf I}_{\Gamma}} \bigr )$, where each group $\Theta _i$ is an ordered group of rank one of the form $\Theta _i = \Theta _{\{\xi\}}$ for an element $\xi$ of $\Gamma$

\vskip .2cm 

Let $\xymatrix{\varphi : \Gamma \ar[r] & \Gamma '}$ be an injective morphism of ordered groups, we define a morphism $\varepsilon _{\varphi}$ from the set of convex subgroups of $\Gamma$ into the set of convex subgroups of $\Gamma '$ in the following way:   
for any $\Delta$ in $\widehat{\bf Cv}(\Gamma )$ its image $\varepsilon _{\varphi} ( \Delta)$ in ${\bf Cv}(\Gamma ')$ is the convex subgroup of $\Gamma '$ generated by $\Delta$.  
  
\begin{lemma}\label{le:injection-de-groupes}
The map $\xymatrix{\varepsilon _{\varphi} : \widehat{\bf Cv}(\Gamma ) \ar[r] & \widehat{\bf Cv}(\Gamma ' )}$ is an injective morphism of totally ordered sets, the image of the principal convex subgroup $\Delta ^+ _{\{\xi\}}$ of $\Gamma$ is the principal convex subgroup $\Delta ^+ _{\{\varphi(\xi )\}}$ of $\Gamma '$, 
and the morphism $\varphi$ induces an injective morphism of ordered groups of rank one: 
$\xymatrix{\varphi _{\{\xi \}} : \Theta _{\{\xi\}} \ar[r] & \Theta _{\{\varphi(\xi )\}} }$. 
\end{lemma}

\begin{preuve} 
By definition the subgroup $\varepsilon _{\varphi} (\Delta )$ is equal to:  
$$\varepsilon _{\varphi} (\Delta )  \ = \ \{ \zeta \in \Gamma ' \ | \ \exists \xi \in \Delta \cap \Gamma _{\geq 0} \  \ \hbox{with} \  - \varphi (\xi ) \leq \zeta \leq \varphi (\xi )  \} \ ,$$  
hence the map $\xymatrix{\varepsilon _{\varphi} : \widehat{\bf Cv}(\Gamma ) \ar[r] & \widehat{\bf Cv}(\Gamma ' )}$ respects the inclusion relation, and it is injective because we have the equality $\Delta = \varepsilon _{\varphi} (\Delta ) \cap \Gamma$.  

The image by $\varepsilon _{\varphi}$ of the convex subgroup $\Delta ^+ _{\{\xi\}}$ of $\Gamma$ is the convex subgroup $\Delta ^+ _{\{\varphi(\xi )\}}$ of $\Gamma '$, and the image of the convex subgroup $\Delta ^- _{\{\xi\}}$ is contained in the convex subgroup $\Delta ^- _{\{\varphi(\xi )\}}$, then we have a morphism $\xymatrix{\varphi _{\{\xi \}} : \Delta ^+ _{\{\xi\}} / \Delta ^- _{\{\xi\}} \ar[r] & \Delta ^+ _{\{\varphi (\xi )\}} / \Delta ^- _{\{\varphi (\xi )\}} }$, and we deduce the injectivity of $\varphi _{\{\xi \}}$ from the equality $\Delta ^- _{\{\xi\}} = \Delta ^- _{\{\varphi (\xi )\}} \cap \Gamma$.  

\hfill$\Box$ 
\end{preuve} 

\begin{definition}
An injective morphism of ordered groups $\xymatrix{\varphi : \Gamma \ar[r] & \Gamma '}$ is called an \emph{immediate extension} if the map $\varepsilon _{\varphi}$ is a bijection such that for any principal convex subgroup $\Delta ^+ _{\{\xi\}}$ the morphism $\varphi _{\{\xi \}}$ is an isomorphism. 
\end{definition} 

\vskip .2cm 

\begin{remark}\label{rmq:immediate=immediate}
An immediate morphism of ordered groups $\xymatrix{\varphi : \Gamma \ar[r] & \Gamma '}$ is an immediate $I$-morphism where $\Gamma$ and $\Gamma '$ are endowed with the $I$-structures defined above with $I = {\bf I} _{\Gamma} \simeq {\bf I} _{\Gamma '}$.  
\end{remark} 

\vskip .2cm

\subsection{Divisible ordered group an discrete ordered group} 

A torsion-free abelian group $\Gamma$ is \emph{divisible} if for any $\xi \in \Gamma$ and for any $n \in \N ^*$, there exists $\zeta \in \Gamma$ such that  $n\zeta =\xi$. 
The \emph{divisible hull} $\Gamma ^{\rm div}$ of a torsion-free abelian group $\Gamma$ is the minimal divisible group that contains it,  it can be defined as the group $\Gamma \otimes _{\Z} \Q$ or as the quotient set of $\{(\xi, n) \in \Gamma \times \N ^* \}$ by the equivalence relation $(\xi ,n) \simeq (\xi ',n') \Longleftrightarrow n'\xi =n\xi '$.

If $\Gamma$ is an abelian totally ordered group, it is torsion-free and we can provide the group $\Gamma ^{\rm div}$ with an order that is compatible with the group law and that extends the one on $\Gamma$ in the following way. 

Let $\xi$ and $\zeta$ be two elements of $\Gamma ^{\rm div}$, let $(\xi _0, m)$ and $(\zeta _0, n)$ be respective representatives of $\xi$ and $\zeta$ in $\Gamma \times \N ^*$, then we define the order by  
$$\xi \leq \zeta \ \Longleftrightarrow \ n \xi _0 \leq m \zeta _0 \ .$$
We verify that we have thus provided the group $\Gamma ^{\rm div}$ with a total order and that the map $\xymatrix{\varphi : \Gamma \ar[r] & \Gamma ^{\rm div}}$ defined by $\varphi (\xi _0) = (\xi _0, 1)$ is an injective morphism of ordered groups.  

\begin{lemma}\label{le:cloture-divisible} 
The map $\xymatrix{\varepsilon _{\varphi} : \widehat{\bf Cv}(\Gamma ) \ar[r] & \widehat{\bf Cv}(\Gamma ^{\rm div} )}$ induced by $\varphi$ is an isomorphism of ordered sets,  in particular the groups $\Gamma$ and $\Gamma  ^{\rm div}$ have the same rank. 

It induces an isomorphism between the subsets ${\bf I}_{\Gamma}$ and ${\bf I}_{\Gamma ^{\rm div}}$ of principal convex subgroups of $\Gamma$ and $\Gamma ^{\rm div}$ respectively. 
For any $\xi \in \Gamma$ the morphism $\xymatrix{\varphi _{\{\xi \}} : \Theta _{\{\xi\}} \ar[r] & \Theta _{\{\varphi(\xi )\}} }$ is the natural morphism from $\Theta _{\{\xi \}}$ into its divisible hull. 
\end{lemma} 

\begin{preuve} 
The convex subgroup of $\Gamma ^{\rm div}$ generated by a convex subgroup $\Delta$ of $\Gamma$ is isomorphic to the group $\Delta ^{\rm div}$, and all the convex subgroups of $\Gamma ^{\rm div}$ are obtained in this way, we thus obtain the isomorphism of ordered sets: $\xymatrix{\varepsilon _{\varphi} : \widehat{\bf Cv}(\Gamma ) \ar[r] & \widehat{\bf Cv}(\Gamma ^{\rm div} )}$. 

We have the equalities  
$$\varepsilon _{\varphi}\bigl (\Delta ^+ _{\{\xi\}}\bigr ) = {\bigl (\Delta ^+ _{\{\xi\}}\bigr )} ^{\rm div} = \Delta ^+ _{\{\xi\}} \otimes _{\Z} \Q \quad\hbox{and}\quad \varepsilon _{\varphi}\bigl (\Delta ^- _{\{\xi\}}\bigr ) = {\bigl (\Delta ^- _{\{\xi\}}\bigr )} ^{\rm div} = \Delta ^- _{\{\xi\}} \otimes _{\Z} \Q \ ,$$
we deduce that the group $\Theta _{\{\varphi(\xi )\}}$ is equal to $\Theta _{\{\xi \}} \otimes _{\Z} \Q = \Theta _{\{\xi \}} ^{\rm div}$.  

\hfill$\Box$ 
\end{preuve}  

\vskip .2cm 

Let $R$ be a commutative ring, an ordered group $\Gamma$ is an \emph{ordered $R$-group} if the group $\Gamma$ is a module over $R$ and if all the convex subgroups of $\Gamma$ are submodules. 
We denote as above ${\bf I}_{\Gamma}$ the totally ordered set of the principal convex subgroups of $\Gamma$, with order opposite at the order induced by the inclusion. 
Then for any $i \in {\bf I}_{\Gamma}$ the convex subgroups $\Delta _i^+$ and $\Delta _i ^-$ are submodules of $\Gamma$ and the ordered group of rank one $\Theta _i$ obtained as the quotient $\Delta _i^+ / \Delta _i ^-$ is a module over $R$.  
We deduce from the foregoing that $\Gamma$ is a ${\bf I}_{\Gamma}$-module over $R$ assuming $\delta (\Gamma )_i^- = \Delta _i^-$ and $\delta (\Gamma )_i^+ = \Delta _i^+$. 

In particular any divisible ordered group $\Gamma$ is an ordered $\Q$-group. 

\vskip .2cm 

\begin{proposition}\label{prop:theoreme-de-hahn} 
Let $\Gamma$  be a divisible ordered group, then the family $\Delta _{{\bf I}_{\Gamma}} = \bigl ( \Delta _i^- , \Delta _i^+ \bigr ) _{i \in I _{\Gamma}}$ defines a non-degenerate ${\bf I}_{\Gamma}$-structure of $\Q$-vector space on $\Gamma$ and there exists a natural immediate ${\bf I}_{\Gamma}$-morphism of ordered $\Q$-vector spaces 
$$\xymatrix{\underline v _{\Gamma} : \Gamma \ar[r] & \varepsilon (\Gamma ) = \prod _{i \in {\bf I}_{\Gamma}} ^{(H)} \varepsilon (\Gamma )_i}$$
where the subgroups $\varepsilon (\Gamma )_i = \Theta _i$ are divisible ordered groups of rank one.  \end{proposition}  
 
\begin{preuve} 
This a consequence of the above and of corollary \ref{cor:theoreme-de-Hahn}. 

\hfill$\Box$ 
\end{preuve} 

\vskip .2cm 

\begin{remark}\label{rmq:cv=ct(pr)} 
We deduce directly from corollary \ref{cor:sous-groupe-isole-du-produit} that we have an isomorphism between the set $\widehat{\bf Cv}(\varepsilon (\Gamma ))$ of convex subgroups of $\varepsilon (\Gamma )$ and the set  ${\bf Cp}({\bf I}_{\Gamma})^{\rm op}$ of cuts of the set of principal convex subgroups of $\Gamma$. 
\end{remark}

\vskip .2cm 

Let $\Gamma$ be a totally ordered group and let $\Delta _{{\bf I}_{\Gamma}} = \bigl ( \Delta _i^-,\Delta _i^+ \bigr )$ be the ${\bf I}_{\Gamma}$-structure on $\Gamma$. We assume that we have an immediate ${\bf I}_{\Gamma}$-morphism $\xymatrix{\underline v _{\Gamma} : \Gamma \ar[r] & \varepsilon (\Gamma ) = \prod _{i \in {\bf I}_{\Gamma}} ^{(H)} \varepsilon (\Gamma )_i}$, and we denote $\xymatrix{\varepsilon _{\underline v _{\Gamma}} : \widehat{\bf Cv}(\Gamma ) \ar[r] & \widehat{\bf Cv}(\Theta )}$ the isomorphism of ordered sets defined above. 

Since the ordered groups $\varepsilon (\Gamma )_i =  \Delta _i^+ / \Delta _i^-$ are of rank one, the convex subgroups of $\Theta$ are the groups associated with the cuts of the set $I$ (cf. corollary \ref{cor:sous-groupe-isole-du-produit}).
More precisely for any convex subgroup $\Delta$ of $\Theta$ there exists a cut ${\bf L}$ of ${\bf I}_{\Gamma}$ such that we have the equality $\Delta = \Theta _{\bf L} = \prod _{i \in {\bf L}_+} ^{(H)} \varepsilon (\Gamma )_i$. 

By proposition \ref{prop:cv=cpct} we have an isomorphism $\xymatrix{f:\widehat{\bf Cv}(\Gamma ) \ar[r] & {\bf Cp}({\bf Pr}(\Gamma ))}$, which associates with any convex subgroup $\Delta$ of $\Gamma$ the cut ${\bf M}^{(\Delta )}$ of the set ${\bf Pr}(\Gamma )$ of principal convex subgroups of $\Gamma$ defined by ${\bf M}_-^{(\Delta )} = \{ \Pi \in {\bf Pr}(\Gamma ) \ | \ \Pi \subset \Delta \}$, and we have the equality $\Delta = \bigcup _{\Pi \in {\bf M}_-^{(\Delta )}} \Pi$. 

The ordered set ${\bf I}_{\Gamma}$ is equal to the set ${\bf Pr}(\Gamma )$ of principal convex subgroups of $\Gamma$ with opposite order, then we have a bijection between the sets of cuts respectively of ${\bf Pr}(\Gamma )$ and of ${\bf I}_{\Gamma}$, which is defined by the following. 
With any cut ${\bf M}$ of ${\bf Pr}(\Gamma )$ we associate the cut ${\bf L}$ of ${\bf I}_{\Gamma}$ defined by ${\bf L}_- = {\bf M}_+$ and ${\bf L}_+ = {\bf M}_-$, and we note this cut ${\bf M}^{op}$.

\begin{proposition}\label{prop:correspondance} 
The image by the isomorphism $\xymatrix{\varepsilon _{\underline v} : \widehat{\bf Cv}(\Gamma ) \ar[r] & \widehat{\bf Cv}(\Theta )}$ of the convex subgroup $\Delta$ of $\Gamma$ is the convex subgroup $\Theta _{\bf L}$ of $\Theta$ associated with the cut ${\bf L} = {\bigl ({\bf M}^{(\Delta )}\bigr )}^{op}$ of ${\bf I}_{\Gamma}$.  
\end{proposition}

\begin{preuve} 
The image by $\underline v_{\Gamma}$ of any elemnt $\xi$ of $\Gamma$ is an element $\underline x =(x_i)$ of $\Theta$ which satisfies 
$$\begin{array}{l} 
\iota _{\underline x} = \iota _{\Delta _{{\bf I}_{\Gamma}}} (\xi ) \quad\hbox{hence} \ x_j =0 \ \hbox{for } \ j< \iota _{\Delta _{{\bf I}_{\Gamma}}}(\xi ) \\ 
x_{\iota _{\underline x}} = {\bf in} _{\Delta _{{\bf I}_{\Gamma}}}(\xi ) \ ,
\end{array}$$ 
where $\iota _{\Delta _{{\bf I}_{\Gamma}}}(\xi )$ is defined as the index $i \in {\bf I}_{\Gamma}$ such that we have $\xi \in \Delta _i^+ \setminus \Delta _i^-$.

Let $\Delta$ be a convex subgroup of $\Gamma$ and let ${\bf L}$ be the cut ${\bigl ({\bf M}^{(\Delta )}\bigr )}^{op}$ of ${\bf I}_{\Gamma}$. 
We deduce from the equality $\Delta = \bigcup _{\Pi \in {\bf M}_-^{(\Delta )}} \Pi = \bigcup _{i \in {\bf L}_+} \Delta _i^+$ that the index $\iota _{\Delta _{{\bf I}_{\Gamma}}} (\xi )$ of any element $\xi$ of $\Delta$ belongs to ${\bf L}_+$, hence the element $\underline v_{\Gamma }(\xi )$ belongs to $\Theta _{\bf L}$. 
From definition of the application $\varepsilon _{\underline v_{\Gamma}}$ we deduce from the above that we have the inclusion $\varepsilon _{\underline v_{\Gamma}}(\Delta ) \subset \Theta _{\bf L}$. 

If we had a strict inclusion there would exist a cut ${\bf L}'$ of ${\bf I}_{\Gamma}$ with $\varepsilon _{\underline v_{\Gamma}}(\Delta ) \subset \Theta _{{\bf L}'}$ and ${\bf L}'<{\bf L}$, ie. ${\bf L}'_+ \subsetneq {\bf L}_+$. 
Then for any $i$ in ${\bf L}_+ \setminus {\bf L}'_+$ and any $\xi$ in $\Delta _i^+ \setminus \Delta _i^-$ we would get $\underline v _{\Gamma} (\xi ) \notin \Theta _{{\bf L}'}$. 

\hfill$\Box$  
\end{preuve}

\vskip .2cm 

Any ordered group of rank one $\Theta$ admits an injective morphism of ordered groups in the additive group of real numbers:  
$\xymatrix{c : \Theta \ar[r] & \R}$, that we call a \emph{completion} of $\Theta$. 
This morphism is not unique but if $\xymatrix{c : \Theta \ar[r] & \R}$ and $\xymatrix{c' : \Theta \ar[r] & \R}$ are two such injective morphisms, there exists a number $\alpha \in \R$, with $\alpha >0$, such that  for any $\xi$ belonging to $\Theta$ we have the equality $c'(\xi ) = \alpha c(\xi )$. 

\vskip .2cm 

Let $\Theta$ be an ordered group and $I$ be a totally ordered set, we define the ordered group $\Theta ^{[I]} _{\rm lex}$ as the Hahn product of \emph{$I$ copies} of $\Theta$: 
$$\Theta ^{[I]} _{\rm lex} = \hbox{$\prod _{i \in I} ^{(H)} \Theta$} \ ,$$ 
equipped with the lexicographic order.  

This group is endowed with a natural structure of totally ordered $I$-group.

\begin{corollary}\label{cor:complete} 
Let $\Gamma$ be an ordered group, of skeleton the pair $\bigl ( {\bf I}_{\Gamma} , (\Theta _i) _{i \in {\bf I}_{\Gamma}} \bigr )$, 
then there exists an injective morphism of totally ordered ${\bf I}_{\Gamma}$-groups  
$$\xymatrix{\psi : \Gamma \ar[r] & \Gamma ^{\rm cp} =  \R ^{[{\bf I}_{\Gamma}]} _{\rm lex}} \ ,$$
that induces for any $i$ in ${\bf I}_{\Gamma}$ a completion $\xymatrix{c_i : \Theta _i \ar[r] & \R}$ of the group $\Theta _i$. 
\end{corollary} 

\begin{preuve} 
The morphism $\psi$ is obtained as the composition of the morphism $\xymatrix{\varphi :\Gamma \ar[r] & \Gamma ^{\rm div}}$, the morphism $\xymatrix{\underline v _{\Gamma ^{\rm div}}: \Gamma ^{\rm div} \ar[r] & \varepsilon (\Gamma ^{\rm div})}$ defined in the proposition \ref{prop:theoreme-de-hahn}, and the morphisms of completion $\xymatrix{c_i : \Theta _i \ar[r] & \Theta _i ^{\rm div} \ar[r] & \R}$. 

\hfill$\Box$ 
\end{preuve} 

\vskip .2cm 

\begin{remark}\label{rmq:non-unicite} 
The ordered group $\Gamma ^{\rm cp} = \R ^{[{\bf I}_{\Gamma}]} _{\rm lex}$ is well defined, but the morphism $\xymatrix{\psi : \Gamma \ar[r] & \Gamma ^{\rm cp}}$ is not unique. 
Indeed, according to proposition \ref{prop:equivalence-decomposition-morphisme}, the morphism $v _{\Gamma ^{\rm div}}$depends essentially on the choice of the $I$-decomposition of the ${\bf I}_{\Gamma}$-structure on $\Gamma ^{\rm div}$. 
\end{remark}

\vskip .2cm

We can describe the structure of an ordered group of rank one. More precisely if $\Gamma$ is an ordered group of rank one, we are in one of the following two cases: 
\begin{enumerate} 
\item \emph{if the set $\Gamma _{>0}$ admits a smallest element $\gamma _0$ the group $\Gamma$ is discrete and is isomorphic to the group $\Z$;}  
\item \emph{otherwise the group $\Gamma$ is \emph{continuous}, i.e. for any $\gamma _1$ and $\gamma _2$ in $\Gamma$ with $\gamma _1 < \gamma _2$, there exists $\gamma _3$ in $\Gamma$ such that $\gamma _1 < \gamma _3 < \gamma _2$.} 
\end{enumerate}  
We want to generalize the notion of discrete group for any ordered group and we get the following definitions.   

\vskip .2cm 

\begin{definition} 
A totally ordered group $\Gamma$ is \emph{well-ranked}, or has a \emph{well-rank}, if the ordered set ${\bf Cv}(\Gamma )$ of convex subgroups of $\Gamma$ is well-ordered. 
\end{definition}  

As the ordered set $\widehat{\bf Cv}(\Gamma )$ is complete, the group $\Gamma$ is well-ranked if and only if any convex subgroup $\Delta$ of $\Gamma$ has an immediate successor $\bar \Delta$ in $\widehat{\bf Cv}(\Gamma )$. 
Hence for any convex subgroup $\Delta$ there exists $\xi$ in $\Gamma$ such that $\Delta = \Delta _{\{\xi\}} ^-$, and we denote $\bar\Delta$ the convex subgroup $\Delta _{\{\xi\}} ^+$, which is the immediate successor of $\Delta$ in $\widehat{\bf Cv}(\Gamma )$. 

We recall that the quotient group $\bar\Delta / \Delta$ is an ordered group of rank one.  

\vskip .2cm

\begin{definition} 
Let $\Gamma$ be a totally ordered group. 
\begin{enumerate} 
\item The group $\Gamma$ is \emph{discrete} if the set $\Gamma _{>0}$  of positive elements has a smallest element $\gamma _0$.  
\item The group $\Gamma$ is \emph{discretely ordered} if it is well-ranked and if for any convex subgroup $\Delta$ of $\Gamma$ the quotient group $\bar\Delta / \Delta$ is discrete. 
\end{enumerate}
\end{definition}  
 
\vskip .2cm 
 
Let $\Gamma$ be a discrete ordered group and let $\gamma _0$ be the smallest element of the set $\Gamma _{>0}$, then $\gamma _0$ is the immediate successor of $0$. Moreover any element $\xi$ of $\Gamma$ has an immediate successor ${\rm suc}(\xi )$ which is equal to $\xi + \gamma _0$, and an immediate predecessor ${\rm pre}(\xi )$ which is equal to $\xi -\gamma _0$.  

\begin{lemma}\label{le:groupe-discret} 
Let $\Gamma$ be a totally ordered group, then the following are equivalent. 
\begin{enumerate}
\item The group $\Gamma$ is discrete. 
\item The group $\Gamma$ has a convex subgroup $\Delta _0$ that is isomorphic to the group $\Z$ of natural integers. 
\item There exists $\sigma$ and $\omega$ in $\Gamma$ such that $\sigma$ is an immediate predecessor of $\omega$. 
\end{enumerate} 
\end{lemma} 

\begin{preuve} 
The convex subgroup $\Delta _0$ is the group generated by the element $\gamma _0$, which is an immediate successor of $0$. 

\hfill$\Box$ 
\end{preuve} 

\vskip .2cm  

\begin{remark}\label{rmq:(0)-predecesseur} 
If the ordred group $\Gamma$ is discrete the subgroup $\Delta _0$ is of rank one, then it is an immediate successor of the subgroup $(0)$, we deduce that the subgroup $(0)$ of a discrete ordered group $\Gamma$ is an immediate predecessor. 
\end{remark} 

\vskip .2cm 
 
There exists also the notion of \emph{generalized discrete} ordered group, introduced by P. Hill and J.L. Mott in \cite{H-M}. An ordered group $\Gamma $ is generalized discrete if for any element $\xi$ in $\Gamma$ the ordered group $\Delta _{\{\xi\}} ^+ / \Delta _{\{\xi\}} ^-$ is isomorphic to $\Z$, where $\Delta _{\{\xi\}} ^+$ is the principal convex subgroup of $\Gamma$ associated with $\xi$.  
This notion is weaker than the notion of discretely ordered group that we have defined.

\vskip .2cm

		      \section{Cuts of an ordered abelian group}
%

\subsection{Invariance subgroup}      
A large part of the results of this section are in the articles of F.V. Kuhlmann, see for instance \cite{Ku1} and \cite{Ku-N} . 
Let $\Gamma$ be a totally ordered group.

\begin{lemma}\label{le:sous-groupe-invariant} 
For any initial or final segment $\Sigma$ of $\Gamma$ the subset $\Delta (\Sigma)$ of $\Gamma$ defined by 
$$\Delta (\Sigma) = \{ \xi \in \Gamma \ | \  \Sigma + \xi = \Sigma \}$$ 
is a convex subgroup of $\Gamma$. 

Moreover we have the equality $\Delta (\Sigma ) = \Delta (\Sigma ^C)$. 
\end{lemma}

\begin{preuve} 
For any $X \subset \Gamma$, the subset $\Delta (X)$ of $\Gamma$ defined by 
$\Delta (X) = \{ \xi \in \Gamma \ | \ X + \xi = X \}$
is a subgroup, which is convex if $X$ is an initial or a final segment of $\Gamma$. 

\hfill$\Box$ 
\end{preuve} 

\vskip .2cm  

\begin{definition} 
Let $\Lambda = \bigl (\Lambda _- , \Lambda _+ \bigr )$ be a cut of $\Gamma$, then the convex subgroup $\Delta (\Lambda)$ defined by 
$$\Delta (\Lambda ) = \Delta (\Lambda _-) = \Delta (\Lambda _+) \ ,$$
is called the \emph{invariance subgroup} or the \emph{breadth} of the cut $\Lambda$ (cf. \cite{Ku1}). 

In the same way for any initial or final segment $\Sigma$ of $\Gamma$ we call the convex subgroup $\Delta (\Sigma )$ defined above the \emph{invariance subgroup} of $\Sigma$. 
\end{definition} 

\vskip .2cm 

\begin{remark}\label{rmq:coupure-triviale}
The invariance subgroup $\Delta (\Lambda )$ of the cut $\Lambda$ is equal to the group $\Gamma$ if and only if the cut $\Lambda$ is trivial. 
\end{remark}

Let $\Theta$ be a convex subgroup of $\Gamma$, and let $\xymatrix{w _{\Theta}: \Gamma \ar[r] & \bar\Gamma}$ be the surjective morphism of totally ordered groups where $\bar\Gamma$ is the quotient group $\Gamma ' = \Gamma / \Theta$. 
Then for any cut $\Lambda = \bigl (\Lambda _- , \Lambda _+ \bigr )$ of $\Gamma$ we deduce from lemma \ref{le:morphisme-surjectif}  that the images $w _{\Theta}({\Lambda} _-)$ and $w _{\Theta}(\Lambda _+)$ are respectively initial and final segments of the ordered group $\bar\Gamma$. 
We have the equality $w _{\Theta}(\Lambda _-) \cup w _{\Theta}(\Lambda _+) = \bar\Gamma$ but in general we don't have  $w _{\Theta}(\Lambda _-) \cap w _{\Theta}(\Lambda _+) = \emptyset$.  

\begin{proposition}\label{prop:image-de-coupure}  
If the set $w _{\Theta}(\Lambda _-) \cap w _{\Theta}(\Lambda _+)$ is non empty then it is reduced to a single element $\{ \bar\delta \}$, and 
if we have $w _{\Theta}(\Lambda _-) \cap w _{\Theta}(\Lambda _+) = \emptyset$ we get a cut $\bigl ( w _{\Theta}(\Lambda _- ), w _{\Theta}(\Lambda _+ )\bigr )$ of $\bar\Gamma$. 

Moreover we have the following equivalences 
$$\begin{array}{rcl}
w _{\Theta}(\Lambda _-) \cap w _{\Theta}(\Lambda _+) = \emptyset 
&  \Longleftrightarrow &  \Theta \subset \Delta (\Lambda ) \\ 
& \Longleftrightarrow & \Lambda _- + \Theta = \Lambda _-  \\ 
& \Longleftrightarrow & \Lambda _+ + \Theta = \Lambda _+  \ .
\end{array}$$
\end{proposition} 

\begin{preuve} 
The first part of the proposition is a consequence of lemma \ref{le:comparaison}. 

The equivalences are direct consequences of definition of the invariance subgroup $\Delta (\Lambda ) = \Delta (\Lambda _-) = \Delta (\Lambda _+)$. 

\hfill$\Box$ 
\end{preuve} 

\vskip .2cm

For any cut $\Lambda = \bigl (\Lambda _- , \Lambda _+ \bigr )$ of $\Gamma$ we define the subset ${\bf L}_- ^{(\Lambda )}$ of ${\bf Cv}(\Gamma )$ made up of the convex subgroups $\Delta$ of $\Gamma$ satisfying  $w _{\Delta}(\Lambda _-) \cap w _{\Delta}(\Lambda _+) = \emptyset$, a convex subgroup $\Delta$ of $\Gamma$ belongs to the subset ${\bf L}_+ ^{(\Lambda )}= {\bf Cv}(\Gamma ) \setminus {\bf L}_-^{(\Lambda )}$ if and only if there exist $\delta _-$ and $\delta _+$ respectively in $\Lambda _-$ and $\Lambda _+$, such that $(\delta _+ - \delta _-)$ belongs to $\Delta$. 
It easy to see that ${\bf L}_-^{(\Lambda )}$ is an initial segment, ${\bf L}_+^{(\Lambda )}$ is a final segment, and we get a cut ${\bf L} ^{(\Lambda )}=\bigl ({\bf L}_- ^{(\Lambda )}, {\bf L}_+^{(\Lambda )}\bigr )$ of ${\bf Cv}(\Gamma )$.  

We always have $\Delta =(0)$ that belongs to ${\bf L}_- ^{(\Lambda )}$, and the cut ${\bf L}^{(\Lambda)}$ is trivial if and only if any proper convex subgroup of $\Gamma$ belongs to ${\bf L}_- ^{(\Lambda )}$. 
When the set ${\bf Cv}(\Gamma )$ of proper convex subgroups admits a greatest element $\Theta$, it is equivalent to say that $\Theta$ belongs to ${\bf L}_- ^{(\Lambda )}$.  

If the cut ${\bf L} ^{(\Lambda )}$ of ${\bf Cv}(\Gamma )$ is non trivial, by proposition \ref{prop:principal-convex-subgroup} we can define two convex subgroups $\Delta _{{\bf L}^{({\Lambda})}}^{-}$ and $\Delta _{{\bf L}^{({\Lambda})}}^{+}$ of $\Gamma$ by 
$$\Delta _{{\bf L}^{({\Lambda})}}^{-} = \bigcup _{\Delta \in {\bf L}_- ^{(\Lambda )}} \Delta \quad\hbox{and}\quad \Delta _{{\bf L}^{({\Lambda})}}^{+} = \bigcap _{\Delta \in {\bf L}_+ ^{(\Lambda )}} \Delta \ .$$

We deduce from proposition \ref{prop:image-de-coupure} that the invariance subroup $\Delta (\Lambda)$ is the greatest element of the initial segment ${\bf L}_-^{(\Lambda )}$, we have the equality $\Delta _{{\bf L}^{({\Lambda})}}^{-} = \Delta (\Lambda )$, and we note $\overline \Delta (\Lambda)$ the subgroup $\Delta _{{\bf L}^{({\Lambda})}}^{+}$.

\vskip .2cm

\begin{proposition}\label{prop:passage-au-quotient}  
Let $\Sigma$ be a non trivial initial segment (resp. non trivial final segment) of $\Gamma$, let $\Theta$ be a convex subgroup of $\Gamma$, and let $\xymatrix{w_{\Theta}: \Gamma \ar[r] & \bar\Gamma = \Gamma / \Theta}$ be the quotient morphism. 
Then $\bar\Sigma = w_{\Theta}(\Sigma )$ is an initial segment (resp. a final segment) of $\bar\Gamma$ and we have the equivalence 
$$w_{\Theta} ^{-1} (\bar\Sigma ) = \Sigma \ \Longleftrightarrow \  \Theta \subset \Delta (\Sigma ) \ .$$  
 
In this case the initial segment (resp. final segment) $\bar\Sigma$ is a non trivial initial segment (resp. non trivial final segment) of $\bar\Gamma$ and for any element $\xi$ of $\Gamma$ we have the equivalence 
$$\xi \in \Sigma \ \Longleftrightarrow \  w_{\Theta}(\xi ) \in \bar\Sigma \ .$$ 
\end{proposition}

\begin{preuve} 
For any convex subgroup $\Theta$ of $\Gamma$, the image by the surjective morphism $w_{\Theta}$ of a non empty initial segment $\Sigma$ is a non empty initial segment $\bar\Sigma$, we have the equality $w_{\Theta}^{-1}(w_{\Theta}(\Sigma )) = \Sigma + \Theta$, and from the definition of the invariance subgroup $\Delta (\Sigma )$
we deduce the equivalence 
$$w_{\Theta} ^{-1} (\bar\Sigma ) = \Sigma \ \Longleftrightarrow \  \Theta \subset \Delta (\Sigma ) \ .$$ 

If we have $\Theta \subset \Delta (\Sigma ) = \Delta(\Sigma ^C)$ the final segment $w_{\Theta}(\Sigma ^C)$ is equal to $(w_{\Theta}(\Sigma ))^C$, then the initial segment $\bar\Sigma$ is not equal to $\bar\Gamma$, and it is non trivial. 
  
\hfill$\Box$ 
\end{preuve} 

\vskip .2cm

\begin{remark}\label{rmq:interval-invariance-group}
The complementary of the invariance subgroup $\Delta (\Lambda)$ in $\Gamma$ is equal to the union $(\Lambda _+ - \Lambda _-) \cup (\Lambda _- - \Lambda _+)$. 
Indeed if $\xi$ belongs to $(\Lambda _+ - \Lambda _-)$ there exist $\delta _- \in \Lambda _-$ and $\delta _+ \in \Lambda _+$ such that $\xi + \delta _- = \delta _+$, hence $\xi$ doesn't belong to $\Delta (\Lambda _-) = \Delta (\Lambda)$, and in the same way if $\xi$ belongs to $(\Lambda _- - \Lambda _+)$ it doesn't belong to $\Delta (\Lambda _+) = \Delta (\Lambda)$. 
Conversely if $\xi$ doesn't belong to $\Delta (\Lambda)$, we deduce from the relation $\xi + \Lambda _+ \not= \Lambda _+$ that $\xi$ belongs to $(\Lambda _+ - \Lambda _-)$ or to $(\Lambda _- - \Lambda _+)$. 

By construction the cut $\Pi = (\Pi _- , \Pi _+)$ of $\Gamma$ defined by $\Pi _+ = \Lambda _+ - \Lambda _-$ is positive, and the symmetric interval $\Xi _{\Pi}$ associated with $\Pi$ is equal to the invariance subgroup $\Delta (\Lambda)$.      
\end{remark}

\vskip .2cm 

\begin{definition} 
Two cuts ${\Lambda} = \bigl ({\Lambda} _- , {\Lambda} _+ \bigr )$ and ${\Lambda} '= \bigl ({\Lambda} '_- , {\Lambda} '_+ \bigr )$ of the totally ordered group $\Gamma$ are \emph{linearly equivalent} if there exists an element $\xi \in \Gamma$ with $\Lambda ' = \Lambda + \xi$. 
\end{definition}

We have thus defined an equivalence relation on the set ${\bf Cp}(\Gamma )$ of the cuts of $\Gamma$, which we note $\simeq$, and we note ${\bf Cp}^{[lin]}(\Gamma )$ the quotient set. 
We can also consider this relation $\simeq$ as an equivalence relation on the sets ${\bf IS}(\Gamma )$ and ${\bf FS}(\Gamma )$ of initial and final segments of the group $\Gamma$, and in the same way we define the quotient sets ${\bf IS}^{[lin]}(\Gamma )$ and ${\bf FS}^{[lin]}(\Gamma )$.
 
We can notice that the only cut linearly equivalent to the trivial cut $\Lambda = (\emptyset ,\Gamma )$ is itself, and likewise for the trivial cut $\Lambda = (\Gamma, \emptyset )$.

\begin{remark}\label{rmq:meme-groupe-d-invariance}
If $\Lambda$ and $\Lambda '$ are two linearly equivalent cuts, the subsets ${\bf L}^{({\Lambda})}_{-}$ and ${\bf L}^{({\Lambda '})}_{-}$ are the same and we have the equalities $\Delta (\Lambda) = \Delta (\Lambda ')$ and $\overline\Delta (\Lambda) = \overline\Delta (\Lambda ')$. 
\end{remark} 

\vskip .2cm 

If we have the strict inclusion $\Delta (\Lambda) \subsetneq \overline\Delta (\Lambda )$, i.e. $\Delta _{{\bf L}^{({\Lambda})}}^{-} \subsetneq \Delta _{{\bf L}^{({\Lambda})}}^{+}$, the subgroup $\Delta _{{\bf L}^{({\Lambda})}}^{+}$ is the immediate successor of the subgroup $\Delta _{{\bf L}^{({\Lambda})}}^{-}$, and we deduce from corollary \ref{successeur-predecesseur} that  there exists $\xi$ in $\Gamma \setminus \{0\}$ such that  $\Delta _{{\bf L}^{({\Lambda})}}^{+} = \Delta _{\{\xi\}} ^+$ and $\Delta _{{\bf L}^{({\Lambda})}}^{-} = \Delta _{\{\xi\}} ^-$, i.e. the cut ${\bf L} ^{(\Lambda )}$  is defined by  
${\bf L} ^{(\Lambda )} _+ =  {\bf Cv}_{\{\xi\}}(\Gamma )$ and ${\bf L} ^{(\Lambda )} _- = {\bf Cv}_{\{\xi\}}(\Gamma ) ^C$, and the quotient group $\Theta _{{\bf L}^{({\Lambda})}} = \Delta _{{\bf L}^{({\Lambda})}}^+ / \Delta _{{\bf L}^{({\Lambda})}}^-$ is an ordered group of rank one. 

If we have the equality $\Delta (\Lambda) = \overline\Delta (\Lambda )$, i.e. $\Delta _{{\bf L}^{({\Lambda})}}^{-} = \Delta _{{\bf L}^{({\Lambda})}}^{+}$, the subgroup $\Delta _{{\bf L}^{({\Lambda})}}^{+}$ belongs to the initial segment ${\bf L} ^{(\Lambda )} _-$, and the final segment ${\bf L} ^{(\Lambda )} _+$ has no smallest element. 

\vskip .2cm 

The cut ${\bf L} ^{(\Lambda )}=\bigl ({\bf L}_- ^{(\Lambda )}, {\bf L}_+^{(\Lambda )}\bigr )$ of ${\bf Cv}(\Gamma )$ associated with the cut $\Lambda$ of the group $\Gamma$ may also be defined in the following way. 
A convex subgroup $\Theta$ of $\Gamma$ belongs to the final segment ${\bf L}_+^{(\Lambda )}$ (resp. to the initial segment ${\bf L}_-^{(\Lambda )}$) if and only if we have $\Pi _+ \cap \Theta \not= \emptyset$ (resp. $\Theta _{\geq 0} \subset \Pi _-$), where $\Pi = (\Pi _- , \Pi _+)$ is the cut defined in remark \ref{rmq:interval-invariance-group}. 

\vskip .2cm

Let $\Lambda$ be a cut of the group $\Gamma$, let $\Delta (\Lambda )$ be the invariance subgroup of $\Lambda$ and let $\Theta$ be a convex subgroup of $\Gamma$. 
We consider the morphism $\xymatrix{w_{\Theta}: \Gamma \ar[r] & \bar\Gamma = \Gamma /\Theta}$, and the pair $w_{\Theta}(\Lambda )  = (w_{\Theta}(\Lambda _-) , w _{\Theta}(\Lambda _+))$ of subsets of the quotient group $\Gamma '$, then we have the equivalence: 

\centerline{$w_{\Theta}(\Lambda )$is a cut of $\Gamma '$ if and only if we have the inclusion $\Theta \subset \Delta (\Lambda )$.}

\begin{lemma}\label{le:sous-groupe-invariant-quotient} 
We assume that we have the inclusion $\Theta \subset \Delta (\Lambda )$, then the invariance subgroup $\Delta (w_{\Theta}(\Lambda ))$ of the cut $w_{\Theta}(\Lambda )$ of $\bar\Gamma$ is equal to the quotient group $\Delta (\Lambda )/\Theta$. 
\end{lemma}  

\begin{preuve} 
This is a direct consequence of the previous characterization of the invariance subgroup of a cut. 

\hfill$\Box$ 
\end{preuve} 

\vskip .2cm

Let ${\Lambda} = ({\Lambda} _- ,{\Lambda} _+ )$ be a cut of an ordered group $\Gamma$, and let $\Theta$ be a convex subgroup of $\Gamma$. 
We can define for any element $\delta$ of $\Gamma$ a pair $\bigl (\Lambda ^{\delta}_- (\Theta ), \Lambda ^{\delta}_+ (\Theta ) \bigr )$ of subsets of the convex subgroup $\Theta$ by 
$$\begin{array}{l}
\Lambda ^{\delta}_- (\Theta ) = \{ \zeta \in \Theta \ \hbox{such that } \ \delta + \zeta \in \Lambda _- \}  = (\Lambda _- -{\delta}) \cap \Theta \\
\Lambda ^{\delta}_+ (\Theta ) = \{ \zeta \in \Theta \ \hbox{such that } \ \delta + \zeta \in \Lambda _+ \}  = (\Lambda _+ -{\delta}) \cap \Theta\ , 
\end{array}$$
and we have le following result. 

\begin{lemma}\label{le:coupure} 
The pair $\bigl (\Lambda ^{\delta}_- (\Theta ), \Lambda ^{\delta}_+ (\Theta ) \bigr )$ is a cut of $\Theta$. Moreover this cut is non trivial if and only if $\delta '=w _{\Theta}(\delta )$ belongs to $w _{\Theta}(\Lambda _-) \cap w _{\Theta}(\Lambda _+)$. 
\end{lemma}

\begin{preuve} 
As $\Lambda$ is a cut of $\Gamma$ it is easy to deduce that $\bigl (\Lambda ^{\delta}_- (\Theta ), \Lambda ^{\delta}_+ (\Theta ) \bigr )$ is a cut of $\Theta$. 
This cut $\bigl (\Lambda ^{\delta}_- (\Theta ), \Lambda ^{\delta}_+ (\Theta ) \bigr )$ is non trivial if and only if there exist $\delta _1$ and $\delta _2$ in $T=w _{\Theta}^{-1}(\delta ')$ such that $\delta _1$ belongs to $\Lambda _-$ and $\delta _2$ belongs to $\Lambda _+$.

\hfill$\Box$ 
\end{preuve} 

\vskip .2cm

We consider now a convex subgroup $\Theta$ of $\Gamma$ with $\Delta (\Lambda ) \subsetneq \Theta$. 
The pair $(w_{\Theta}(\Lambda _-) , w _{\Theta}(\Lambda _+))$ of subsets of the quotient group $\Gamma '= \Gamma / \Theta$ is not a cut, there exists a unique element $\delta '$ in the intersection $ w_{\Theta}(\Lambda _-) \cap w _{\Theta}(\Lambda _+)$, and for any $\delta \in \Gamma$ with $w_{\Theta}(\delta ) = \delta '$ we have defined a non trivial cut $\Lambda ^{\delta} (\Theta ) = \bigl (\Lambda ^{\delta} _- (\Theta) , \Lambda ^{\delta} _+ (\Theta ) \bigr)$ of the group $\Theta$ by $\Lambda ^{\delta} (\Theta ) = (\Lambda -{\delta}) \cap \Theta$.  

The image $w_{\Theta}(\delta )$ of an element $\delta$ of $\Gamma$ belongs to the intersection $ w_{\Theta}(\Lambda _-) \cap w _{\Theta}(\Lambda _+)$ if and only if $\delta$ belongs to the intersection $ (\Lambda _- + \Theta) \cap (\Lambda _+ + \Theta)$. 

\begin{lemma}\label{le:sous-groupe-invariant-sous-groupe} 
We assume that we have the stict inclusion $\Delta (\Lambda ) \subsetneq \Theta$, then 
all the cuts $\Lambda ^{\delta}(\Theta )$ of $\Theta$ are linearly equivalent for $\delta$ that belongs to $ (\Lambda _- + \Theta) \cap (\Lambda _+ + \Theta)$, and the invariance subgroup $\Delta (\Lambda ^{\delta}(\Theta ))$ of the cut $\Lambda ^{\delta}(\Theta )$ of $\Theta$ is equal to the subgroup $\Delta (\Lambda )$. 
\end{lemma}  

\begin{preuve} 
This a consequence of the fact that the invariance subgroup $\Delta (\Lambda )$ of a cut $\Lambda$ may be defined as the subgroup made up of elements $\xi$ of $\Gamma$ satisfying $\Lambda _- + \xi = \Lambda _-$. 

\hfill$\Box$ 
\end{preuve} 

\vskip .2cm

We consider two convex subgroups $\Theta _1$ and $\Theta _2$ of $\Gamma$ with $\Delta (\Lambda ) \subsetneq \Theta _1 \subset \Theta _2$, and the morphisms $\xymatrix{w_{\Theta _1}: \Gamma \ar[r] & \bar\Gamma _1= \Gamma /\Theta _1}$ and $\xymatrix{w_{\Theta _2}: \Gamma \ar[r] & \bar\Gamma _2= \Gamma /\Theta _2}$. The image of the unique element $\bar\delta _1$ in $w_{\Theta _1} (\Lambda _-) \cap w _{\Theta _1}(\Lambda _+)$ by the morphism $\xymatrix{\bar\Gamma _1 \ar[r] & \bar\Gamma _2}$ is the element $\bar\delta _2$ in $w_{\Theta _2}(\Lambda _-) \cap w _{\Theta _2}(\Lambda _+)$. 
So we can choose $\delta$ in $\Gamma$ whose images by the morphisms $w_{\Theta _1}$ and $w_{\Theta _2}$ are respectively the elements $\bar\delta _1$ and $\bar\delta _2$.

\begin{proposition}\label{prop::existence-d-un-delta}   
If the final segment ${\bf L}_+^{(\Lambda )}$ of ${\bf Cv}(\Gamma )$ has a smallest element, that is, if we have $\Delta (\Lambda) \subsetneq \overline\Delta (\Lambda)$, there exists an element $\delta$ in $\Gamma$ such that for any convex subgroup $\Theta$ in ${\bf L}_+^{(\Lambda )}$, its image $\delta ' = w_{\Theta}(\delta )$ in the quotient group $\Gamma /\Theta$ is the element of  $w_{\Theta} (\Lambda _-) \cap w _{\Theta}(\Lambda _+)$. 
\end{proposition} 

\begin{preuve} 
If the final segment ${\bf L}_+^{(\Lambda )}$ of ${\bf Cv}(\Gamma )$ has a smallest element $\Theta _0$ it is enough to choose an element $\delta$ in $\Gamma$ such that its image in $\Gamma /\Theta _0$ belongs to the intersection $w_{\Theta _0} (\Lambda _-) \cap w _{\Theta _0}(\Lambda _+)$.

\hfill$\Box$ 
\end{preuve} 

\vskip .2cm

\begin{definition}
Let $\Theta _1$ and $\Theta _2$ be a convex subgroups of $\Gamma$, then we can define the following subsets of the set ${\bf Cp}(\Gamma )$ of cuts of $\Gamma$ by the following: 
\begin{enumerate}
\item ${\bf Cp}_{\scriptscriptstyle\Theta _1 \subset .}(\Gamma )$ is the set of cuts $\Lambda$ of $\Gamma$ such that we have the inclusion $\Theta _1 \subset \Delta (\Lambda )$; 

\item ${\bf Cp}_{\scriptscriptstyle . \subsetneq \Theta _2}(\Gamma )$ is the set of cuts $\Lambda$ of $\Gamma$ such that we have the strict inclusion $\Delta (\Lambda ) \subsetneq \Theta _2$;  

\item ${\bf Cp}_{\scriptscriptstyle \Theta _1 \subset . \subsetneq \Theta _2}(\Gamma )$ is the set of cuts $\Lambda$ of $\Gamma$ such that we have the inclusions $\Theta _1 \subset \Delta (\Lambda ) \subsetneq \Theta _2$.
\end{enumerate}
\end{definition}

As two linearly equivalent cuts of $\Gamma$ have the same invariance subgroup, these subsets are stable by the equivalence relation $\simeq$, and we can define the subsets ${\bf Cp}^{[lin]}_{\scriptscriptstyle \Theta _1 \subset .}(\Gamma )$, ${\bf Cp}^{[lin]}_{\scriptscriptstyle . \subsetneq \Theta _2}(\Gamma )$, and ${\bf Cp}^{[lin]}_{\scriptscriptstyle \Theta _1 \subset . \subsetneq \Theta _2}(\Gamma )$ of the set ${\bf Cp}^{[lin]}[\Gamma )$. 

\vskip .2cm 

In the same way we can define the sets ${\bf IS}_{\scriptscriptstyle \Theta _1 \subset .}(\Gamma )$, ${\bf IS}_{\scriptscriptstyle . \subsetneq \Theta _2}(\Gamma )$ and ${\bf IS}_{\scriptscriptstyle \Theta _1 \subset . \subsetneq \Theta _2}(\Gamma )$ of initial segments of $\Gamma$, and the sets ${\bf FS}_{\scriptscriptstyle \Theta _1 \subset .}(\Gamma )$, ${\bf FS}_{\scriptscriptstyle . \subsetneq \Theta _2}(\Gamma )$ and ${\bf FS}_{\scriptscriptstyle \Theta _1 \subset . \subsetneq \Theta _2}(\Gamma )$ of final segments of $\Gamma$

\vskip .2cm

Let $\Lambda$ be a non trivial cut of the ordered group $\Gamma$, let $\Delta = \Delta (\Lambda )$ be its invariance subgroup, for any convex subgroups $\Theta _1$ and $\Theta _2$ of $\Gamma$ with $\Theta _1 \subset \Delta \subsetneq \Theta _2$, the cut $\Lambda$ induces a non trivial cut $\Lambda '$ of the ordered group $\Gamma ' = \Theta _2 / \Theta _1$ in the following way: 
let $\delta$ be an element of $\Gamma$ that belongs to $(\Lambda _- + \Theta _2) \cap (\Lambda _+ + \Theta _2)$, and let $\Lambda ^{\delta}(\Theta _2) = \bigl (\Lambda ^{\delta} _- (\Theta _2) , \Lambda ^{\delta} _+ (\Theta _2) \bigr)$ be the non trivial cut of the group $\Theta _2$ defined above, then the cut $\Lambda '$ is the cut induced on $\Gamma ' = \Theta _2 / \Theta _1$. 
The cut $\Lambda '$ is not uniquely defined, it depends on the element $\delta$ chosen, but all cuts thus defined are linearly equivalent. 
We have then defined a map 
$$\xymatrix @C=10mm 
{T_{\Theta _1}^{\Theta _2} : {\bf Cp}^{[lin]}_{\scriptscriptstyle \Theta _1 \subset . \subsetneq \Theta _2}(\Gamma )) \ar[r] & {\bf Cp}^{[lin]}(\Theta _2 / \Theta _1)}$$
from the set of linear equivalence classes of cuts $\Lambda$ of the group $\Gamma$ with $\Theta _1 \subset \Delta (\Lambda ) \subsetneq \Theta _2$, and the set of linear equivalence classes of non-trivial cuts of the group $\Theta _2 / \Theta _1$. 

\begin{proposition}\label{prop:bijection} 
The application $\xymatrix{T_{\Theta _1}^{\Theta _2} : {\bf Cp}^{[lin]}_{\scriptscriptstyle \Theta _1 \subset . \subsetneq \Theta _2}(\Gamma )) \ar[r] & {\bf Cp}^{[lin]}(\Theta _2 / \Theta _1)}$ is a bijection. 
\end{proposition}

\begin{preuve} 
It is enough to prove that the applications $\xymatrix{T_{\Theta _1} : {\bf Cp}^{[lin]}_{\scriptscriptstyle \Theta _1 \subset .}(\Gamma )) \ar[r] & {\bf Cp}^{[lin]}(\Gamma / \Theta _1)}$ and $\xymatrix{T^{\Theta _2} : {\bf Cp}^{[lin]}_{\scriptscriptstyle . \subsetneq \Theta _2}(\Gamma )) \ar[r] & {\bf Cp}^{[lin]}(\Theta _2)}$ are bijections.

1) The application $T_{\Theta _1}$ is induced by the map $\xymatrix{Is_{\Theta _1} : {\bf IS}_{\scriptscriptstyle \Theta _1 \subset .}(\Gamma )) \ar[r] & {\bf IS}(\Gamma / \Theta _1)}$ which sends the initial segment $\Sigma$ of $\Gamma$ to the initial segment $\bar\Sigma = w_{\Theta _1}(\Sigma )$ of $\Gamma /\Theta _1$, and from proposition \ref{prop:passage-au-quotient} we deduce from the inclusion $\Theta _1 \subset \Delta (\Sigma )$ the equality $w_{\Theta _1}^{-1}(\bar\Sigma ) = \Sigma$. Consequently the application $Is_{\Theta _1}$ is a bijection, hence the application $T_{\Theta _1}$ is also a bijection.  

\vskip .2cm 

2) Let $\Lambda '=(\Lambda '_-,\Lambda '_+)$ be a non-trivial cut of the subgroup $\Theta _2$, then the subsets $\Lambda _-$ and $\Lambda _+$ of $\Gamma$ defined by $\Lambda _- = \{ \xi \in \Gamma \ | \ \xi \leq \Lambda '_- \}$ and $\Lambda _+ =\{ \xi \in \Gamma \ | \ \xi \geq \Lambda '_+ \}$ are respectively initial and final segments of $\Gamma$, and they define a cut $\Lambda = (\Lambda _-,\Lambda _+)$, such that we have the strict inclusion $\Delta (\Lambda ) \subsetneq \Theta _2$. 
Then for any $\zeta$ in $\Gamma$ which belongs to $(\Lambda _- + \Theta _2 ) \cap (\Lambda _+ + \Theta _2)$, the non trivial cut $\Lambda ^{\zeta}(\Theta _2)$ is linearly equivalent to the cut $\Lambda '$. 

Let $\Lambda _a$ and $\Lambda _b$ two cuts of $\Gamma$ with $\Delta (\Lambda _a) \subsetneq \Theta _2$ and $\Lambda _b) \subsetneq \Theta _2$, and let $\zeta _a$ and $\zeta _b$ belonging respectively to the sets $({\Lambda _a} _- + \Theta _2 ) \cap ({\Lambda _a} _+ + \Theta _2)$ and $({\Lambda _b} _- + \Theta _2 ) \cap ({\Lambda _b} _+ + \Theta _2)$. We assume that the cuts ${\Lambda _a} ^{\zeta _a}(\Theta _2)$ and ${\Lambda _b} ^{\zeta _b}(\Theta _2)$ of $\Theta _2$ are linearly equivalent, i.e. there exists $\delta \in \Theta _2$ such that we have the equality ${\Lambda _a} ^{\zeta _a}(\Theta _2) = {\Lambda _b} ^{\zeta _b}(\Theta _2) + \delta$.

Then by replacing the cuts $\Lambda _a$ and $\Lambda _b$ by the linearly equivalent cuts $\Lambda _a - \zeta _a$ and $\Lambda _b - \zeta _b - \delta$, we can assume that the traces of these cuts on $\Theta _2$ are equal non-trivial cuts, then there exist an element $\xi _-$ in ${\Lambda _a}_- \cap \Theta _2 = {\Lambda _b}_- \cap \Theta _2$ and an element $\xi _+$ in ${\Lambda _a}_+ \cap \Theta _2 = {\Lambda _b}_+ \cap \Theta _2$. 
Suppose that the cuts $\Lambda _a$ and $\Lambda _b$ are different and that there exists an element $\xi$ of $\Gamma$ that belongs to ${\Lambda _a}_- \cap {\Lambda _b}_+$, then we would have the inequalities $\xi _- < \xi < \xi _+$. as $\xi _-$ and $\xi _+$ belong to the convex subgroup $\Theta _2$, we deduce that $\xi$ also belongs to $\Theta _2$, which is impossible.

\hfill$\Box$ 
\end{preuve} 

We deduce from remark \ref{rmq:meme-groupe-d-invariance} that for any totally ordered group $\Gamma$ we can define an application $\xymatrix{\Delta _{\Gamma}^{[lin]} :  {\bf Cp}^{[lin]} \ar[r] & \widehat{\bf Cv}(\Gamma )}$, that sends the class of linear equivalence of a cut $\Lambda$ of $\Gamma$ on the subgroup of invariance $\Delta (\Lambda )$ of $\Lambda$. 
Then we have the following result. 

\begin{proposition}\label{prop:sous-groupe-d-invariance-de-l-image} 
Let $\Theta _1$ and $\Theta _2$ be two convex subgroups of $\Gamma$ and let $\Lambda$ be a cut of $\Gamma$ with the inclusions $\Theta _1 \subset \Delta (\Lambda ) \subsetneq \Theta _2$, then we have the equality 
$$\Delta _{\Theta _2 / \Theta _1}^{[lin]}\bigl ( T_{\Theta _1}^{\Theta _2}(\Lambda )\bigr ) = \Delta _{\Gamma}^{[lin]} (\Lambda ) / \Theta _1 \ .$$ 
\end{proposition}

\begin{preuve} 

This is a consequence of lemmas \ref{le:sous-groupe-invariant-quotient} and \ref{le:sous-groupe-invariant-sous-groupe}. 

\hfill$\Box$ 
\end{preuve}

\vskip .2cm

Similarly we can also define bijections

$$\xymatrix @R=1mm @C=10mm
{Is_{\Theta _1}^{\Theta _2} : {\bf IS}^{[lin]}_{\scriptscriptstyle \Theta _1 \subset . \subsetneq \Theta _2}(\Gamma )) \ar[r] & {\bf IS}^{[lin]}(\Theta _2 / \Theta _1) \\  
Fs_{\Theta _1}^{\Theta _2} : {\bf FS}^{[lin]}_{\scriptscriptstyle \Theta _1 \subset . \subsetneq \Theta _2}(\Gamma )) \ar[r] & {\bf FS}^{[lin]}(\Theta _2 / \Theta _1)}$$
between sets of initial segments, or final segments.

\vskip .2cm

\subsection{Cuts associated with a subgroup}

Let $\Delta$ be a proper convex subgroup of $\Gamma$, and let $\zeta$ be an element of $\Gamma$, then we can define as in section 1.1 the cuts $(\zeta +\Delta )^+ = \Lambda ^{\scriptscriptstyle \leq \zeta + \Delta} = \bigl ( \Lambda ^{\scriptscriptstyle \leq \zeta + \Delta}_- , \Lambda ^{\scriptscriptstyle \leq \zeta + \Delta} _+ \bigr )$ and $(\zeta +\Delta )^- = \Lambda ^{\scriptscriptstyle \geq \zeta + \Delta} = \bigl ( \Lambda ^{\scriptscriptstyle \geq \zeta + \Delta}_- , \Lambda ^{\scriptscriptstyle \geq \zeta + \Delta} _+ \bigr )$ of $\Gamma$ by the following: 
$$ \Lambda ^{\scriptscriptstyle \leq \zeta + \Delta}_- = \{ \xi \in \Gamma \ | \ \exists \, \delta \in \Delta \ \hbox{such that} \ \xi - \zeta \leq \delta \}  \ \hbox{and} \  \Lambda ^{\scriptscriptstyle \leq \zeta + \Delta} _+ = \{ \xi \in \Gamma \ | \ \xi - \zeta >\Delta \} \phantom{\ ,}$$ 
$$\Lambda ^{\scriptscriptstyle \geq \zeta + \Delta} _- = \{ \xi \in \Gamma \ | \ \xi - \zeta <\Delta \}  \ \hbox{and}  \ \Lambda ^{\scriptscriptstyle \geq \zeta + \Delta}_+ = \{ \xi \in \Gamma \ | \ \exists \, \delta \in \Delta \ \hbox{such that} \ \xi - \zeta \geq \delta \} \ ,$$
and in the same way we can define the cuts  $\Lambda ^{\scriptscriptstyle < \zeta + \Delta}$ and $\Lambda ^{\scriptscriptstyle > \zeta + \Delta}$. 
As a subgroup $\Delta$ is never empty, and as we assume $\Delta \not= \Gamma$, none of the cuts $\Lambda ^{\scriptscriptstyle > \zeta + \Delta}$, $\Lambda ^{\scriptscriptstyle < \zeta + \Delta}$, $\Lambda ^{\scriptscriptstyle \geq \zeta + \Delta}$ and $\Lambda ^{\scriptscriptstyle \leq \zeta + \Delta}$ are trivial. 

\vskip .2cm 

We can also define the initial segment $ \Lambda ^{\scriptscriptstyle \leq \zeta + \Delta}_-$ and the final segment $\Lambda ^{\scriptscriptstyle \geq \zeta + \Delta}_+$ by the following: 
$$ \Lambda ^{\scriptscriptstyle \leq \zeta + \Delta}_- = \{ \xi \in \Gamma \ | \  w_{\Delta}(\xi ) \leq w_{\Delta}(\zeta ) \}  \ \hbox{and} \  \Lambda ^{\scriptscriptstyle \geq \zeta + \Delta}_+ = \{ \xi \in \Gamma \  | \  w_{\Delta}(\xi ) \geq w_{\Delta}(\zeta ) \} \ ,$$ 
where $w_{\Delta}$ is the natural morphism $\xymatrix{w_{\Delta} : \Gamma \ar[r] & \Gamma / \Delta}$.

\vskip .2cm

\begin{remark}\label{rmq:egalite-des-coupures} 
For any elements $\zeta$ and $\zeta '$ and any convex subgroups $\Delta$ and $\Delta '$ we have the following relations: 
\begin{enumerate} 
\item $\zeta \leq \zeta ' \ \Longrightarrow \ \Lambda ^{\scriptscriptstyle \leq \zeta + \Delta} \leq \Lambda ^{\scriptscriptstyle \leq \zeta ' + \Delta}$ and $ \Lambda ^{\scriptscriptstyle \geq \zeta + \Delta} \leq \Lambda ^{\scriptscriptstyle \geq \zeta ' + \Delta}$;      
\item $\Delta \subset \Delta ' \ \Longrightarrow \ \Lambda ^{\scriptscriptstyle \leq \zeta + \Delta} \leq \Lambda ^{\scriptscriptstyle \leq \zeta + \Delta '}$ and $ \Lambda ^{\scriptscriptstyle \geq \zeta + \Delta '} \leq \Lambda ^{\scriptscriptstyle \geq \zeta + \Delta}$;   
\item $\Lambda ^{\scriptscriptstyle \leq \zeta + \Delta} = \Lambda ^{\scriptscriptstyle \leq \zeta ' + \Delta '} \ \Longleftrightarrow \ \Lambda ^{\scriptscriptstyle \geq \zeta + \Delta} = \Lambda ^{\scriptscriptstyle \geq \zeta ' + \Delta '} \ \Longleftrightarrow \ \Delta = \Delta '$ and $\zeta - \zeta ' \in \Delta\ $.
\end{enumerate}  
\end{remark} 

\vskip .2cm 

If the convex subgroup $\Delta$ is non reduced to $(0)$, $\zeta + \Delta$ has no greatest element and no smallest element, hence we deduce from remark \ref{rmq:<=leq->=geq} the equalities $\Lambda ^{\scriptscriptstyle > \zeta + \Delta} = \Lambda ^{\scriptscriptstyle \geq \zeta + \Delta}$ and $\Lambda ^{\scriptscriptstyle < \zeta + \Delta} = \Lambda ^{\scriptscriptstyle \leq \zeta + \Delta}$. 
And for $\Delta =(0)$ we get the cuts $\Lambda ^{\scriptscriptstyle > \zeta } = \Lambda ^{\scriptscriptstyle \leq \zeta } = \bigl (\Gamma _{\leq \zeta}, \Gamma _{> \zeta} \bigr ) $ and $\Lambda ^{\scriptscriptstyle \geq \zeta } = \Lambda ^{\scriptscriptstyle < \zeta} = \bigl (\Gamma _{< \zeta }, \Gamma _{\geq \zeta} \bigr )$. 

\begin{lemma}\label{le:plus-grand-sous-groupe} 
For any convex subgroup $\Delta$ of $\Gamma$ we have the strict inequality $\Lambda ^{\scriptscriptstyle \geq \zeta + \Delta } < \Lambda ^{\scriptscriptstyle \leq \zeta + \Delta}$ and the subset $\Lambda ^{\scriptscriptstyle \geq \zeta + \Delta} _+ \cap \Lambda ^{\scriptscriptstyle \leq \zeta + \Delta} _-$ is equal to $\zeta + \Delta$. 

Moreover $\Delta$ is the greatest subgroup of $\Gamma$ included in $\Lambda ^{\scriptscriptstyle \geq \Delta} _+$, and is also greatest subgroup of $\Gamma$ included in $\Lambda ^{\scriptscriptstyle \leq \Delta} _-$. 
\end{lemma} 

\begin{preuve} 
The first part of the lemma is obvious. 

Let $\Theta$ be a subgroup of $\Gamma$ included in $\Lambda ^{\scriptscriptstyle \geq \Delta} _+$, then for any $\xi$ in $\Theta$ there exist $\delta _1$ and $\delta _2$ in $\Delta$ such that $\xi \geq \delta _1$ and $- \xi \geq  \delta _2$. 
As $\Delta$ is a convex subgroup of $\Gamma$ we deduce that $\xi$ belongs to $\Delta$. 

\hfill$\Box$ 
\end{preuve} 

\vskip .2cm 

\begin{proposition}\label{prop:sous-groupe-invariant} 
(1) For any convex subgroup $\Delta$ of $\Gamma$ and any element $\zeta$ in $\Gamma$, we have the equalities: 
$$\Delta \ = \ \Delta (\Lambda ^{\scriptscriptstyle \leq \zeta + \Delta}) \ = \ \Delta (\Lambda ^{\scriptscriptstyle \geq \zeta + \Delta}) \ .$$ 

(2) Moreover for any cut $\Lambda$ of $\Gamma$ that satisfies the strict inequalities $\Lambda ^{\scriptscriptstyle \geq \zeta + \Delta } < \Lambda < \Lambda ^{\scriptscriptstyle \leq \zeta + \Delta}$, the convex subgroup $\Delta$ belongs to the final segment ${\bf L}_+ ^{(\Lambda )}$, i.e. we have the strict inclusion 
$$\Delta (\Lambda ) \subsetneq \Delta \ .$$ 
\end{proposition} 

\begin{preuve} 
For any convex subgroup $\Delta$ all the cuts $\Lambda ^{\scriptscriptstyle \leq \zeta + \Delta}$, and all the cuts $\Lambda ^{\scriptscriptstyle \geq \zeta + \Delta}$, when $\zeta$ goes through $\Gamma$, are linearly equivalent. Then we can assume $\zeta =0$. 

(1) The proposition is obvious for $\Delta = (0)$, so we assume in the following that the subgroup $\Delta$ is non reduced to $(0)$. 

We have the equalities  
$\Lambda ^{\scriptscriptstyle \geq \Delta} _- = - (\Lambda ^{\scriptscriptstyle \leq \Delta} _+)$ and 
$\Lambda ^{\scriptscriptstyle \geq \Delta} _+ = - (\Lambda ^{\scriptscriptstyle \leq \Delta} _-)$, we deduce from proposition \ref{prop:image-de-coupure} that the sets ${\bf L} _- ^{(\Lambda ^{\scriptscriptstyle \geq \Delta})}$ and  ${\bf L} _- ^{(\Lambda ^{\scriptscriptstyle \leq \Delta})}$ are equal, hence the equality 
$\Delta (\Lambda ^{\scriptscriptstyle \leq \Delta}) = \Delta (\Lambda ^{\scriptscriptstyle \geq \Delta})$. 

If $\xi$ belongs to the initial segment $\Lambda ^{\scriptscriptstyle \geq\Delta} _-$, for any $\delta$ belonging to the subgroup $\Delta$ it is obvious that $\xi + \delta$ belongs also to $\Lambda ^{\scriptscriptstyle \geq \Delta} _-$, hence the subgroup $\Delta$ is in the subset ${\bf L} _- ^{(\Lambda ^{\scriptscriptstyle \geq \Delta})}$, and we have the inclusion $\Delta \subset \Delta (\Lambda ^{\scriptscriptstyle \geq \Delta})$.  

Conversely let $\Delta '$ be a convex subgroup of $\Gamma$ strictly containing $\Delta$, we want to show that the intersection $( \Lambda ^{\scriptscriptstyle \geq \Delta} _- + \Delta ') \cap \Lambda ^{\scriptscriptstyle \geq \Delta} _+$ is non empty, hence that $\Delta '$ doesn't belong to ${\bf L} _- ^{(\Lambda ^{\scriptscriptstyle \geq \Delta})}$. 

We deduce from lemma \ref{le:plus-grand-sous-groupe} that $\Delta '$ is not included in $\Lambda ^{\scriptscriptstyle \geq \Delta} _+$, then we can choose an element $\zeta$ belonging to $\Lambda ^{\scriptscriptstyle \geq \Delta} _- \cap \Delta '$, and in an analogous way we can find an element $\xi$ belonging to $\Lambda ^{\scriptscriptstyle \leq \Delta}  _+ \cap \Delta '$, hence belonging to $\Lambda ^{\scriptscriptstyle \geq \Delta} _+ \cap \Delta '$.  
Hence the element $\zeta + (\xi - \zeta) = \xi$ belongs to the intersection $( \Lambda ^{\scriptscriptstyle \geq \Delta} _- + \Delta ') \cap \Lambda ^{\scriptscriptstyle \geq \Delta} _+$.   

\vskip .2cm 

(2) From the strict inequalities $\Lambda ^{\scriptscriptstyle \geq \Delta } < \Lambda < \Lambda ^{\scriptscriptstyle \leq \Delta}$ we deduce that there exist $\xi _1$ and $\xi _2$ that belong respectively to $\Lambda ^{\scriptscriptstyle \geq \Delta}_+ \cap \Lambda_-$ and $\Lambda _+ \cap \Lambda ^{\scriptscriptstyle \leq \Delta}_-$, and we get the following relations: 
\begin{enumerate} 
\item \emph{$\xi _1 \in \Lambda _-$ and $\xi _2 \in \Lambda _+$, then $0 < \xi _2 - \xi _1$;}  

\item \emph{$\xi _1 \in \Lambda ^{\scriptscriptstyle \geq \Delta}_+$ then there exists $\delta _1 \in \Delta$ such that $\xi _1 \geq \delta _1$;} 
 
\item \emph{$\xi _2 \in \Lambda ^{\scriptscriptstyle \leq \Delta}_- $ then there exists $\delta _2 \in \Delta$ such that $\xi _2 \leq \delta _2$.}  
\end{enumerate}  

Then we obtain the inequalities $0 < \xi _2 - \xi _1 \leq \delta _2 - \delta _1$, hence $\xi _2 - \xi _1$ belongs to $\Delta$, and $\xi _2 - \xi _1$ doesn't belong to $\Delta (\Lambda )$ because we have the equality $(\xi _2 - \xi _1 )+\xi _1 = \xi _2$ with $\xi _1 \in \Lambda _-$ and $\xi _2 \in \Lambda _+$.  

\hfill$\Box$ 
\end{preuve}

\vskip .2cm 

\subsection{Symmetric interval associated with a cut} 

Let $\Sigma$ be an initial segment of $\Gamma$, then for any $\sigma \in \Sigma$ we define the new initial segment $\Sigma ^{(\sigma )}$ by the following: 
$$\Sigma ^{(\sigma )} = \Sigma -{\sigma} = \{ \xi \in \Gamma \ | \ \xi + \sigma \in \Sigma \} \ ,$$
the element $0$ always belongs to $\Sigma ^{(\sigma )}$, we define the subset $T _{(\Sigma ,\sigma )}$ of $\Gamma$ by 
$$T _{(\Sigma ,\sigma )} = \bigl ( \Sigma ^{(\sigma )} \cap \Gamma _{\geq 0} \bigr ) = \{ \xi \in \Gamma _{\geq 0} \ | \ \xi + \sigma \in \Sigma \} \ ,$$ 
and we define the symmetric interval $S _{(\Sigma ,\sigma )}$ of $\Gamma$ by  
$$S _{(\Sigma ,\sigma )} = T _{(\Sigma ,\sigma )} \cup -(T _{(\Sigma ,\sigma )}) \ .$$ 
By construction $S_{(\Sigma ,\sigma )}$ is the greatest symmetric interval of $\Gamma$ which is contained in $\Sigma ^{(\sigma )}$.  

We have associated with the symmetric interval $S _{(\Sigma ,\sigma )}$ a cut ${\bf L}^{(S _{(\Sigma ,\sigma )})}$ of the set ${\bf Cv}(\Gamma )$ of convex subgroups of $\Gamma$, that we denote ${\bf C} ^{(\Sigma ,\sigma )}=\bigl ({\bf C}_- ^{(\Sigma ,\sigma )}, {\bf C}_+^{(\Sigma ,\sigma )}\bigr )$, by the following (c.f. proposition \ref{prop:coupure-intervalle-symetrique}):
$$ \Delta \in  {\bf C}_- ^{(\Sigma ,\sigma )} \ \Longleftrightarrow \ \Delta \subset S _{(\Sigma ,\sigma )}  \ \Longleftrightarrow \ \Delta \subset \Sigma ^{(\sigma )}  \ \Longleftrightarrow \  \sigma + \Delta \subset \Sigma  \ .$$
If the initial segment $\Sigma$ is non trivial, the cut ${\bf C} ^{(\Sigma ,\sigma )}$ of ${\bf Cv}(\Gamma )$ is non trivial, and we denote by $\Phi ^-(\Sigma ,\sigma )$ ans $\Phi ^+(\Sigma, \sigma )$ the two convex subgroups of $\Gamma$ defined in proposition \ref{prop:principal-convex-subgroup}:  
$$\Phi ^-(\Sigma ,\sigma ) = \Delta _{{\bf C}^{({\Sigma ,\sigma })}}^{-} = \bigcup _{\Delta \in {\bf C}_- ^{(\Sigma ,\sigma )}} \Delta \quad\hbox{and}\quad \Phi ^+(\Sigma ,\sigma ) = \Delta _{{\bf C}^{({\Sigma ,\sigma })}}^{+} = \bigcap _{\Delta \in {\bf C}_+ ^{(\Sigma ,\sigma )}} \Delta \ .$$  
By construction $\Phi ^-(\Sigma ,\sigma )$ belongs to ${\bf C}_- ^{(\Sigma ,\sigma )}$, it is the greatest convex subgroup of $\Gamma$ which is contained in $\Sigma ^{(\sigma )}$, and we have the equivalence $\Delta \in {\bf C}_- ^{(\Sigma ,\sigma )} \Longleftrightarrow \Delta \subset \Phi ^-(\Sigma ,\sigma )$. 
Moreover we deduce from lemma \ref{le:phi-et-delta+}  and remark \ref{rmq:initial-segment-convex-subgroup} that we have the equality 
$$\Phi ^-(\Sigma ,\sigma ) = (S _{\Sigma ,\sigma )}) _{\rm tor} =  \varprojlim (1/n) S _{(\Sigma ,\sigma )} \ .$$

\vskip .2cm 

From remark \ref{rmq:contientT} we can define another cut ${\bf D} ^{(\Sigma ,\sigma )}=\Bigl ({\bf D}_- ^{(\Sigma ,\sigma )}, {\bf D}_+^{(\Sigma ,\sigma )}\Bigr )$ of the set ${\bf Cv}(\Gamma )$ by 
$${\bf D}_+^{(\Sigma ,\sigma )} = {\bf Cv}_{T _{(\Sigma ,\sigma )}} (\Gamma ) = \{ \Delta \in {\bf Cv}(\Gamma ) \ | \ T _{_{(\Sigma ,\sigma )}}  \subset \Delta \} = \{ \Delta \in {\bf Cv}(\Gamma ) \ | \ S _{_{(\Sigma ,\sigma )}}  \subset \Delta \} \ ,$$ 
and we denote by $\Psi ^-(\Sigma ,\sigma )$ ans $\Psi ^+(\Sigma, \sigma )$ the two convex subgroups of $\Gamma$ defined in proposition \ref{prop:principal-convex-subgroup}: 
$$\Psi ^-(\Sigma ,\sigma ) = \Delta _{{\bf D}^{({\Sigma ,\sigma })}}^{-} = \bigcup _{\Delta \in {\bf D}_- ^{(\Sigma ,\sigma )}} \Delta \quad\hbox{and}\quad \Psi ^+(\Sigma ,\sigma ) = \Delta _{{\bf D}^{({\Sigma ,\sigma })}}^{+} = \bigcap _{\Delta \in {\bf D}_+ ^{(\Sigma ,\sigma )}} \Delta \ .$$ 
By construction $\Psi ^+(\Sigma, \sigma )$ belongs to the final segment ${\bf D}_+^{(\Sigma ,\sigma )}$, and is the smallest convex subgroup of $\Gamma$ which contains $S _{(\Sigma ,\sigma )}$, moreover we have the inclusion $\Psi ^-(\Sigma ,\sigma ) \subset \Psi ^+(\Sigma ,\sigma )$ with strict inequality if and only if $\Psi ^+(\Sigma, \sigma )$ is equal to a principal convex subgroup $\Delta ^+_{\{\xi\}}$ for some element $\xi$. For any convex subgroup $\Delta$ we have the equivalences: 

$$ \Delta \in  {\bf D}_+ ^{(\Sigma ,\sigma )} \ \Longleftrightarrow \  S _{(\Sigma ,\sigma )} \subset \Delta  \ \Longleftrightarrow \ \Psi ^+(\Sigma ,\sigma ) \subset \Delta   \  .$$

\vskip .2cm 

\begin{proposition}\label{prop:comparaison-C-D} 
Let $\Sigma$ be an initial segment and let $\sigma$ be an element in $\Sigma$, then the two cuts ${\bf C} ^{(\Sigma ,\sigma )}$ and ${\bf D} ^{(\Sigma ,\sigma )}$ of ${\bf Cv}(\Gamma )$ defined above satisfy the inequality ${\bf D} ^{(\Sigma ,\sigma )}  \leq {\bf C} ^{(\Sigma ,\sigma )}$, and we have a strict inequality if and only if the symmetric interval $S _{(\Sigma ,\sigma )}$ is a convex subgroup of $\Gamma$.  
\end{proposition}

\begin{preuve} 
The two cuts ${\bf C} ^{(\Sigma ,\sigma )}$ and ${\bf D} ^{(\Sigma ,\sigma )}$ are defined by 
$${\bf C} ^{(\Sigma ,\sigma )} _- = \{ \Delta  \in {\bf Cv}(\Gamma ) \ | \ \Delta  \subset S _{(\Sigma ,\sigma )} \} \quad\hbox{and}\quad
{\bf D} ^{(\Sigma ,\sigma )}_+ = \{ \Delta  \in {\bf Cv}(\Gamma ) \ | \ S _{(\Sigma ,\sigma )} \subset \Delta \}  \ , $$
then we have the equalities ${\bf D} ^{(\Sigma ,\sigma )} _+ \cap {\bf C} ^{(\Sigma ,\sigma )} _- = \{ \Delta  \in {\bf Cv}(\Gamma ) \ | \ S _{(\Sigma ,\sigma )} = \Delta  \}$ and ${\bf D} ^{(\Sigma ,\sigma )} _- \cap {\bf C} ^{(\Sigma ,\sigma )} _+ = \emptyset$, hence the desired result.  

\hfill$\Box$ 
\end{preuve} 

\vskip .2cm 

\begin{remark}\label{rmq:Delta=Phi}
If the two cuts ${\bf C} ^{(\Sigma ,\sigma )}$ and ${\bf D} ^{(\Sigma ,\sigma )}$ are not equal, the convex subgroup $\Delta$ of $\Gamma$ such that we have the equality $S _{(\Sigma ,\sigma )} = \Delta $ is the subgroup $\Phi ^- (\Sigma ,\sigma )$. 
\end{remark} 

\vskip .2cm 

\begin{lemma}\label{le:C-D}
Let ${\bf C}_-$ and ${\bf D}_+$ be respectively an initial segment and a final segment of a totally ordered set $I$ such that we have ${\bf D}_+ \cap {\bf C}_- = \{ x \}$ for some element $x$ in $I$. 
Then we have  
$${\bf D}_+ = I_{\geq x} \quad \hbox{and} \quad {\bf C}_- = I_{\leq x} \ .$$
\end{lemma} 

\begin{preuve} 
By symmetry it is enough to show ${\bf D}_+ = I_{\geq x}$. 

(i) $y \geq x \ \Longrightarrow \  y \in {\bf D}_+$, because $x \in {\bf D}_+$ and ${\bf D}_+$ is a final segment. 

(ii) $y \in {\bf D}_+$ and $y \not= x \ \Longrightarrow \  y \notin {\bf C}_-  \ \Longrightarrow \  y  > x$, because $x \in {\bf C}_-$ and ${\bf C}_-$ is an inital segment. 

\hfill$\Box$ 
\end{preuve} 

\vskip .2cm 

\begin{corollary}\label{cor:Phi-Psi} 
(1) If the two cuts ${\bf C} ^{(\Sigma ,\sigma )}$ and ${\bf D} ^{(\Sigma ,\sigma )}$ are equal, we have the following relations: 
$$\Psi ^-(\Sigma ,\sigma ) = \Phi ^-(\Sigma ,\sigma ) \subset \Psi ^+(\Sigma ,\sigma ) = \Phi ^+(\Sigma ,\sigma ) \ .$$ 

(2) If the two cuts ${\bf C} ^{(\Sigma ,\sigma )}$ and ${\bf D} ^{(\Sigma ,\sigma )}$ are not equal, we have the following relations: 
$$\Psi ^-(\Sigma ,\sigma ) \subset \Phi ^-(\Sigma ,\sigma ) = \Psi ^+(\Sigma ,\sigma ) \subset \Phi ^+(\Sigma ,\sigma ) \ .$$  
\end{corollary} 

\begin{preuve} 
We deduce from proposition \ref{prop:comparaison-C-D} that we have the inclusions $\Psi ^-(\Sigma ,\sigma ) \subset \Phi ^-(\Sigma ,\sigma )$ and $ \Psi ^+(\Sigma ,\sigma ) \subset \Phi ^+(\Sigma ,\sigma )$. 
We have equalities if the two cuts ${\bf C} ^{(\Sigma ,\sigma )}$ and ${\bf D} ^{(\Sigma ,\sigma )}$ are equal. 

In the case where the two cuts ${\bf C} ^{(\Sigma ,\sigma )}$ and ${\bf D} ^{(\Sigma ,\sigma )}$ are not equal, the result is a consequence of lemma \ref{le:C-D}. 

\hfill$\Box$ 
\end{preuve}

\vskip .2cm

We consider now the initial segment associated with the cut $\Lambda ^{\scriptscriptstyle \leq \zeta + \Delta}$:
$$\Lambda ^{\scriptscriptstyle \leq \zeta + \Delta} _- =  \{ \xi \in \Gamma \ | \ \exists \, \delta \in \Delta \ \hbox{such that} \ \xi - \zeta \leq \delta \} \ . $$ 

\vskip .2cm 

\begin{lemma}\label{le:lambda=phi} 
Let $\Sigma$ be an initial segment and $\sigma$ be an element of  $\Sigma$, then for any convex subgroup $\Delta$ we have the equivalences: 

$$\begin{array}{lrl}
1)\ \Lambda ^{\scriptscriptstyle \leq \sigma +\Delta} _- \subset \Sigma  &  \Longleftrightarrow &  \Delta \in {\bf C} _- ^{(\Sigma ,\sigma )} \\ 
2)\ \Sigma \subset \Lambda ^{\scriptscriptstyle \leq \sigma +\Delta} _- &  \Longleftrightarrow & \Delta \in {\bf D} _+ ^{(\Sigma ,\sigma )} \ .
\end{array}$$
\end{lemma}

\begin{preuve} 
1) It is obvious from the equivalence $\Delta \in {\bf C} _- ^{(\Sigma ,\sigma )} \Longleftrightarrow \sigma + \Delta \subset \Sigma$ and from the equality $\Lambda ^{\scriptscriptstyle \leq \sigma + \Delta} _- =  \{ \xi \in \Gamma \ | \ \exists \, \delta \in \Delta \ \hbox{such that} \ \xi \leq \sigma + \delta \}$.    

2) We recall that we have the equivalence $ \Delta \in  {\bf D}_+ ^{(\Sigma ,\sigma )} \Longleftrightarrow  T_{(\Sigma ,\sigma )} \subset \Delta$ where $T _{(\Sigma ,\sigma )}$ is defined by $T _{(\Sigma ,\sigma )}= \{ \xi \in \Gamma _{\geq 0} \ | \ \xi + \sigma \in \Sigma \}$.

First we assume the inclusion $\Sigma  \subset \Lambda ^{\scriptscriptstyle \leq \sigma +\Delta} _-$, then we have:  
$$\begin{array}{rcl}
\xi \in T_{(\Sigma ,\sigma )} & \Longrightarrow & \xi \geq 0 \ \hbox{and} \ \xi + \sigma \in \Sigma \\ 
& \Longrightarrow & \xi \geq 0 \ \hbox{and} \ \xi + \sigma \in \Lambda ^{\scriptscriptstyle \leq \sigma +\Delta} _- \\ 
& \Longrightarrow & \xi \geq 0 \ \hbox{and} \ \ \exists \, \delta \in \Delta \ \hbox{such that} \ \xi + \sigma \leq \delta + \sigma \\ 
& \Longrightarrow & \exists \, \delta \in \Delta \ \hbox{such that} \  0 \leq \xi \leq \delta \\ 
& \Longrightarrow & \xi \in \Delta \ .
\end{array}$$

Conversely we assume the inclusion $T_{(\Sigma ,\sigma )} \subset \Delta$, it is obvious that  we have $\zeta \in \Lambda ^{\scriptscriptstyle \leq \sigma +\Delta} _-$ for any $\zeta \leq \sigma$, then we can assume the inequality $\zeta \geq \sigma$ and we have: 

$$\begin{array}{rcl}
\zeta \in \Sigma \ \hbox{and} \ \zeta \geq \sigma & \Longrightarrow & \zeta - \sigma  \in T_{(\Sigma ,\sigma )}   \\ 
& \Longrightarrow & \zeta - \sigma \in \Delta \\ 
& \Longrightarrow & \zeta \in \Lambda ^{\scriptscriptstyle \leq \sigma +\Delta} _-  \ .
\end{array}$$

\hfill$\Box$ 
\end{preuve} 

\vskip .2cm 

\begin{proposition}\label{prop:comparaison-C-D_2}
Let $\Sigma$ be an initial segment and let $\sigma$ be an element in $\Sigma$, then the two cuts ${\bf C} ^{(\Sigma ,\sigma )}$ and ${\bf D} ^{(\Sigma ,\sigma )}$ of ${\bf Cv}(\Gamma )$ are different if and only if  the initial segment $\Sigma$ is equal to $\Lambda ^{\scriptscriptstyle \leq \sigma +\Delta} _-$ for some convex subgroup $\Delta$ of $\Gamma$.  
\end{proposition} 

\begin{preuve} 
This is a consequence of proposition \ref{prop:comparaison-C-D}, remark \ref{rmq:Delta=Phi}, and lemma \ref{le:lambda=phi}.

\hfill$\Box$ 
\end{preuve} 

\vskip .2cm

Let $\sigma '$ be an other element of $\Sigma$ with $\sigma ' \geq \sigma$, then we have the inclusions $\Sigma ^{(\sigma ')} \subset \Sigma ^{(\sigma )}$ and $S _{(\Sigma ,\sigma ')} \subset S _{(\Sigma ,\sigma )}$. 
We deduce from the definitions of ${\bf C} ^{(\Sigma ,\sigma )} _-$ and of ${\bf D} ^{(\Sigma ,\sigma )} _+$ that we have the inclusions ${\bf C} ^{(\Sigma ,\sigma ')} _- \subset {\bf C} ^{(\Sigma ,\sigma )} _-$ and ${\bf D} ^{(\Sigma ,\sigma )} _+  \subset {\bf D} ^{(\Sigma ,\sigma ')} _+ $, hence the inequalities 
${\bf C} ^{(\Sigma ,\sigma ')} \leq {\bf C} ^{(\Sigma ,\sigma )}$ and ${\bf D} ^{(\Sigma ,\sigma ')} \leq {\bf D} ^{(\Sigma ,\sigma )}$, 
and we get the inclusions
\begin{enumerate} 
\item $\Psi ^-(\Sigma ,\sigma ') \subset \Psi ^-(\Sigma ,\sigma ) \, , \  
\Psi ^+(\Sigma ,\sigma ') \subset \Psi ^+(\Sigma ,\sigma )$
\item $\Phi ^-(\Sigma ,\sigma ') \subset \Phi ^-(\Sigma ,\sigma ) \, , \ 
\Phi ^+(\Sigma ,\sigma ') \subset \Phi ^+(\Sigma ,\sigma )$ .
\end{enumerate} 
We define the convex subgroups $\Psi ^-(\Sigma )$, $\Psi ^+(\Sigma )$, $\Phi ^-(\Sigma )$, and $\Phi ^+(\Sigma )$ by the following: 
\begin{enumerate}
\item \emph{$\displaystyle{\Psi ^-(\Sigma ) = \bigcap _{\sigma \in \Sigma} \Psi ^-(\Sigma ,\sigma ) \quad \hbox{and} \quad  
\Psi ^+(\Sigma )= \bigcap _{\sigma \in \Sigma} \Psi ^+(\Sigma ,\sigma )}$;} 
\item \emph{$\displaystyle{\Phi ^-(\Sigma ) = \bigcap _{\sigma \in \Sigma} \Phi ^-(\Sigma ,\sigma ) \quad \hbox{and} \quad  
\Phi ^+(\Sigma )= \bigcap _{\sigma \in \Sigma} \Phi ^+(\Sigma ,\sigma )}$.}  
\end{enumerate} 

\vskip .2cm

\begin{remark}\label{rmq:segments-equivalents} 
If $\Sigma '$ is an initial segment equivalent to $\Sigma$, i.e. if we have $\Sigma '=\zeta + \Sigma$ for an element $\zeta \in \Gamma$, we have the equalities $\Psi ^-(\Sigma ,\sigma ) = \Psi ^-(\Sigma ', \zeta + \sigma )$, $\Psi ^+(\Sigma ,\sigma ) = \Psi ^+(\Sigma ', \zeta + \sigma )$, $\Phi ^-(\Sigma ,\sigma ) = \Phi ^-(\Sigma ', \zeta + \sigma )$, and $\Phi ^+(\Sigma ,\sigma ) = \Phi ^+(\Sigma ', \zeta + \sigma )$ for any $\sigma \in \Sigma$, hence the equalities $\Psi ^-(\Sigma ) = \Psi ^-(\Sigma ')$, $\Psi ^+(\Sigma ) = \Psi ^+(\Sigma ')$, $\Phi ^-(\Sigma ) = \Phi ^-(\Sigma ')$, and $\Phi ^+(\Sigma ) = \Phi ^+(\Sigma ')$.  
\end{remark} 

\vskip .2cm 

\begin{proposition}\label{prop:delta=phi} 
For any initial segment $\Sigma$ of the ordered group $\Gamma$, we have the equality 
$$\Delta (\Sigma ) = \Phi ^-(\Sigma ) \ .$$  
\end{proposition}

\begin{preuve} 
For any $\sigma \in \Sigma$, the invariance subgroup $\Delta (\Sigma )$ of  $\Sigma$ satisfies $\sigma + \Delta (\Sigma ) \subset \Sigma$, we deduce that we have the inclusion $\Delta (\Sigma ) \subset \Phi ^- (\Sigma ,\sigma)$ for any $\sigma$, hence the inclusion $\Delta (\Sigma ) \subset \Phi^- (\Sigma )$.   

Conversely, let $\zeta$ be a positive element of $\Gamma$ that doesn't belong to the invariance subgroup $\Delta (\Sigma )$, then we have $\zeta + \Sigma \not= \Sigma$, and as $\zeta > 0$, there exists $\sigma \in \Sigma$ with $\zeta + \sigma \notin \Sigma$, hence we have $\zeta \notin \Phi ^-(\Sigma ,\sigma )$.  

\hfill$\Box$ 
\end{preuve} 

\vskip .2cm

We deduce from lemma \ref{le:lambda=phi} that $\Phi ^-(\Sigma ,\sigma )$ is the greatest convex subgroup $\Phi$ of $\Gamma$ such that we have the inclusion $\Lambda ^{\scriptscriptstyle \leq \sigma +\Phi} _- \subset \Sigma$. 
We denote by $\Lambda ^{\scriptscriptstyle (\Sigma , \sigma )} _-$ the initial segment $\Lambda ^{\scriptscriptstyle \leq \sigma +\Phi ^-(\Sigma ,\sigma )} _-$, and we define the cut $\Lambda ^{\scriptscriptstyle (\Sigma , \sigma )}$ of $\Gamma$ as the cut associated with this initial segment, i.e. we have $\Lambda ^{\scriptscriptstyle (\Sigma , \sigma )} =  \Lambda ^{\scriptscriptstyle \leq \sigma +\Phi ^-(\Sigma ,\sigma )}$.  

Then for any $\sigma \in \Sigma$, by remark \ref{rmq:egalite-des-coupures}, we have the inclusions 
$\Lambda _- ^{\scriptscriptstyle \leq \sigma + \Delta (\Sigma )} \subset \Lambda _- ^{(\Sigma ,\sigma )} \subset \Sigma$, and the equalities 
$$\Sigma = \bigcup _{\sigma \in \Sigma} \Lambda _- ^{\scriptscriptstyle \leq \sigma + \Delta (\Sigma )} = \bigcup _{\sigma \in \Sigma} \Lambda _- ^{(\Sigma ,\sigma )} \ .$$

\vskip .2cm 

\begin{proposition}\label{prop:T=Delta} 
Let $\Delta$ be a convex subgroup of $\Gamma$ and let $\zeta$ be an element of $\Gamma$, then the initial segment $\Sigma = \Lambda _- ^{\scriptscriptstyle \leq \zeta + \Delta}$ satisfies the equality $S_{(\Sigma ,\zeta )} = \Delta$. 

\end{proposition} 

\begin{preuve} 
From the definition of the initial segment $\Lambda _- ^{\scriptscriptstyle \leq \zeta + \Delta}$ we deduce that $T_{(\Sigma ,\zeta )}$ is the subset of the elements $\xi$ in $\Gamma$ such that there exists $\delta \in \Delta$ with $0 \leq \xi \leq \delta$. 
As $\Delta$ is a convex subgroup of $\Gamma$ we obtain the equality $T_{(\Sigma ,\zeta )} = \Delta _{\geq 0}$. 

\hfill$\Box$ 
\end{preuve} 

\vskip .2cm

\begin{corollary}\label{cor:Phi=Delta} 
For any $\sigma$ in $\Sigma = \Lambda _- ^{\scriptscriptstyle \leq \zeta + \Delta}$ with $\sigma \geq \zeta$ the two cuts ${\bf C}^{(\Sigma ,\sigma )}$ and ${\bf D}^{(\Sigma ,\sigma )}$ are not equal, and we have the equality $\Phi ^- (\Sigma ,\sigma ) = \Delta$. 
\end{corollary}

\begin{preuve} 
From proposition \ref{prop:T=Delta} the result is obvious for $\sigma = \zeta$, and it is enough to see that any $\sigma$ in $\Sigma$ with $\sigma \geq \zeta$ satisfies $\sigma - \zeta \in \Delta$, then the initial segment $\Sigma$ is also equal to $\Lambda _- ^{\scriptscriptstyle \leq \sigma + \Delta}$. 

\hfill$\Box$ 
\end{preuve} 

\vskip .2cm

From proposition \ref{prop:delta=phi} and corollary \ref{cor:Phi=Delta} we recover that the invariance subgroup $\Delta (\Lambda _- ^{\scriptscriptstyle \leq \zeta + \Delta})$ is equal to $\Delta$ (cf. proposition \ref{prop:sous-groupe-invariant}).

\vskip .2cm

\begin{proposition}\label{prop:Lambda-notequal-Sigma} 
For any $\sigma$ in $\Sigma$ we have the equivalence: 
$$\Lambda ^{\scriptscriptstyle (\Sigma , \sigma )} _- \subsetneq \Sigma \ \Longleftrightarrow \ {\bf D} ^{(\Sigma ,\sigma )} = {\bf C} ^{(\Sigma ,\sigma )} \ .$$ 
\end{proposition} 

\begin{preuve} 
The initial segments ${\bf D} ^{(\Sigma ,\sigma )}$ and ${\bf C} ^{(\Sigma ,\sigma )}$ are equal if and only if there exists no convex subgroup $\Delta$ of $\Gamma$ such that we have $S_{(\Sigma ,\sigma )} = \Delta$, i.e. if and only if we have $\Phi ^- (\Sigma ,\sigma ) \subsetneq S_{(\Sigma ,\sigma )}$. 
Then the proposition is a consequence of the fact that $\Sigma \cap \Lambda ^{\scriptscriptstyle (\Sigma , \sigma )} _+$ is the set of the elements $\xi$ of $\Gamma$ such that  $\varphi = \xi -\sigma$ satisfies $\varphi \in \Sigma ^{(\sigma )} \cap \Gamma _{\geq 0} = T_{(\Sigma ,\sigma )}$ 
and $\varphi \notin \Phi ^- (\Sigma ,\sigma)$.   

\hfill$\Box$ 
\end{preuve} 

\vskip .2cm 

\subsection{Types of cuts} 

We now want to use the results from the previous section to describe the different types of initial segments of a totally ordered group $\Gamma$, and in the same way the different types of final segments and of cuts of $\Gamma$. 

\begin{theorem}\label{th:segment-initial}
Let $\Sigma$ be an initial segment of $\Gamma$ and let $\Delta (\Sigma )$ be its invariance subgroup. Then we are in one of the following cases: 
\begin{enumerate} 
\item There exists $\sigma \in \Sigma$ which satisfies the equality $\Phi ^- (\Sigma ,\sigma ) =  S_{(\Sigma ,\sigma )}$, then the initial segment $\Sigma$ is equal to $\Lambda _- ^{\scriptscriptstyle \leq \sigma + \Phi ^- (\Sigma ,\sigma )}$, and we have the equality $\Phi ^- (\Sigma ,\sigma ) = \Delta (\Sigma )$.    
\item There exists $\sigma \in \Sigma$ which satisfies the equality $\Phi ^- (\Sigma ,\sigma ) = \Delta (\Sigma )$, but for any $\zeta \in \Sigma$ we have a strict inclusion $\Phi ^- (\Sigma ,\zeta ) \subsetneq S_{(\Sigma ,\zeta )}$, then the invariance subgroup $\Delta (\Sigma )$ is an immediate predecessor in ${\bf Cv}(\Gamma )$.   
\item For any $\sigma$ in $\Sigma$ the subgroup $\Phi ^- (\Sigma ,\sigma )$ is different from $\Delta (\Sigma )$, then the invariance subgroup $\Delta (\Sigma )$ is not an immediate predecessor in ${\bf Cv}(\Gamma )$. 
\end{enumerate} 
\end{theorem}  

\begin{preuve} 
(1) If we have the equality $\Phi ^- (\Sigma ,\sigma ) =  S_{(\Sigma ,\sigma )}$, the two cuts ${\bf D} ^{(\Sigma ,\sigma )}$ and ${\bf C} ^{(\Sigma ,\sigma )}$ are different, then the result is a consequence of proposition \ref{prop:Lambda-notequal-Sigma}.  

(2) If we have the strict inclusion $\Phi ^- (\Sigma ,\sigma ) \subsetneq S_{(\Sigma ,\sigma )}$, the symmetric interval $S_{(\Sigma ,\sigma )}$ is not a subgroup and we deduce from proposition \ref{prop:coupure-intervalle-symetrique} that the cut ${\bf C} ^{(\Sigma ,\sigma )} = {\bf L}^{(S _{(\Sigma ,\sigma )})}$ of ${\bf Cv}(\Gamma )$ is a jump and that the group $\Phi ^-(\Sigma ,\sigma ) = \Delta _{{\bf C}^{({\Sigma ,\sigma })}}^-$ is an immediate predecessor. 

(3) It is a consequence of proposition \ref{prop:delta=phi} which states $\Delta (\Sigma ) = \Phi ^-(\Sigma )  = \bigcap _{\sigma \in \Sigma} \Phi ^-(\Sigma ,\sigma )$. 

\hfill$\Box$ 
\end{preuve} 

\vskip .2cm  

\begin{remark} 
If we are in the case where any subgroup $\Phi ^- (\Sigma ,\sigma )$ is different from $\Delta (\Sigma )$, we deduce from part (1) of the theorem \ref{th:segment-initial} that we have for any $\sigma$ in $\Sigma$ the strict inclusion $\Phi ^- (\Sigma ,\sigma ) \subsetneq S_{(\Sigma ,\sigma )}$, hence from proposition \ref{prop:coupure-intervalle-symetrique} the subgroup $\Phi ^- (\Sigma ,\sigma )$ is an immediate predecessor and the subgroup $\Phi ^+ (\Sigma ,\sigma )$ is its immediate successor and is principal. 

Then the invariance subgroup $\Delta (\Sigma )$ is equal to $\bigcap _{\sigma \in \Sigma} \Phi ^-(\Sigma ,\sigma )$, where each $\Phi ^-(\Sigma ,\sigma )$ is an immediate successor in ${\bf Cv}(\Gamma )$, and it is also equal to $\bigcap _{\sigma \in \Sigma} \Phi ^+(\Sigma ,\sigma )$ where each $\Phi ^+(\Sigma ,\sigma )$ is a principal convex subgroup.  
\end{remark}

\vskip .2cm

For a final segment $\Omega$ and an element $\omega \in \Omega$ we have analogous definitions of the symmetric interval $S_{(\Omega ,\omega )}$, of the cut ${\bf C}^{(\Omega ,\omega )}$ of the ${\bf Cv}(\Gamma )$, and of the subgroup $\Phi ^-(\Omega, \omega )$: 
\begin{enumerate}
\item \emph{the symmetric interval $S_{(\Omega ,\omega )}$ is the greatest symmetric interval contained in $\Omega ^{(\omega )} = \Omega - \omega = \{ \zeta \ | \ \zeta + \omega \in \Omega \}$;}   
\item \emph{the cut ${\bf C}^{(\Omega ,\omega )} = \bigr ({\bf C}^{(\Omega ,\omega )} _- , {\bf C}^{(\Omega ,\omega )} _+\bigl )$ is the cut ${\bf L}^{(S_{(\Omega ,\omega )})}$ associated with the symmetric interval $S_{(\Omega ,\omega )}$, it is defined by 
${\bf C}^{(\Omega ,\omega )} _- = \{ \Delta \ | \ \Delta \subset S_{(\Omega ,\omega )} \}$ and ${\bf C}^{(\Omega ,\omega )} _+ = \{ \Delta \ | \ S_{(\Omega ,\omega )} \subsetneq \Delta \}$;} 
\item \emph{the subgroup $\Phi ^-(\Omega, \omega )$ is defined by 
$ \Phi^-(\Omega, \omega ) = \Delta _{{\bf C}^{({\Omega ,\omega })}}^{-} = \bigcup _{\Delta \in {\bf C}_- ^{(\Omega ,\omega )}} \Delta$,  
it is the greatest convex subroup contained in the symmetric interval $S_{(\Omega ,\omega )}$.}   
\end{enumerate} 

\vskip .2cm 

We then obtain a result analogous to theorm \ref{th:segment-initial}. 

\begin{theorem}\label{th:segment-final}
Let $\Omega$ be a final segment of $\Gamma$ and let $\Delta (\Omega )$ be its invariance subgroup. Then we are in one of the following cases: 
\begin{enumerate} 
\item There exists $\omega \in \Omega$ which satisfies the equality $\Phi ^- (\Omega ,\omega ) =  S_{(\Omega ,\omega )}$, then the final segment $\Omega$ is equal to $\Lambda _+ ^{\geq \omega + \Phi ^- (\Omega ,\omega )}$, and we have the equality $\Phi ^- (\Omega ,\omega ) = \Delta (\Omega )$.    
\item There exists $\omega \in \Omega$ which satisfies the equality $\Phi ^- (\Omega ,\omega ) = \Delta (\Omega )$, but for any $\zeta \in \Omega$ we have a strict inclusion $\Phi ^- (\Omega ,\zeta ) \subsetneq S_{(\Omega ,\zeta )}$, then the invariance subgroup $\Delta (\Omega )$ is an immediate predecessor in ${\bf Cv}(\Gamma )$.   
\item For any $\omega$ in $\Omega$ the subgroup $\Phi ^- (\Omega ,\omega )$ is different from $\Delta (\Omega )$, then the invariance subgroup $\Delta (\Omega )$ is not an immediate predecessor in ${\bf Cv}(\Gamma )$. 
\end{enumerate} 
\end{theorem}

\begin{preuve} 
It is enough to consider the initial segment $\Sigma$ defined by $\Sigma = - \Omega = \{ \sigma \ | \ - \sigma \in \Omega \}$, then for any $\omega \in \Omega$ the symmetric interval $S_{(\Omega ,\omega )}$ is equal to the symmetric interval$S_{(\Sigma , -\omega )}$, and the theorem is a consequence of theorem \ref{th:segment-initial}.  

\hfill$\Box$ 
\end{preuve} 

\vskip .2cm 

\begin{definition} 
Let $\Sigma$ be an initial segment of the ordered group $\Gamma$, with invariance subgroup $\Delta$. 
\begin{enumerate} 
\item If there exists $\sigma$ in $\Sigma$ such that we have the equality $\Sigma = \Lambda _- ^{\scriptscriptstyle \leq \sigma + \Delta}$, we say that the initial segment $\Sigma$ is \emph{relatively principal}. 
\item If $\Sigma$ is not of the form $\Lambda _- ^{\scriptscriptstyle \leq \sigma + \Delta}$, we say that the initial segment $\Sigma$ is 
\begin{enumerate}
\item \emph{gapped} if the invariance group $\Delta$ is an immediate predecessor in ${\bf Cv}(\Gamma )$; 
\item \emph{tightened} if the invariance group $\Delta$ is not an immediate predecessor in ${\bf Cv}(\Gamma )$. 
\end{enumerate}
\end{enumerate} 

Let $\Omega$ be a final segment of the ordered group $\Gamma$, with invariance subgroup $\Delta$. 
\begin{enumerate} 
\item If there exists $\omega$ in $\Omega$ such that we have the equality $\Omega = \Lambda _+ ^{\scriptscriptstyle \geq \omega + \Delta}$, we say that the final segment $\Omega$ is \emph{relatively principal}. 
\item If $\Omega$ is not of the form $\Lambda _+ ^{\scriptscriptstyle \geq \omega + \Delta}$, we say that the final segment $\Omega$ is 
\begin{enumerate}
\item \emph{gapped} if the invariance group $\Delta$ is an immediate predecessor in ${\bf Cv}(\Gamma )$; 
\item \emph{tightened} if the invariance group $\Delta$ is not an immediate predecessor in ${\bf Cv}(\Gamma )$. 
\end{enumerate}
\end{enumerate} 
\end{definition} 

\vskip .2cm  

\begin{lemma}\label{le:jump-d-un-groupe}
Let $\Lambda$ be a cut of an ordered group $\Gamma$, if the cut is a jump the group $\Gamma$ is discrete. 
\end{lemma} 

\begin{preuve} 
If the cut $\Lambda$ is a jump, we deduce from remark \ref{rmq:jump} that there exist $\sigma$ and $\omega$ in $\Gamma$ such that $\sigma$ is an immediate predecessor of $\omega$, then the group $\Gamma$ is discrete from lemma \ref{le:groupe-discret}.  

\hfill$\Box$ 
\end{preuve} 

\vskip .2cm 

Let $\Lambda = (\Lambda _- , \Lambda _+)$ be a cut of the ordered group $\Gamma$. 
For any $\xi \in \Gamma$ we can define the symmetric interval $S_{(\Lambda , \xi )}$ and the convex subgroup $\Phi ^- {(\Lambda , \xi )}$ by $S_{(\Lambda , \xi )} = S_{(\Lambda _- , \xi )}$ and $\Phi ^- {(\Lambda , \xi )} = \Phi ^- {(\Lambda _-, \xi )}$ if $\xi$ belongs to the initial segment $\Lambda _-$ and by $S_{(\Lambda , \xi )} = S_{(\Lambda _+ , \xi )}$ and $\Phi ^- {(\Lambda , \xi )} = \Phi ^- {(\Lambda _+, \xi )}$ if $\xi$ belongs to the final segment $\Lambda _+$. In any case $\Phi ^- {(\Lambda , \xi )}$ is the greatest convex subgroup included in the symmetric interval $S_{(\Lambda , \xi )}$. 

\vskip .2cm 

\begin{proposition}\label{prop:segment-initial-final-principal}
Let $\Lambda$ be a cut of the ordered group $\Gamma$ with invariance subgroup $\Delta$, if both initial segment $\Lambda _-$ and final segment $\Lambda _+$ of the cut are relatively principal, the quotient group $\bar\Gamma = \Gamma / \Delta$ is discrete. 
In particular the group $\Delta$ is an immediate predecessor in ${\bf Cv}(\Gamma )$. 
\end{proposition} 

\begin{preuve} 
Let $\xymatrix{w_{\Delta} : \Gamma \ar[r] & \bar\Gamma}$ be the natural morphisme, then the cut $\Lambda$ induces a cut $\bar\Lambda = (\bar\Lambda _-,\bar\Lambda _+)$ of the group $\bar\Gamma$ defined by $\bar\Lambda _-= w_{\Delta}(\Lambda _-)$ and $\bar\Lambda _+= w_{\Delta}(\Lambda _+)$. 
The images by $w_{\Delta}$ of the initial segment $\Lambda _-^{\scriptscriptstyle \leq \sigma + \Delta}$ and of the final segment $\Lambda _+^{\scriptscriptstyle \geq \omega + \Delta}$ of $\Gamma$ are respectively the initial segment $\Lambda _-^{\leq \bar\sigma} = \bar\Gamma _{\leq \bar\sigma}$ and of the final segment $\Lambda _+^{\geq \bar\omega} =\bar\Gamma _{\geq \bar\omega}$ of $\bar\Gamma$, where $\bar\sigma = w_{\Delta} (\sigma )$ and $\bar\omega = w_{\Delta}(\omega )$. 

Then the cut $\bar\Lambda$ is a gap of the group $\bar\Gamma$, from lemma \ref{le:jump-d-un-groupe} we deduce that $\bar\Gamma$ is discrete, and the subgroup $(0)$ is an immediate predecessor in ${\bf Cv}(\Gamma /\Delta )$, hence $\Delta$ is an immediate predecessor in ${\bf Cv}(\Gamma )$. 

\hfill$\Box$ 
\end{preuve}

\vskip .2cm 

\begin{theorem}\label{th:structure-d-une-coupure} 
Let $\Lambda = (\Lambda _- , \Lambda _+)$ be a cut of the ordered group $\Gamma$, and let $\Delta = \Delta (\Lambda )$ be its invariance subgroup.  
\begin{enumerate} 
\item If there exists $\xi \in \Gamma$ which satisfies the equality $\Phi ^- (\Lambda ,\xi ) =  S_{(\Lambda ,\xi )}$, the invariance subgroup $\Delta$ is equal to $\Phi ^- (\Lambda ,\xi )$, and the cut $\Lambda$ is equal to $\Lambda ^{\scriptscriptstyle \leq \xi + \Delta}$ if $\xi$ belongs to the initial segment $\Lambda _-$ or is equal to $\Lambda ^{\scriptscriptstyle \geq \xi + \Delta}$ if $\xi$ belongs to the final segment $\Lambda _+$.   
\item If for any $\xi \in \Gamma$ we have the strict inequality $\Phi ^- (\Lambda ,\xi ) \subsetneq  S_{(\Lambda ,\xi )}$ and if the invariance subgroup is an immediate predecessor in $\widehat{\bf Cv}(\Gamma )$, there exist $\sigma \in \Lambda _-$ and $\omega \in \Lambda _+$ with $\Phi ^- (\Lambda ,\sigma ) =  \Phi ^- (\Lambda ,\omega ) =  \Delta$.  
\item If for any $\xi \in \Gamma$ we have the strict inequality $\Phi ^- (\Lambda ,\xi ) \subsetneq  S_{(\Lambda ,\xi )}$ and if the invariance subgroup is not an immediate predecessor in $\widehat{\bf Cv}(\Gamma )$, then for any $\xi \in \Gamma$ we have the strict inequality $\Delta \subsetneq \Phi ^- (\Lambda ,\xi )$. 
\end{enumerate} 
\end{theorem} 

\begin{preuve} 
This is a consequence of theorems \ref{th:segment-initial} and \ref{th:segment-final}. 

\hfill$\Box$ 
\end{preuve} 

\vskip .2cm  

\begin{corollary}\label{cor:groupe-bien-range} 
If the ordered group $\Gamma$ is well-ranked, for any cut $\Lambda$ of $\Gamma$ there exists $\xi \in \Gamma$ with the equality $\Delta (\Lambda ) = \Phi ^- (\Lambda ,\xi )$. 
\end{corollary} 

\begin{preuve} 
If the group $\Gamma$ is well-ranked any convex subgroup of $\Gamma$ is an immediate predecessor in $\widehat{\bf Cv}(\Gamma )$, then this is a consequence of theorem \ref{th:structure-d-une-coupure}. 

\hfill$\Box$ 
\end{preuve} 

\vskip .2cm  

\begin{definition}
Let $\Lambda$ be a cut of the totally ordered group $\Gamma$, with invariance subgroup $\Delta (\Lambda )$, then the cut is of one of the following \emph{types}. 
\begin{enumerate} 
\item If the cut $\Lambda$ is equal to $\Lambda ^{\scriptscriptstyle \leq \xi + \Delta}$, we say that the cut $\Lambda$ is \emph{relatively principal below}, and if the cut is equal to $\Lambda ^{\scriptscriptstyle \geq \xi + \Delta}$, we say that the cut $\Lambda$ is \emph{relatively principal above}. 
\item If the cut is relatively principal below or relatively principal above we say that the cut is \emph{relatively principal}.  
\item If the cut is both relatively principal below and relatively principal above we say that the cut is \emph{a relative jump}.
\item If the cut $\Lambda$ is not relatively principal and if the subgroup $\Delta (\Lambda )$ is an immediate predecessor in $\widehat{\bf Cv}(\Gamma )$, we say that the cut $\Lambda$ is \emph{gapped}. 
\item If the cut $\Lambda$ is not relatively principal and if the subgroup $\Delta (\Lambda )$ is not an immediate predecessor in $\widehat{\bf Cv}(\Gamma )$, we say that the cut $\Lambda$ is \emph{tightened}. 
\end{enumerate} 
\end{definition}

\vskip .2cm 

\begin{remark} 
\begin{enumerate}
\item A cut $\Lambda$ is relatively principal if and only the cut $\bar\Lambda$ induced in the quotient group $\Lambda / \Delta (\Lambda )$ is principal. 
\item A cut $\Lambda$ is a relative jump if and only if the cut $\bar\Lambda$ induced in the quotient group $\Lambda / \Delta (\Lambda )$ is a jump. 
\item If $\Lambda$ and $\Lambda '$ are two linearly equivalent cuts of $\Gamma$, they are of the same type. 
\end{enumerate}
\end{remark} 

\vskip .2cm 

We find in \cite{Ku-N} another description of the different types of cuts for an ordered group.
The authors organize the cuts into six different types, which can be compared with the classification we have just given.

The relatively principal cuts are called \emph{ball cuts} and they are noted $(b+G)^+$ and $(b+NG)^+$ if they are relatively principal below, and $(b+G)^-$ and $(b+NG)^-$ if they are relatively principal above. 
The gapped cuts are noted $nb+G$ and the tightened ones are noted $nb+NG$. 

The authors also introduce the notion of \emph{convexity gap}: a cut $\Lambda$ possesses a convexity gap if its invariance subgroup $\Delta (\Lambda )$ is an immediate predecessor in ${\bf Cv}(\Gamma )$. 

\vskip .2cm 

\begin{corollary} 
The cut $\Lambda$ has a convexity gap if and only if there exists $\sigma$ in the initial segment $\Lambda _-$ and $\omega$ in the final segment $\Lambda _+$ such that we have the equalities 
$$\Phi ^- (\Lambda _-,\sigma ) = \Phi ^- (\Lambda _+ ,\omega ) = \Delta (\Lambda ) \ .$$
\end{corollary} 

\begin{preuve} 
This is a consequence of theorems \ref{th:segment-initial} and \ref{th:segment-final}. 

\hfill$\Box$ 
\end{preuve}

\vskip .2cm 

Let $\Lambda$ be a cut of the ordered group $\Gamma$, let $\Delta = \Delta (\Lambda )$ be its invariance subgroup, and let let $\Theta$ be a convex subgroup of $\Gamma$ with $\Delta (\Lambda ) \subsetneq \Theta$. 
We recall that for any $\delta \in \Gamma$ which belongs to $(\Lambda _- + \Theta) \cap (\Lambda _+ + \Theta )$ we can define a non trivial cut $\Lambda ^{\delta}(\Theta ) = \bigl (\Lambda ^{\delta} _- (\Theta) , \Lambda ^{\delta} _+ (\Theta ) \bigr)$ of the group $\Theta$ by the following:  

$$\begin{array}{l}
\Lambda ^{\delta}_- (\Theta ) = \{ \zeta \in \Theta \ \hbox{such that } \ \delta + \zeta \in \Lambda _- \} \\
\Lambda ^{\delta}_+ (\Theta ) = \{ \zeta \in \Theta \ \hbox{such that } \ \delta + \zeta \in \Lambda _+ \} \ . 
\end{array}$$
We can replace the cut $\Lambda$ by the equivalent cut $\Lambda - \delta$, then we get a non trivial cut $\Lambda (\Theta)$ of the group $\Theta$ which is defined as the trace of the cut $\Lambda$ on the subgroup $\Theta$, i.e. we have 
$$\Lambda _- (\Theta ) = \Lambda _- \cap \Theta \quad\hbox{and}\quad \Lambda _+ (\Theta ) = \Lambda _+ \cap \Theta \ .$$

\vskip .2cm 

\begin{lemma}\label{le:invariant-par-trace}
Let $\Sigma$ be an initial segment of $\Gamma$ such that the trace $\Sigma (\Theta ) = \Sigma \cap \Theta$ is non trivial. 
Then for any $\xi$ in $\Sigma$ there exists $\sigma$ in $\Sigma (\Theta )$ with $\xi \leq \sigma$, and for any $\sigma$ in $\Sigma (\Theta )$ the subsets $S_{(\Sigma ,\sigma )}$ and $S_{(\Sigma (\Theta ),\sigma )}$ are equal. 
Moreover the convex subgroups $\Phi ^- (\Sigma ,\sigma )$ and $\Phi ^- (\Sigma (\Theta ),\sigma )$ are equal. 
\end{lemma}

\begin{preuve} 
As the initial segment $\Sigma (\Theta )$ is non trivial there exist $\zeta _1$ and $\zeta _2$ in $\Theta$ such that $\zeta _1$ belongs to $\Sigma$ and $\zeta _2$ doesn't belong to $\Sigma$. 

Let $\xi$ be an element of the initial segment $\Sigma$ with $\xi > \zeta _1$, then we have $\zeta _1 < \xi < \zeta _2$ and as $\Theta$ is a convex subgroup we deduce that $\xi$ belongs to $\Theta$. 

We recall that the subsets $T_{(\Sigma ,\sigma )}$ and $T_{(\Sigma (\Theta ),\sigma )}$ are defined by the following
$$T_{(\Sigma ,\sigma )} = \{ \xi \in \Gamma \ | \ \xi \geq 0 \ \hbox{and} \ \sigma + \xi \in \Sigma \} 
\quad\hbox{and}\quad 
T_{(\Sigma (\Theta ),\sigma )} = \{ \xi \in \Theta \ | \ \xi \geq 0 \ \hbox{and} \ \sigma + \xi \in \Sigma (\Theta )\}  \ .$$ 
Let $\xi \in T_{(\Sigma ,\sigma )}$, then the element $\zeta = \sigma + \xi$ satisfies $\zeta \in \Sigma$ and $\zeta \geq \sigma$ with $\sigma \in \Sigma (\Theta )$. We deduce from the above that $\zeta$ belongs to $\Sigma (\Theta )$ and that $\xi = \zeta - \sigma$ belongs to $\Theta$. 
Then we get the equality $T_{(\Sigma ,\sigma )} = T_{(\Sigma (\Theta ),\sigma )}$, hence the equality $S_{(\Sigma ,\sigma )} = S_{(\Sigma (\Theta ),\sigma )}$. 

As the subgroups $\Phi ^- (\Sigma ,\sigma )$ and $\Phi ^- (\Sigma (\Theta ),\sigma )$ are defined as the greatest convex subgroups contained respectively in $S_{(\Sigma ,\sigma )}$ and $S_{(\Sigma (\Theta ), \sigma )}$, we get the equality $\Phi ^- (\Sigma ,\sigma ) = \Phi ^- (\Sigma (\Theta ),\sigma )$. 

\hfill$\Box$ 
\end{preuve}

\vskip .2cm

Let $\Theta$ be a convex subgroup of $\Gamma$ and let $\xymatrix{w_{\Theta} : \Gamma \ar[r] & \bar\Gamma = \Gamma / \Theta}$ be the natural morphism. 
We consider an initial segment $\Sigma$ of $\Gamma$ and an element $\sigma$ in $\Sigma$, the image $\bar\Sigma = w_{\Theta}(\Sigma )$ of $\Sigma$ is an initial segment of $\bar\Gamma$ and the image $\bar\sigma = w_{\Theta}(\sigma )$ of $\sigma$ is an element of $\bar\Sigma$. 

\begin{lemma}\label{le:invariant-par-projection}
If the convex subgroup $\Theta$ is contained in the invariance subgroup $\Delta (\Sigma )$, for any $\sigma$ in $\Sigma$ the symmetric interval $S_{(\bar\Sigma ,\bar\sigma )}$ of $\bar\Gamma$ is equal to the image of the symmetric interval $S_{(\Sigma ,\sigma )}$ of $\Gamma$ by the morphism $w_{\Theta}$, and the convex subgroup $\Phi ^- (\bar\Sigma ,\bar\sigma )$ of $\bar\Gamma$ is the image of the convex subgroup $\Phi ^- (\Sigma ,\sigma )$ by the morphism $w_{\Theta}$. 
\end{lemma}

\begin{preuve} 
(i) To get the first part of the lemma it is enough to show that we have the equality $T_{(\bar\Sigma ,\bar\sigma )} = w_{\theta} \bigl ( T_{(\Sigma ,\sigma )} \bigr )$.  

 Let $\xi$ be an element of $T_{(\Sigma ,\sigma )}$, then we have $\xi \geq 0$ and $\xi + \sigma \in \Sigma$, we deduce that we have $w_{\Theta}(\xi ) \geq 0$ and $w_{\Theta}(\xi ) + \bar\sigma \in \bar\Sigma$, hence $w_{\Theta}(\xi )$ belongs to $T_{(\bar\Sigma ,\bar\sigma )}$. 

Conversely let $\bar\xi$ be an element of $T_{(\bar\Sigma ,\bar\sigma )}$, then $\bar\xi + \bar\sigma$ belongs to $\bar\Sigma$ and there exists $\zeta$ in $\Sigma$ with $w_{\Theta}(\zeta ) = \bar\xi + \bar\sigma$. 
Let $\xi$ be equal to $\zeta - \sigma$, then we get $w_{\Theta}(\xi ) = \bar\xi$, and by hypothesis we have $w_{\Theta}(\xi )\geq 0$, hence there exists $\theta \in \Theta$ such that $\xi \geq \theta$.  
Then, as we have the inclusion $\Theta \subset \Delta (\Sigma )$,  the element $\xi ' = \xi - \theta$ satisfies $\xi ' \geq 0$, $\xi ' + \sigma = \zeta - \theta \in \Sigma$, and $w_{\theta}(\xi ')=\bar\xi$. 
 
We also proved that we have the equality $S_{(\Sigma ,\sigma )} = w_{\theta} ^{-1} \bigl ( S_{(\bar\Sigma ,\bar\sigma )} \bigr )$.  Let $\xi$ be an element of $w_{\theta} ^{-1} \bigl ( T_{(\bar\Sigma ,\bar\sigma )} \bigr )$ with $\xi \geq 0$, then there exists $\xi '$ in $T_{(\Sigma ,\sigma )}$ such that we have $w_{\theta} (\xi ') = w_{\theta}(\xi )$, hence there exists $\theta$ in $\Theta$ such that $(\xi - \theta ) + \sigma$ belongs to $\Sigma$, and as above we deduce from the inclusion $\Theta \subset \Delta (\Sigma )$ that $\xi$ belongs to $T_{(\Sigma ,\sigma )}$.  
 
(ii) The image $w_{\theta}\bigl ( \Phi ^- (\Sigma ,\sigma ) \bigr )$ is a convex subgroup of $\bar\Gamma$ which is contained in the symmetric interval $S_{(\bar\Sigma ,\bar\sigma )} = w_{\theta} \bigl ( S_{(\Sigma ,\sigma )} \bigr )$, hence we have the inclusion $w_{\theta}\bigl ( \Phi ^- (\Sigma ,\sigma ) \bigr ) \subset \Phi ^- (\bar\Sigma ,\bar\sigma )$. 

The inverse image $w_{\Theta} ^{-1} \bigl ( \Phi ^- (\bar\Sigma ,\bar\sigma ) \bigr )$ of $\Phi ^- (\bar\Sigma ,\bar\sigma )$ is a convex subgroup $\Psi$ of $\Gamma$ which satisfies $\Psi \subset S_{(\Sigma ,\sigma )} = w_{\theta} ^{-1} \bigl ( S_{(\bar\Sigma ,\bar\sigma )} \bigr )$, then we get the inclusion $\Psi \subset \Phi ^- (\Sigma ,\sigma )$, hence the result. 

\hfill$\Box$ 
\end{preuve} 

\vskip .2cm

\begin{proposition} \label{prop:image-de-principal}
(i) Let $\Delta$ and $\Theta$ be two convex subgroups of $\Gamma$ with the inclusion $\Theta \subset \Delta$, and let $\bar\Delta$ be the convex subgroup of $\bar\Gamma$ image of $\Delta$ by the natural morphism $\xymatrix{\Gamma \ar[r] & \bar\Gamma = \Gamma / \Theta}$.  
\begin{enumerate}
\item Let $\Lambda = \Lambda ^{\scriptscriptstyle \leq \zeta + \Delta}$ (resp. $\Lambda = \Lambda ^{\scriptscriptstyle \geq \zeta + \Delta}$) be the relatively principal cut of $\Gamma$ associated with $\Delta$ and with an element $\zeta$ of $\Gamma$, then its image $\bar\Lambda
= w_{\Theta}(\Lambda )$ of the cut $\Lambda$ is the relatively principal cut $\Lambda ^{\scriptscriptstyle \leq \bar\zeta + \bar\Delta}$ (resp. $\Lambda ^{\scriptscriptstyle \geq \bar\zeta + \bar\Delta}$) of $\bar\Gamma$ associated with the convex subgroup $\bar\Delta$ and with the element $\bar\zeta = w_{\Theta}(\zeta )$. 

\item Conversely a cut $\Lambda$ of $\Gamma$ whose image $\bar\Lambda = w_{\Theta}(\Lambda )$ is the relatively principal cut $\Lambda ^{\scriptscriptstyle \leq \bar\zeta + \bar\Delta}$ (resp. $\Lambda ^{\scriptscriptstyle \geq \bar\zeta + \bar\Delta}$) of $\bar\Gamma$ associated with the convex subgroup $\bar\Delta$ and with an element $\bar\zeta$ of $\bar\Gamma$, is equal to the relatively principal cut $\Lambda ^{\scriptscriptstyle \leq \zeta + \Delta}$ (resp. $\Lambda ^{\scriptscriptstyle \geq \zeta + \Delta}$) for any element $ \zeta$ of $\Gamma$ such that we have $w_{\Theta}(\zeta )=\bar\zeta$.   
\end{enumerate} 

\noindent (ii) Let $\Delta$ and $\Theta$ be two convex subgroups of $\Gamma$ with the strict inclusion $\Delta \subsetneq \Theta$.
\begin{enumerate} 
\item Let $\Lambda = \Lambda ^{\scriptscriptstyle \leq \zeta + \Delta}$  (resp. $\Lambda = \Lambda ^{\scriptscriptstyle \geq \zeta + \Delta}$) be the relatively principal cut of $\Gamma$ associated with $\Delta$ and with an element $\zeta$ of $\Gamma$, then the element $\zeta$ belongs to $(\Lambda _- + \Theta ) \cap (\Lambda _+ + \Theta)$, and the non trivial cut $\Lambda ^{\zeta}(\Theta )$ of the group $\Theta$ defined above is the relatively principal cut $\Lambda ^{\scriptscriptstyle \leq \Delta}$ (resp. $\Lambda ^{\scriptscriptstyle \geq \Delta}$) of $\Theta$, where $\Delta$ is view as a convex subgroup of $\Theta$. 

\item  Conversely any relatively principal cut $\Lambda$ of the group $\Theta$ is the trace of a relatively principal cut of $\Gamma$. 
\end{enumerate}
\end{proposition}

\begin{preuve} 
We will consider the relatively principal cuts of the form $\Lambda ^{\scriptscriptstyle \leq \zeta + \Delta}$, the demonstration is identical for those of the form $\Lambda ^{\scriptscriptstyle \geq \zeta + \Delta}$. 

(i)(1) Let $\bar\xi$ be an element of $\bar\Gamma$ and let $\xi$ be an element of $\Gamma$ with $w_{\Theta}(\xi ) = \bar\xi$,  
then we have the equivalences 
$$\begin{array}{rcl}
\bar\xi \in  \Lambda ^{\scriptscriptstyle \leq \bar\zeta + \bar\Delta}_+ &\Longleftrightarrow & \bar\xi - \bar\zeta  > \bar\delta \quad \forall \, \bar\delta \in \bar\Delta \\ 
&\Longleftrightarrow & w_{\Theta}(\xi - \zeta ) > w_{\Theta}(\delta ) \quad \forall \, \delta \in \Delta \\
&\Longleftrightarrow & \xi - \zeta > \delta  \quad \forall \, \delta \in \Delta \quad\hbox{because}\quad \Theta \subset \Delta \\ 
& \Longleftrightarrow & \xi \in  \Lambda ^{\scriptscriptstyle \leq \zeta + \Delta}_+ \ . 
\end{array}$$ 

(i)(2) This is a consequence of (i)(1) and of proposition \ref{prop:passage-au-quotient}.

\vskip .2cm 

(ii)(1) The non trivial cut $\Lambda ^{\zeta}(\Theta ) = \bigl (\Lambda ^{\zeta} _- (\Theta) , \Lambda ^{\zeta} _+ (\Theta ) \bigr)$ of the group $\Theta$ is given by the following:  

$$\begin{array}{l}
\Lambda ^{\zeta}_- (\Theta ) = \{ \xi \in \Theta \ \hbox{\rm such that } \ \zeta + \xi \in \Lambda _- \} = \{ \xi \in \Theta \ | \ \exists \, \delta \in \Delta \ \hbox{\rm such that} \ \xi \leq \delta \} \\
\Lambda ^{\zeta}_+ (\Theta ) = \{ \xi \in \Theta \ \hbox{\rm such that } \ \zeta + \xi \in \Lambda _+ \} = \{ \xi \in \Gamma \ | \ \xi >\Delta \}  \ , 
\end{array}$$
then we get the equality
$\Lambda ^{\zeta}(\Theta ) = \Lambda ^{\scriptscriptstyle \leq \Delta}$, where $\Delta$ is view as a convex subgroup of $\Theta$.  

(ii)(2) The relatively principal cut $\Lambda ^{\scriptscriptstyle \leq \zeta + \Delta}$ of $\Theta$, where $\zeta$ is an element of $\Theta$ and $\Delta$ is a convex subgroup of $\Theta$, is the trace on $\Theta$ of the relatively principal cut $\Lambda ^{\scriptscriptstyle \leq \zeta + \Delta}$ of $\Gamma$, where we consider now $\zeta$ as an element of $\Gamma$ and $\Delta$ as a convex subgroup of $\Gamma$.

\hfill$\Box$ 
\end{preuve} 

\vskip .2cm

Let $\Lambda$ be a non trivial cut of the ordered group $\Gamma$, with invariance subgroup $\Delta = \Delta (\Lambda )$, and let $\Theta _1$ and $\Theta _2$ be two convex subgroups of $\Gamma$ with $\Theta _1 \subset \Delta \subsetneq \Theta _2$. 
We recall that any cut $\Lambda '$ of the group $\Gamma ' =\Theta _2 / \Theta _1$ whose linear equivalence class is the image of the class of $\Lambda$ by the map $T_{\Theta _1}^{\Theta _2}$, is defined as follows:
let $\delta$ be an element of $\Gamma$ that belongs to $(\Lambda _- + \Theta _2) \cap (\Lambda _+ + \Theta _2)$, and let $\Lambda ^{\delta}(\Theta _2) = \bigl (\Lambda ^{\delta} _- (\Theta _2) , \Lambda ^{\delta} _+ (\Theta _2) \bigr)$ be the non trivial cut of the group $\Theta _2$ defined above, then the cut $\Lambda '$ is the cut induced on $\Gamma ' = \Theta _2 / \Theta _1$.

\vskip .2cm  

\begin{proposition}\label{prop:sous-quotient} 
The cut $\Lambda '$ of the group $\Gamma '$ is of the same type as the cut $\Lambda$ of $\Gamma$. 
\end{proposition} 

\begin{preuve} 
The invariance subgroup $\Delta ' = \Delta (\Lambda ')$ of the cut $\Lambda '$ is equal to the subgroup $\Delta / \Theta _1$ of $\Gamma '$, then $\Delta '$ is an immediate predecessor in $\widehat{\bf Cv}(\Gamma ')$ if and only if $\Delta$ is an immediate predecessor in $\widehat{\bf Cv}(\Gamma )$. 
Then the proposition is a consequence of lemmas \ref{le:invariant-par-trace}, \ref{le:invariant-par-projection}, and of proposition \ref{prop:image-de-principal}.  

\hfill$\Box$ 
\end{preuve}

\vskip .2cm

Let $\Lambda$ be a cut of the totally ordered group $\Gamma$ such that its invariance subgroup $\Delta$ is an immediate predecessor, and we denote by $\bar\Delta$ its immediate successor in $\widehat{\bf Cv}(\Gamma )$. 
We consider the previous situation with $\Theta _1 = \Delta \subsetneq \bar\Delta = \Theta _2$, and the group $\Gamma _{\Lambda}= \bar\Delta / \Delta$ is an ordered group of rank one. 

Then we want to study the cuts of an ordered group $\Theta$ of rank one. 
Let $\Theta$ be such an ordered group, we consider one completion $\xymatrix{c: \Theta \ar[r] & \R}$ of $\Theta$, and we identify $\Theta$ with its image in $\R$.  

\begin{proposition}\label{prop:groupe-de-rang-un} 
Let $\Theta$ be a discrete ordered group of rank one, with $c(\Theta )= \Z$, then any non trivial cut $\Lambda = (\Lambda _- , \Lambda _+)$ of $\Theta$ is a jump, i.e. there exists an integer $n \in \Z$ with $\Lambda _- = \Theta _{\leq n}$ and $\Lambda _+ = \Theta _{\geq (n +1)}$. 

\vskip .2cm 

Let $\Theta$ be a continuous ordered group of rank one, then any cut $\Lambda = (\Lambda _- , \Lambda _+)$ of $\Theta$ is of one of the following form: 
\begin{enumerate} 
\item the initial segment $\Lambda _-$ is equal to $\Theta _{\leq \theta}$ for some $\theta \in \Theta$ and the final segment $\Lambda _+$ has no smallest element; 
\item the final segment $\Lambda _+$ is equal to $\Theta _{\geq \theta}$ for some $\theta \in \Theta$ and the initial segment $\Lambda _-$ has no greatest element; 
\item the initial segment $\Lambda _-$ has no greatest element, the final segment $\Lambda _+$ has no smallest element, and there exists a unique $\delta \in \R \setminus \Theta$ such that $\Lambda _- = \R _{\scriptscriptstyle \geq \delta} \cap \Theta$ and $\Lambda _+ = \R _{\scriptscriptstyle \leq \delta} \cap \Theta$. 
\end{enumerate}
\end{proposition} 

\begin{preuve}   
Obvious. 

\hfill$\Box$ 
\end{preuve} 

\vskip .2cm  

If the group $\Gamma$ is discrete, or if the group $\Gamma$ is continuous and we are in the first two cases of the proposition, the cut is principal. 
Otherwise the cut is gapped.

\vskip .2cm 

\subsection{Image of a cut by a morphism} 

Let $\xymatrix{\varphi: \Gamma \ar[r] & \Gamma '}$ be an injective morphism of totally ordered groups, we have defined a morphism $\xymatrix{\varepsilon_{\varphi}:{\widehat{\bf Cv}}(\Gamma ) \ar[r] & {\widehat{\bf Cv}}(\Gamma ')}$ of ordered sets such that for any convex subgroup $\Delta$ of $\Gamma$, its image $\Delta '= \varepsilon _{\varphi}(\Delta )$ is the smallest convex subgroup of $\Gamma '$ that contains $\varphi(\Delta )$, we recall it is defined by 
$$\varepsilon _{\varphi}(\Delta )= \{ \zeta \in \Gamma ' \ | \ \exists \ \xi \in \Delta  _{\geq 0} \ \hbox{such that} \ -\varphi (\xi ) \leq \gamma \leq \varphi (\xi ) \} \ .$$ 

We can also define a morphism $\xymatrix{\varepsilon ^{\varphi}:{\widehat{\bf Cv}}(\Gamma ') \ar[r] & {\widehat{\bf Cv}}(\Gamma )}$ of ordered sets such that for any convex subgroup $\Delta '$ of $\Gamma '$, its image $\varepsilon^{\varphi}(\Delta ')$ is the inverse image of $\Delta$ by $\varphi$. If we consider the injective morphism $\varphi$ as the inclusion $\Gamma \subset \Gamma '$, we have $\varepsilon^{\varphi}(\Delta ') = \Delta ' \cap \Gamma$.

\begin{remark}\label{rmq:cofinal}
If the convex subgroup $\Delta$ is different from $\Gamma$, its image $\varepsilon_{\varphi}(\Delta )$ is different from $\Gamma '$, then the morphism $\varepsilon_{\varphi}$ induces a morphism $\xymatrix{\varepsilon_{\varphi}:{\bf Cv}(\Gamma ') \ar[r] & {\bf Cv}(\Gamma )}$. 

On the other hand, it may happen that the image of a proper convex subgroup $\Delta '$ of $\Gamma '$ is the entire group $\Gamma$. 
If this is not the case, that is to say if the morphism $\varepsilon^{\varphi}$ induces a morphism $\xymatrix{\varepsilon^{\varphi}:{\bf Cv}(\Gamma ') \ar[r] & {\bf Cv}(\Gamma )}$, we say that the subgroup $\Gamma$ is \emph{cofinal} in $\Gamma '$. 

The subgroup $\Gamma$ is cofinal in $\Gamma '$ if it is not contained in any proper subgroup of $\Gamma '$, which is equivalent to the equality $\varepsilon_{\varphi}(\Gamma )=\Gamma '$. 
\end{remark} 

\vskip .2cm 

Let $\xymatrix{\varphi :\Gamma \ar[r] & \Gamma '}$ be an injective morphism of ordered groups, and as before we assume that $\Gamma$ is a subgroup of $\Gamma '$, and let $\xymatrix{\varepsilon^{\varphi}:\widehat{\bf Cv}(\Gamma ') \ar[r] & \widehat{\bf Cv}(\Gamma )}$ and $\xymatrix{\varepsilon_{\varphi}:\widehat{\bf Cv}(\Gamma ) \ar[r] & \widehat{\bf Cv}(\Gamma ')}$ be the morphisms of ordered sets defined above. 
Then for any convex subgroup $\Delta$ of $\Gamma$ we have the equality $\varepsilon _{\varphi} (\Delta ) \cap \Gamma = \Delta$, hence $\varepsilon ^{\varphi} \circ \varepsilon _{\varphi} = id_{{\bf Cv}(\Gamma )}$, and for any convex subgroup $\Delta '$ of $\Gamma '$ we have only the inclusion $\varepsilon _{\varphi}(\Delta ' \cap \Gamma ) \subset \Delta '$. 

\begin{definition} 
We say that a subgroup $\Gamma$ of a totally ordered set $\Gamma '$ is \emph{convex dense}, or that an injective morphism $\xymatrix{\varphi :\Gamma \ar[r] & \Gamma '}$ of ordered groups is \emph{convex dense}, if the morphism $\varepsilon^{\varphi}$ is an isomorphism of ordered sets, or equivalently if the morphism $\varepsilon_{\varphi}$ is an isomorphism. 

This is also equivalent to say that for any convex subgroup $\Delta '$ of $\Gamma '$ we have the equality $\Delta ' = \varepsilon _{\varphi}(\Delta ' \cap \Gamma )$.
\end{definition} 

We have seen that for any totally ordered set $\Gamma$ the natural morphism $\xymatrix{\varphi :\Gamma \ar[r] & \Gamma ^{\rm div}}$, where $\Gamma ^{\rm div}$ is the divisible hull of $\Gamma$, is convex dense. 
In the same way, if we denote $\bigl ( {\bf I}_{\Gamma} , (\Theta _i) _{i \in {\bf I}_{\Gamma}} \bigr )$ the skeleton of $\Gamma$, any injective morphism of totally ordered ${\bf I}_{\Gamma}$-groups    
$\xymatrix{\psi : \Gamma \ar[r] & \Gamma ^{\rm cp} =  \R ^{[{\bf I}_{\Gamma}]} _{\rm lex}}$ 
as defined in corollary \ref{cor:complete}, is convex dense.

\vskip .2cm 

Let $\Gamma$ be a totally ordered group, then for any initial or final segment $\Sigma$ of $\Gamma$ we have introduced its invariance subgroup $\Delta (\Sigma )$, it is the convex subgroup of $\Gamma$ defined by 
$$\Delta (\Sigma) = \{ \xi \in \Gamma | \Sigma + \xi = \Sigma \} \ .$$ 
If $\Sigma$ is an initial segment (resp. a final segment) we can also defined the subgroup  $\Delta (\Sigma )$ by the following 
$$ 
\Delta (\Sigma) _{\geq 0} = \{ \xi \in \Gamma _{\geq 0} | \Sigma + \xi \subset \Sigma \} \quad  
\Bigl ( \hbox{resp. } \Delta (\Sigma) _{\leq 0} = \{ \xi \in \Gamma _{\leq 0} | \Sigma + \xi \subset \Sigma \} \Bigr ) \ . $$

\begin{proposition}\label{prop:image-directe-image-inverse}
Let $\xymatrix{\varphi : \Gamma \ar[r] & \Gamma '}$ be an injective morphism of totally ordered groups, then for any initial segment $\Sigma '$ of $\Gamma '$ we have the inclusion:
$$\varepsilon ^{\varphi}(\Delta (\Sigma ')) \ \subset \ \Delta (\varphi ^*(\Sigma '))  \ ,$$ 
and for any initial segment $\Sigma$ of $\Gamma$ we have the inclusion:
$$\varepsilon _{\varphi}(\Delta (\Sigma )) \ \subset \ \Delta (\varphi _!(\Sigma ))  \ .$$ 

Moreover if the morphism $\varphi$ is convex dense we have the equality: 
$$\varepsilon _{\varphi}(\Delta (\Sigma )) \ = \ \Delta (\varphi _!(\Sigma )) \ .$$
\end{proposition}

\begin{preuve}   
(i) By definition we have the equalities 
$\Delta (\varphi ^*(\Sigma ')) _{\geq 0} = \{ \zeta \in \Gamma _{\geq 0} \ | \ \zeta + \varphi ^{-1}(\Sigma ') \subset \varphi ^{-1}(\Sigma ') \}$ and 
$\varepsilon ^{\varphi}(\Delta (\Sigma ')) _{\geq 0} = \{ \zeta \in \Gamma _{\geq 0} \ | \ \varphi (\zeta ) + \Sigma ' \subset \Sigma '\}$. 
Then if $\zeta$ is an element of $\varepsilon ^{\varphi}(\Delta (\Sigma ')) _{\geq 0}$, for any $\xi$ in $\varphi ^{-1}(\Sigma ')$ we have 
$$\varphi (\xi ) \in \Sigma ' \ \Longrightarrow \ \varphi (\zeta ) + \varphi (\xi ) \in \Sigma ' \ \Longrightarrow \ \varphi (\zeta + \xi) \in \Sigma ' \ \Longrightarrow \ \zeta + \xi \in \varphi ^{-1}(\Sigma ') \ ,$$   
hence $\zeta $ belongs to $\Delta (\varphi ^*(\Sigma ')) _{\geq 0}$. 

\vskip .2cm

(ii) Let $\Sigma$ be an initial segment of $\Gamma$, we recall the equality   
$\varphi _!(\Sigma ) =\bigcup _{\xi \in \Sigma} \Gamma '_{\leq \varphi (\xi )}$,  
then an element $\zeta '$ of $\Gamma ' _{\geq 0}$ belongs to $\Delta (\varphi _!(\Sigma ) _{\geq 0}$ if and only if for any $\xi _1 \in \Sigma$ there exists $\xi _2 \in \Sigma$ such that $\zeta '+ \varphi (\xi _1) \leq \varphi (\xi _2)$. 

If $\zeta '$ be an element of $\Gamma ' _{\geq 0}$ that belongs to $\varepsilon _{\varphi}(\Delta (\Sigma )) _{\geq 0}$, there exists $\zeta \in \Delta (\Sigma ) _{\geq 0}$ such that $\zeta '\leq \varphi (\zeta )$. 
Then for any $\xi _1 \in \Sigma$ we have $\xi _2 =\zeta + \xi _1$ that belongs to $\Sigma$ and $\zeta '+ \varphi (\xi _1) \leq \varphi (\xi _2)$, and $\zeta '$ belongs to $\Delta (\varphi _!(\Sigma ) _{\geq 0}$.

\vskip .2cm

(iii) We assume now that the morphism $\varphi$ is convex dense, then for any initial segment $\Sigma$ of $\Gamma$ we deduce from the equalities $\varphi ^* \varphi _! (\Sigma ) = \Sigma$ and $\varepsilon _{\varphi} \circ \varepsilon ^{\varphi} = id_{{\bf Cv}(\Gamma ')}$ the inclusion
$$\Delta (\varphi _! (\Sigma )) = \varepsilon _{\varphi} \varepsilon ^{\varphi} (\Delta (\varphi _! (\Sigma ))) \subset \varepsilon _{\varphi} (\Delta (\varphi ^* \varphi _! (\Sigma ))) =  \varepsilon _{\varphi} (\Delta (\Sigma )) \ .$$

\hfill$\Box$ 
\end{preuve} 

\vskip .2cm 

\begin{remark}\label{rmq:image-directe-exceptionnelle}
For any totally ordered group $\Gamma$ we can consider the group $\Gamma ^{op}$, which the same group with the opposite order. 
Then $\Gamma ^{op}$ is also a totally ordered group and the identity induces a bijection of ordered sets between ${\bf Cv}(\Gamma )$ and ${\bf Cv}(\Gamma ^{op})$. 
The set ${\bf FS}(\Gamma )$ of final segments of $\Gamma$ is equal to the set ${\bf IS}(\Gamma ^{op})$ of initial segments of $\Gamma ^{op}$, then the application $\xymatrix{(.)^C_{\Gamma} : \Sigma \ar@{|->} [r] & \Sigma ^C}$ is a bijection between the sets ${\bf IS}(\Gamma )$ and ${\bf IS}(\Gamma ^{op})$, with the equality $\bigl ((.)^C_{\Gamma}\bigr ) ^{-1} = (.)^C_{\Gamma ^{op}}$. 

Any morphism $\xymatrix{\varphi : \Gamma \ar[r] & \Gamma '}$ of totally ordered groups induces a morphism $\xymatrix{\varphi ^{op}: \Gamma ^{op} \ar[r] & {\Gamma '}^{op}}$, and we have the equality 
$\varphi _* = (.)^C_{{\Gamma '}^{op}} \circ  \bigl ( \varphi ^{op}\bigr ) _! \circ (.)^C_{\Gamma}$. 

\vskip .2cm 

Then we have the following result, which is the analogous of proposition \ref{prop:image-directe-image-inverse} for the exceptionnal image. 
Let $\xymatrix{\varphi : \Gamma \ar[r] & \Gamma '}$ be an injective morphism of totally ordered groups, then for any initial segment $\Sigma$ of $\Gamma$ we have the inclusion:
$$\varepsilon _{\varphi}(\Delta (\Sigma )) \ \subset \ \Delta (\varphi _*(\Sigma ))  \ ,$$ 
and if the morphism $\varphi$ is convex dense we have the equality: 
$$\varepsilon _{\varphi}(\Delta (\Sigma )) \ = \ \Delta (\varphi _*(\Sigma )) \ .$$
\end{remark}

\vskip .2cm  

\begin{remark}\label{rmq:inegalite} 
In general, even for a convex dense subgroup $\Gamma$ of $\Gamma '$, we have the strict inclusion 
$$\varepsilon ^{\varphi}(\Delta (\Sigma ')) \ \subsetneq \ \Delta (\varphi ^*(\Sigma '))  \ .$$

\vskip .2cm 

Consider the ordered group $\Gamma '$ defined as the group $\R ^2$ with the lexicographic order, and the injective morphism $\xymatrix{\varphi :\Gamma \ar[r] & \Gamma '}$, where $\Gamma$ is the convex dense subgroup $(\Z \oplus \R) _{\rm lex}$ of $\Gamma '$. 

Let $\Sigma '$ be the initial segment of $\Gamma '$ defined by $\Sigma ' = \Gamma ' _{\leq \zeta}$ where $\zeta = (z_1,z_2)$ is an element of $\Gamma '$, then its invariance subgroup $\Delta (\Sigma ')$ is the convex subgroup $(0)$ of $\Gamma '$, and let $\Sigma$ be the initial segment of $\Gamma$ defined by $\Sigma = \varepsilon ^{\varphi} (\Sigma ') = \Sigma ' \cap \Gamma$. 

\begin{enumerate} 
\item If $\zeta$ belongs to $\Gamma$, i.e. if $z_1$ belongs to $\Z$, the initial segment $\Sigma$ is equal to $\Gamma _{\leq \zeta}$, and its invariance subgroup $\Delta (\Sigma )$ is still equal to the convex subgroup $(0)$ of $\Gamma$. Then we have the equality $\varepsilon ^{\varphi}(\Delta (\Sigma ')) = \Delta (\varphi ^*(\Sigma '))$.   

\item If $\zeta$ doesn't belong to $\Gamma$, let $y_1$ the integral part of $z_1$, i.e. $y_1$ is the element of $\Z$ with $y_1 < z_1 < y_1 +1$, an element $\xi = (x_1,x_2)$ in $\Z \times \R$ satisfies $\xi \leq \zeta$ if and only if we have $x_1 \leq y_1$, then for any element $\zeta '$ with $\zeta ' = (y_1, z'_2)$, the initial segment $\Sigma$ is equal to 
$ \Lambda ^{\scriptscriptstyle \leq \zeta '+ \Delta _1}_- = \{ \xi \in \Gamma \ | \ \exists \, \delta \in \Delta _1 \ \hbox{such that} \ \xi - \zeta ' \leq \delta \} $ 
where $\Delta _1$ is the convex subgroup $\bigl ((0) \oplus \R \bigr )$ of $\Gamma$, and the invariance subgroup $\Delta (\Sigma )$ is equal to $\Delta _1$. 
Then we have the strict inclusion $\varepsilon ^{\varphi}(\Delta (\Sigma ')) = (0) \subsetneq \Delta _1 = \Delta (\varphi ^*(\Sigma '))$.
\end{enumerate}
\end{remark}

\vskip .2cm

We want now to study the images of the initial segments $\Lambda _-^{>\zeta + \Delta}$ and $\Lambda _-^{<\zeta + \Delta}$ associated to relatively principal cuts of $\Gamma$, respectively by the morphisms $\varphi _!$ and $\varphi _*$. 

\begin{proposition}\label{prop:image-de-principal-2} 
Let $\xymatrix{\varphi : \Gamma \ar[r] & \Gamma '}$ be an injective morphism of ordered groups, and let $\xymatrix{\varepsilon_{\varphi}:{\widehat{\bf Cv}}(\Gamma ) \ar[r] & {\widehat{\bf Cv}}(\Gamma ')}$ be the morphism of ordered sets defined above. 
Then for any element $\xi$ of $\Gamma$ and any convex subgroup $\Delta$ of $\Gamma$ we have the equalities 
$$\begin{array}{l}
\varphi _!(\Lambda _-^{<\zeta + \Delta}) = \Lambda _-^{<\varphi(\zeta )+ \varepsilon _{\varphi}(\Delta )} \\ 
\varphi _*(\Lambda _-^{>\zeta + \Delta}) = \Lambda _-^{>\varphi(\zeta )+ \varepsilon _{\varphi}(\Delta )} \ .
\end{array}$$ 

\end{proposition} 

\begin{preuve} 

Let $\xi '$ be an element of $\Gamma '$ then we have the following equivalences: 
$$\begin{array}{rcl} 
\xi ' \in \Lambda _-^{<\varphi(\zeta )+ \varepsilon _{\varphi}(\Delta )} & \Longleftrightarrow & \exists \, \delta ' \in \varepsilon _{\varphi}(\Delta ) \quad \hbox{such that} \quad \xi ' \leq \varphi (\zeta ) + \delta ' \\ 
& \Longleftrightarrow & \exists \, \delta \in \Delta \quad \hbox{such that} \quad \xi ' \leq \varphi (\zeta + \delta ) \\ 
& \Longleftrightarrow & \exists \, \xi \in \Lambda _-^{<\zeta + \Delta} \quad \hbox{such that} \quad \xi ' \leq \varphi (\xi ) \ ,
\end{array}$$ 
hence the equality $\varphi _!(\Lambda _-^{<\zeta + \Delta}) = \Lambda _-^{<\varphi(\zeta )+ \varepsilon _{\varphi}(\Delta )}$ by definition of the functor $\varphi _!$. 

We get the equality $\varphi _*(\Lambda _-^{>\zeta + \Delta}) = \Lambda _-^{>\varphi(\zeta )+ \varepsilon _{\varphi}(\Delta )}$ by considering as in remark \ref{rmq:image-directe-exceptionnelle} the morphism $\xymatrix{\varphi ^{op}: \Gamma ^{op} \ar[r] & {\Gamma '}^{op}}$. 

\hfill$\Box$ 
\end{preuve} 

For any morphism $\xymatrix{\varphi : \Gamma \ar[r] & \Gamma '}$ of totally ordered groups we have defined the two functors $\xymatrix{\varphi _!: {\bf Cp}(\Gamma ) \ar[r] &{\bf Cp}(\Gamma ')}$ and $\xymatrix{\varphi _*: {\bf Cp}(\Gamma ) \ar[r] &{\bf Cp}(\Gamma ')}$, and we deduce from proposition \ref{prop:image-de-principal-2} that in the case of an injective morphisms $\varphi$, the images of a relatively principal cut $\Lambda$ of $\Gamma$ by the functors $\varphi _!$ and $\varphi _*$ are relatively principal. 

\begin{remark}\label{rmq:image-de-non-principal}
It may happen that a cut $\Lambda$ of the group $\Gamma$ has its image by $\varphi _!$ or by $\varphi _*$ that is relatively principal, but that the cut $\Lambda$ is not, even if the injective morphism $\xymatrix{\varphi : \Gamma \ar[r] & \Gamma '}$ is convex dense. 

For example, let $\xymatrix{\varphi : \Q \ar[r] & \R}$ the natural embedding, let $\alpha$ in $\R \setminus \Q$, and let $\Sigma$ be the non principal initial segment of $\Q$ defined by $\Sigma = \Q _{\leq \alpha} = \Q _{< \alpha}$. 
Then the initial segments $\varphi _!(\Sigma ) = \R _{\leq \alpha}$ and $\varphi _*(\Sigma ) = \R _{< \alpha}$ are principal initial segments of $\R$. 

But if we know that the image by $\varphi _!$ (resp. by $\varphi _*$)of a non principal cut $\Lambda$ of $\Gamma$ is still non relatively principal, the cuts $\Lambda$ and $\varphi _!(\Lambda )$ (resp. the cuts $\Lambda$ and $\varphi _*(\Lambda )$) are of the same type. 

Indeed, in this case the type is determined only by the fact that the invariance subgroup is or is not an immediate predecessor. And we know that we have the equality $\varepsilon _{\varphi}(\Delta (\Sigma )) \ = \ \Delta (\varphi _!(\Sigma ))$, and that the morphism $\varepsilon _{\varphi}$ is a bijection of ordered sets. 
\end{remark}

\vskip .2cm 

Let $\Theta = \prod _{i \in I} ^{(H)} \Theta _i$ be the Hahn product of a family of ordered groups, 
and for any $i$ in $I$ let $\xymatrix{\vartheta _i : \Theta _i \ar@<-1.5pt>@{^{(}->}[r] & \Theta}$ be the injective morphism of ordered groups defined by $\vartheta _i (x_i) = \underline y$ where the element $\underline y = ( y_j) _{j \in I}$ satisfies $y_j = 0$ for $j \not= i$ and $y_i = x_i$. 

Let $\Sigma$ be an initial segment of $\Theta$, we define the subset ${\bf L}_+^{(\Sigma )}$ of $I$ by the following 
$${\bf L}_+^{(\Sigma )} \ = \ \{ i \in I \ | \  \vartheta _i(\Theta _i ) + \Sigma = \Sigma \} \ .$$

\begin{lemma}\label{le:segment-final-de-I}
The subset ${\bf L}_+^{(\Sigma )}$ is a final segment of the ordered set $I$. 
\end{lemma} 

\begin{preuve}   
Let $i$ and $j$ in $I$ with $j > i$, we have $\vartheta _j(\Theta _j) < \bigl (\vartheta _i(\Theta _i) \bigr )_{>0}$, i.e. for any element $\underline y$ in $\vartheta _i(\Theta _i)$ with $y_i >0$ and any element $\underline z$ in $\vartheta _i(\Theta _j)$ we have the inequality $\underline z < \underline y$. 

Hence the result because $\Sigma$ is an initial segment of $\Theta$. 

\hfill$\Box$ 
\end{preuve} 

\vskip .2cm  

We define a cut ${\bf L}^{(\Sigma )} = \bigl ({\bf L}_-^{(\Sigma )},{\bf L}_+^{(\Sigma )}\bigr )$ of $I$, where ${\bf L}_-^{(\Sigma )} = I \setminus {\bf L}_+^{(\Sigma )}$ may also be defined by 

$${\bf L}_-^{(\Sigma )} \ = \ \{ i \in I \ | \ \exists \, x_i \in \Theta _i \ \hbox{and} \ \exists \, \underline y  \in \Sigma \  \hbox{such that} \  \vartheta _i(x_i) + \underline y \notin \Sigma \} \ .$$

\vskip .2cm 

We recall that with any convex subgroup $\Delta$ of the group $\Theta = \prod _{i \in I} ^{(H)} \Theta _i$ we have associated a cut ${\bf L} ^{\scriptscriptstyle \geq (\Delta)}$ of $I$ defined by 
$${\bf L} ^{\scriptscriptstyle \geq (\Delta)} _- = \{ i \in I \ | \ \Delta \subsetneq \Theta _{{\bf L} ^{\geq i}} \} \quad \hbox{and}\quad  
{\bf L} ^{\scriptscriptstyle \geq (\Delta)} _+ = \{ i \in I \ | \  \Theta _{{\bf L} ^{\geq i}}  \subset \Delta \} \ .$$


\vskip .2cm 

\begin{lemma}\label{le:cas-produit-de-Hahn} 
Let $\Delta (\Sigma )$ be the invariance subgroup of the initial segment $\Sigma$, then we have the equality 
$${\bf L}^{(\Sigma )} = {\bf L}^{\scriptscriptstyle \geq (\Delta (\Sigma ))} \ .$$  
\end{lemma} 

\begin{preuve}   
The invariance subgroup of $\Sigma$ is defined by $\Delta (\Sigma ) = \{ \underline x \in \Theta \ | \ \underline x + \Sigma = \Sigma \}$, and for any convex subgroup $\Delta$ of $\Theta$ we have the equivalence $\Delta \subset \Delta (\Sigma ) \ \Longleftrightarrow \ \Delta + \Sigma = \Sigma$.  
Then we have the equalities 
$$\begin{array}{l} 
{\displaystyle {\bf L}_+^{\scriptscriptstyle \geq (\Delta (\Sigma ))} \ = \ \{ i \in I \ | \ \Theta _{{\bf L}^{\geq i}} + \Sigma = \Sigma \}  \phantom{\bigcup _{i\in I}}} \\ 
{\bf L}_-^{\scriptscriptstyle \geq (\Delta (\Sigma ))} \ = \ \{ i \in I \ | \ \exists \, \underline x \in \Theta _{{\bf L}^{\geq i}} \ \hbox{and} \ \exists \, \underline y  \in \Sigma \  \hbox{such that} \  \underline x + \underline y \notin \Sigma \} \ .
\end{array}$$

We deduce the inclusion ${\bf L}^{\scriptscriptstyle \geq (\Delta (\Sigma ))} _+ \subset {\bf L}^{(\Sigma )}_+$ from the inclusion $\vartheta _i (\Theta _i) \subset \Theta _{{\bf L}^{\geq i}}$. 

Conversely, we suppose that $i$ belongs to ${\bf L}^{\scriptscriptstyle \geq (\Delta (\Sigma ))} _-$ and we consider an element $\underline x = (x_l)$ in $\Theta _{{\bf L}^{\geq i}}$ and $\underline y$ in $\Sigma$ such that $\underline x + \underline y$ doesn't belong to $\Sigma$. Let $j = \iota _{\underline x}$ be the smallest element of $Supp(\underline x)$, for any $x'_j$ in $\Theta _j$ with $x'_j > x_j$ we have the inequality $\underline x < \vartheta _j(x' _j)$, hence $\vartheta _j(x'_j) + \underline y$ doesn't belong to $\Sigma$ and $j$ belongs to ${\bf L}^{(\Sigma )}_-$. 
As $\underline x$ is in $\Theta _{{\bf L}^{\geq i}}$, we have the inequality $i \leq j$, then $i$ belongs to ${\bf L}^{(\Sigma )}_-$. 

\hfill$\Box$ 
\end{preuve} 

\vskip .2cm

We recall that any convex subgroup $\Delta$ of the Hahn product $\Theta = \prod _{i \in I} ^{(H)} \Theta _i$ can be describe as in remark \ref{rmq:sous-groupe-isole-du-produit}. 
  
\begin{proposition}\label{prop:sous-groupe-de-theta} 
\emph{(a)} If the cut ${\bf L}^{(\Sigma )}$ is not of the form ${\bf L}^{> i}$, the invariance subgroup $\Delta (\Sigma )$ is the convex subgroup $\Theta _{{\bf L}^{(\Sigma )}}$ associated with this cut. 

\emph{(b)} If there exists an element $i \in I$ such that the cut ${\bf L}^{(\Sigma )}$ is the cut ${\bf L}^{> i}$, i.e. if the set ${\bf L}_-^{(\Sigma )}$ has a greatest element $i \in I$, the invariance subgroup $\Delta (\Sigma )$ corresponds to a proper convex subgroup $\Delta _i(\Sigma)$ of the group $\Theta _i$, which is defined by 
$$\Delta _i(\Sigma ) \ = \ \{ x \in \Theta _i \  |  \  \vartheta _i(x) + \Sigma = \Sigma \} \ .$$ 
\end{proposition}

\begin{preuve} 
We deduce from proposition \ref{prop:coupures-associees} and lemma \ref{le:cas-produit-de-Hahn} the inclusion $\Theta _{{\bf L}^{(\Sigma )}} \subset \Delta (\Sigma )$.  

Let $\underline x$ be an element of $\Theta$ that doesn't belong to $\Theta _{{\bf L}^{(\Sigma )}},$ with $\underline x > 0$, the smallest element $j= \iota _{\underline x}$ of $Supp(\underline x)$ satisfies $j \in {\bf L}_-^{(\Sigma )}$. 
If there exists $i \in {\bf L}_-^{(\Sigma )}$ with $i > j$, there exists $x_i$ in $\Theta _i$ and $\underline y$ in $\Sigma$ such that $\vartheta _i(x_i) + \underline y$ is not in $\Sigma$, we deduce from the inequality $\underline x > \vartheta _i(x_i)$ that $\underline x + \underline y$ is not in $\Sigma$, hence the element $\underline x$ doesn't belong to $\Delta (\Sigma )$. 

a) If the cut  ${\bf L}^{(\Sigma )}$ is not of the form ${\bf L}^{> i}$, i.e. if the initial segment ${\bf L}^{(\Sigma )} _-$ is not of the form $I _{\leq i}$, we deduce from the above the equality $\Theta _{{\bf L}^{(\Sigma )}} = \Delta (\Sigma )$. 

b) We assume now that the set ${\bf L}_-^{(\Sigma )}$ has a greatest element $i$, i.e. ${\bf L}^{(\Sigma )} = {\bf L}^{> i}$, we deduce from the above that any element $\underline x$ of $\Delta (\Sigma ) \setminus \Theta _{{\bf L}^{(\Sigma )}}$ satisfies $\iota _{\underline x} =i$, then we have the inclusions $\Theta _{{\bf L}^{> i}} \subset \Delta (\Sigma ) \subsetneq \Theta _{{\bf L}^{\geq i}}$.   

The subgroup $\Delta (\Sigma )$ is equal to ${\bf F}_i(\Delta _i(\Sigma ))$ where $\Delta _i(\Sigma )$ is a proper convex subgroup of $\Theta _i$ defined by $\Delta _i(\Sigma ) = \vartheta _i^{-1}(\Delta (\Sigma )$, hence the equality  $\Delta _i(\Sigma ) \ = \ \{ x \in \Theta _i \  |  \  \vartheta _i(x) + \Sigma = \Sigma \}$.  

\hfill$\Box$ 
\end{preuve} 

\vskip .2cm  

\begin{corollary}\label{cor:sous-groupe-non-predecesseur} 
Let $\Sigma$ be an initial segment of $\Theta$ such that the cut ${\bf L}^{(\Sigma )}$ is not of the form ${\bf L}^{> i}$, then the invariance subgroup of $\Sigma$ is not an immediate predecessor. 
\end{corollary} 

\begin{preuve} 
This is a consequence of proposition \ref{prop:sous-groupe-de-theta} and of corollary \ref{cor:non-predecesseur}. 

\hfill$\Box$ 
\end{preuve} 

\vskip .2cm

\begin{corollary}\label{cor:sous-groupe-de-theta} 
Let $\Sigma$ be an initial segment of $\Theta$ such that the cut ${\bf L}^{(\Sigma )}$ is equal to ${\bf L}^{> i}$, then an element $\underline y=(y_j)$ belongs to $\Delta (\Sigma )$ if and only if $\underline y$ belongs to $\Theta _{{\bf L}^{\geq i}}$, i.e. with $y_j = 0$ for $j<i$, and $y_i$ belongs to $\Delta _i(\Sigma )$. 
\end{corollary} 

\begin{preuve} 
From remark \ref{rmq:sous-groupe-isole-du-produit}\emph{(1)} we deduce the equality $\Delta (\Sigma ) = \Delta _i (\Sigma ) \times \prod _{j > i} ^{(H)} \Theta _j$, hence the result. 

\hfill$\Box$ 
\end{preuve} 

\vskip .2cm

Let $\Sigma$ be a non-trivial initial segment of a totally ordered group $\Theta$, with invariance subgroup $\Delta (\Sigma )$, and let $\Delta _1$ and $\Delta _2$ be two convex subgroups of $\Theta$ with $\Delta _1 \subset \Delta (\Sigma ) \subsetneq \Delta _2$. Then we can define a non-trivial initial segment $\overline{\Sigma ^{(\delta )}}$ of the group $\Delta _2 / \Delta _1$, which is unique up to linear equivalence from proposition \ref{prop:bijection}, and we can recall from lemma \ref{le:sous-groupe-invariant-quotient} and lemma \ref{le:sous-groupe-invariant-sous-groupe} how it is defined: 

- from the strict inclusion $\Delta (\Sigma ) \subsetneq \Delta _2$, we deduce that there exists $\delta \in \Theta$ such that the initial segment $\Sigma ^{(\delta )}$ of $\Delta _2$ defined by $\Sigma ^{(\delta )} \ = \ \{ \xi \in \Delta _2  \  |  \   \xi + \delta  \in \Sigma \} = (\Sigma - \delta ) \cap \Delta _2$ is non trivial, and its invariance subgroup $\Delta (\Sigma ^{(\delta )})$ is independant from the element $\delta$, and is equal to $\Delta (\Sigma )$; 

- from the inclusion $\Delta _1 \subset \Delta (\Sigma ) = \Delta (\Sigma ^{(\delta )})$, we deduce that the image $\overline{\Sigma ^{(\delta )}}$ of $\Sigma ^{(\delta )}$ by the morphism $\xymatrix{w_{\Delta _2 / \Delta _1}: \Delta _2 \ar[r] & \Delta _2 / \Delta _1}$ is a non trivial initial segment of the group $\Delta _2 / \Delta _1$, and its invariance subgroup $\Delta (\overline{\Sigma ^{(\delta )}})$ is equal to $\Delta (\Sigma ) / \Delta _1$.

\vskip .2cm

If the group $\Theta$ is the Hahn product $\Theta = \prod _{i \in I} ^{(H)} \Theta _i$ and if the initial segment $\Sigma$ of $\Theta$ is such that the cut ${\bf L}^{(\Sigma )}$ is equal to the cut ${\bf L}^{> i}$ for some $i \in I$, we deduce from the proof of proposition \ref{prop:sous-groupe-de-theta} that we have the inclusions $\Theta _{{\bf L}^{> i}} \subset \Delta (\Sigma ) \subsetneq \Theta _{{\bf L}^{\geq i}}$. 
Then we deduce from the above that we can define an initial segment $\overline{\Sigma ^{(\underline x)}}$ of the group $\Theta _i =  \Theta _{{\bf L}^{\geq i}} / \Theta _{{\bf L}^{> i}}$, for some element $\underline x$ of $\Theta$. 
By hypothesis $i$ belongs to ${\bf L}_-^{(\Sigma )}$, then there exists $x_i^0$ in $\Theta _i $ and $\underline x =(x_j)$ in $\Sigma$ such that $\vartheta _i(x_i^0) + \underline x$ doesn't belong to $\Sigma$, hence the initial segment $\Sigma ^{(\underline x )} = (\Sigma -\underline x ) \cap \Theta _{{\bf L}^{\geq i}}$ is non trivial. 

\begin{proposition}\label{prop:sous-groupe-du-quotient} 
Let $\Delta (\overline{\Sigma ^{(\underline x)}})$ be the invariance subgroup of the initial segment $\overline{\Sigma ^{(\underline x)}}$  and let $\Delta _i(\Sigma )$ be the convex subgroup of $\Theta _i$ defined by $\Delta _i(\Sigma ) = \vartheta _i^{-1}(\Delta (\Sigma )$, then we have the equality 
$$\Delta _i(\Sigma ) = \Delta (\overline{\Sigma ^{(\underline x)}}) \ .$$ 
\end{proposition}

\begin{preuve} 
It is a consequence of corollary \ref{cor:sous-groupe-de-theta} and of the equality $\Delta (\overline{\Sigma ^{(\underline x)}}) = \Delta (\Sigma ) / \Theta _{{\bf L}^{> i}}$. 

\hfill$\Box$ 
\end{preuve} 

\vskip .2cm

For any cut $\Lambda$ of the Hahn product $\Theta = \prod _{i \in I} ^{(H)} \Theta _i$, we denote by ${\bf L}^{(\Lambda )}$ the cut $\Lambda ^{\scriptscriptstyle \geq (\Delta (\Lambda ))}$ of the set $I$, by lemma \ref{le:cas-produit-de-Hahn} it also equal to the cut ${\bf L}^{(\Lambda _-)}$ associated with the initial segment $\Lambda _-$.  
If there exists $i$ in $I$ such that the cut ${\bf L}^{(\Lambda )}$ is equal to ${\bf L}^{>i}$, we denote by $\Lambda _i$ the cut of the ordered group $\Theta _i$, which is defined previously only up to linear equivalence. 

\begin{proposition} 
Let $\Lambda$ be a cut of the ordered group $\Theta = \prod _{i \in I} ^{(H)} \Theta _i$. 
\begin{enumerate}
\item If the cut ${\bf L}^{(\Lambda )}$ is not of the form ${\bf L}^{>i}$, the cut $\Lambda$ is either relatively principal or tightened. 

\item If there exists $i$ in $I$ such that we have the equality ${\bf L}^{(\Lambda )} = {\bf L}^{>i}$, the cut $\Lambda$ is of the same type as the cut $\Lambda _i$.    

\end{enumerate}
\end{proposition}

\begin{preuve} 
(1) It is a consequence of corollary \ref{cor:sous-groupe-non-predecesseur} and of theorem \ref{th:segment-initial}. 

(2) It is a consequence of proposition \ref{prop:sous-quotient}. 

\hfill$\Box$ 
\end{preuve} 

\vskip .2cm

		  

\end{document}